\documentclass [a4paper, twoside]{article}
\usepackage{latexsym}
\title{ 
Convergence in distribution  for filtering processes associated to Hidden
Markov Models with densities} 
\author{ Thomas Kaijser 
\\ 
\it\small Department of Mathematics,
\it\small Link\"{o}ping University,
\\ \it\small S-581 83 Link\"{o}ping, Sweden
\/{\rm ;} \it\small thkai@mai.liu.se }

\date{}
\begin{document}
\maketitle
\newtheorem{lem}{Lemma}[section]
\newtheorem{thm}{Theorem}[section]
\newtheorem{prop}{Proposition}[section]
\newtheorem{corr}{Corollary}[section]
\newtheorem{conjecture}{Conjecture}[section]
\newtheorem{definition}{Definition}[section]
\newtheorem{example}{Example}[section]
\newtheorem{condition}{Condition}[section]
\newtheorem{observation}{Observation}[section]

\begin{abstract}

Consider a filtering process associated to a hidden Markov model with
densities for which  both the  state space and the observation space  
are  complete, separable, metric 
spaces.
 If the underlying, hidden Markov chain is strongly ergodic and 
the filtering process fulfills a certain coupling condition 
we prove that,  in the limit, the distribution of the 
filtering process is independent of the initial distribution of 
the hidden  Markov chain. If  furthermore 
the  hidden  Markov chain is uniformly ergodic, then we prove that the 
filtering process converges in distribution.

\vspace{.5cm}

\vspace{.5cm}

\end{abstract}

\section{Introduction}

A Hidden Markov Model (HMM)  is a mathematical concept which usually is 
defined in such a way that it  consists of 
a state space, an observation space,  a transition probability
function 
(tr.pr.f) or a 
transition probability matrix (tr.pr.m) on the 
state space and a tr.pr.f or a tr.pr.m from the state space to the observation
space.

In the present  paper we shall
consider HMMs 
for which both the state space and observation space are
complete, separable, metric  spaces  with topologies and \(\sigma-algebras\) 
determined by the given metrics. To both spaces we will associate a 
\(\sigma-finite\) measure,  
which we  call \(\lambda\) and \(\tau \) respectively.
Our definition of a HMM 
(see Definition \ref{HMMdefinition1} below)  will be
slightly more general than what is usual, 
and will simply be based on a tr.pr.f
from the state space to the product space of 
the state space and the observation space.   We shall denote this 
tr.pr.f by \(M\),
and shall often assume that 
the tr.pr.f \(M\) has a probability density kernel \(m\) with respect to  
the product measure of
 the measures \(\lambda\) and \(\tau\).

A HMM generates two basic stochastic processes, 
a Markov chain, which is called the hidden Markov chain and which takes
its values in the state space, and an observation sequence  taking its values
in the observation space.
The filtering process of 
a HMM
 is, loosely speaking, 
the sequence 
 of conditional distributions of the hidden  Markov chain that is
obtained as new observations are received.

It is well-known, that the filtering process  itself, is also a 
Markov chain.
A classical, theoretical problem 
is to
find 
conditions 
 such  that 
the  filtering process, has a unique, invariant 
probability measure.

This   problem goes back to the paper  \cite{Bla57} from 1957 by  
D Blackwell for 
the case when  the hidden   Markov chain takes its values in a finite set 
and it goes back to the paper \cite{Kun71} from 1971 by H Kunita
for the case when 
the state space of the hidden Markov chain is a compact, separable, Hausdorff  space.

Blackwell studied HMMs with finite state space for which 
 the connection between the hidden Markov 
chain \(\{X_n\}\) and the observation sequence  \(\{Y_n\}\) is simply
\begin{equation}\label{nonrandom}
 Y_n=g(X_n),
\end{equation}
where thus \(g\) is a mapping from the state space to the observation space, 
and he proved that 
if the tr.pr.m of the Markov chain is ergodic and has
``rows which are nearly identical and no element  which is very small",
then  there is a
unique invariant probability measure for the filtering process.  
Blackwell also made the {\em conjecture} 
that there is unique invariant probability  measure  if  the tr.pr.m of 
the hidden Markov chain is indecomposable.

In \cite{Kun71},  Kunita considered  two coupled stochastic processes which 
one could regard as   the hidden Markov process and the observation process 
of a  {\em continuous
 time} HMM.
The hidden  Markov process \(\{X_t\}\)  was supposed to be a Feller 
process and to take its values 
in a compact, separable, Hausdorff space.
The observation process \(\{Y_t\}\) was defined by the equation 
\[
Y_t-Y_0= \int_{\tau=0}^t h(X_{\tau})dW_{\tau} + W_t-W_0,\]
where \( h\) is a continuous map from the state space to \({\tt R}^N\)
and \(\{W_t\}\) is an N-dimensional Wiener process.

In the proof of Theorem 3.3 of \cite{Kun71} - the main
theorem -, Kunita  proves
 the existence of a unique invariant probability measure, if
  the  hidden Markov process  has a unique invariant probability measure
\(\pi\) such that 
\begin{equation}\label{Kunitacondition}
\limsup_{t \rightarrow \infty}
\int_S|E[u(X_t(s))]-\langle u, \pi\rangle|\pi(ds) \; = 0, 
\forall \;real \;continuous \;\;u \;
\end{equation}
where thus \(\langle u, \pi\rangle\) means 
the integral of \(u\) with respect to \(\pi\),
\(S\) denotes the state space and \(X_t(s)\) denotes the hidden 
Markov process at time \(t\) when starting
at \(s\). (See \cite{Kun71}, formula (5).)

Kunita's proof is  based  on the observation 
that probabilities on a compact, convex set 
is partially ordered with respect 
to integration of convex functions, (see e.g \cite{Cho69}, section 26).
Kunita  considers the  two {\em extremal}  measures among the
set of probabilities on the set of probabilities on \(S\) 
which have the invariant probability measure \(\pi\) as the
barycenter.
 The smallest   
is simply \(\delta_{\pi}\) - the Dirac measure at \(\pi\). 
The other extremal measure,
the largest, is more abstract, and loosely speaking,  
it is the measure that "puts" mass \(\pi(ds)\) at the 
Dirac measure \(\delta_s\), where thus \(s\) denotes an arbitrary point in 
the given set.
By proving that the filtering process with the smallest extremal measure 
as initial measure
gives rise to a process of {\em increasing} 
probability measures with barycenter \(\pi\),  
and  the largest extremal measure  gives rise 
to a process of {\em decreasing} probability  measures with barycenter \(\pi\), 
 Kunita proves - by using (\ref{Kunitacondition}) - 
that the two limiting measures, 
both of which are invariant with respect to the Markov kernel  of the filtering
process, are equal.

Unfortunately,  approximately  30 years after its
publication,  it was 
found that there is a gap in the proof
of Theorem 3.3 of \cite{Kun71}; this gap is still not fully 
resolved but  in the paper \cite{vH09} the conclusions of 
Theorem 3.3 
are proved under slightly stronger assumptions than those made in \cite{Kun71}.
(For en extensive discussion regarding the gap in the proof of
Theorem 3.3 of 
\cite{Kun71}, 
see \cite{BCL04}.)

A  problem  closely  related to the problem of the existence of a 
unique invariant probability measure is  the following 
 convergence problem: 
{\em When does there exist a unique limit distribution towards which  the
distributions of the filtering process converge irrespectively of the
initial distribution of the hidden Markov chain?} When such a limit
distribution exists, 
then  we simply say  that 
{\em the filtering process converges in distribution} 
or that the {\em convergence property} holds. Of course, 
if the convergence property holds, then
 there also exists a unique 
invariant probability measure for the filtering process.

In the paper \cite{Kai75} from 1975 the convergence property 
was proved for a HMM with finite state
space,
under a condition called Condition A. In order to present Condition A we need 
to introduce the two notions "stepping matrix" and "subrectangular matrix".   

Consider a HMM for which both the state space and the observation space 
are finite - or denumerable -,  and 
let \(P\) be the tr.pr.m that governs the hidden Markov chain.
To every observation \(a\) one can associate a nonnegative matrix 
\(M(a)\), of the same format as \(P\), 
called the {\em stepping matrix}. 
An element \((M(a))_{i,j}\) of the stepping matrix \(M(a)\) expresses the 
probability that  the
next state of the hidden Markov chain will be the state \(j\) and the next 
observation will be  \(a\), 
given that the hidden Markov chain is in state \(i\). 

Note that 
\begin{equation}\label{partition}
\;\sum_a M(a) = P.\;
\end{equation}

A nonnegative matrix \(M\) is called  {\em subrectangular}, if 
\[
(M)_{i_1,j_1}(M)_{i_2,j_2} > 0 \Rightarrow (M)_{i_1,j_2}(M)_{i_2, j_1} > 0,
\]
where thus \((M)_{i,j}\) denotes the \((i,j)th\) element of the matrix \(M\).

In the paper 
\cite{Kai75} from 1975 
the convergence property
was proved for a HMM for which both the state space and the observation space 
are finite,  the hidden Markov chain is ergodic, the relation between
the hidden Markov chain and the observation sequence is given by 
(\ref{nonrandom}) and the following condition holds.
\newline
{\em Condition A: 
There exists a finite sequence \(\{a_1, a_2, ..., a_N\}\)   of observations such that 
the product \(\prod_{n=1}^N M(a_n)\) of stepping matrices is a nonzero, subrectangular matrix.}

The restriction to the case, when the relation between the hidden
Markov chain and the observation sequence is given by
(\ref{nonrandom}), is not a serious restriction since, 
as was first pointed out by L Baum and T Petrie  (see
\cite{BP66}), 
any   HMM  with "random observations" 
can be represented by another HMM for which (\ref{nonrandom}) holds, simply by 
1) enlarging  the original state space 
to the product space of the state space and the observation space, 
2) enlarging the tr.pr.m and 3)  
defining \(g((x,a))= a\). 
This was also pointed out in \cite{Kai75}

Also in \cite{Kai75}, a simple {\em counterexample} to Blackwell's conjecture
was given, an example which was not difficult to construct, once
condition A was found.

In the paper \cite{KR06} from 2006,  Kochman and Reeds formulated a
slightly 
weaker condition than Condition A, - a condition 
also formulated by using the stepping matrices associated to the 
elements of the observation space.

 Thus, 
consider a HMM with finite state space and finite observation space \(A\) and 
let \({\cal K}\) be the set of matrices defined by
\[
{\cal K} = 
\{cM(a_1)M(a_2)...M(a_n): n=1,2,..., \;\;a_1, a_2,...\in A, \; c \in
{\tt R},\;
c > 0\}.
\]
The condition introduced by Kochman and Reeds reads simply as follows:
\newline
\begin{equation}\label{kr}
The\;\; closure\;\; of \;\;{\cal K}\;\; contains\;\; a \;\;
 rank\;\; 1\;\; matrix.
\end{equation}
We call  the condition of Kochman and Reeds
  "the rank 1 condition"  or simply Condition KR.

In \cite{KR06}, Kochman and Reeds  proved the existence of a unique invariant 
probability measure, if the hidden Markov chain is irreducible 
 and  Condition KR holds, and, if furthermore, the hidden Markov chain is
aperiodic, they proved that the filtering process converges in distribution.

In \cite{KR06}, the authors also observed that, 
if the hidden Markov chain is irreducible and aperiodic
and Condition A is satisfied, then there exists a finite sequence 
 \(\{b_1, b_2, ..., b_N\}\)   of observations such that 
the product \(\prod_{n=1}^N M(b_n)\) of stepping matrices,  after rearrangement of the 
labelling of the states if necessary,  
can be written in the following block structure:
\begin {equation}\label{Amatrix}
\Lambda=
\left( \begin{array}{cc}
A \;0 \; B\;  0 \\
C \;  0 \;  D \; 0\\
 0\;\, 0\;\, 0\;\,   0\\
 0\;\, 0\;\, 0\;\,   0\\
\end{array}
\right).
\end{equation}
In (\ref{Amatrix})  all the elements of \(A\), \(B\), \(C\) and \(D\) are strictly
positive, 
the elements of 
the \(0- blocks\) are 0
and the formats of the blocks along the diagonal are quadratic.
By applying  Perron's theorem (see e.g. \cite{Gan65}, vol II, Theorem 8.1) to the matrix \(A\), Kochman and Reeds
prove that - after suitable  normalisation - the \(nth\) power of the matrix 
\(\Lambda\) tends to a rank 1 matrix and thereby they have showed that  
Condition A implies Condition KR.

The set of HMMs with finite state space and 
observation space, for which Condition KR holds but not Condition A,
 is probably  
quite small. 
In fact, it seems likely, that the problem of finding the set of 
HMMs which satisfy Condition KR but not Condition A, is equivalent 
to the problem of finding those HMM which 
do not satisfy Condition KR.

In the paper \cite{Kai11} published 2011,   
the convergence theorem for  HMMs with finite state space proved 
in \cite{KR06},
is generalised to  HMMs with denumerable state space. The starting
point of \cite{Kai11} is the relation (\ref{partition}) between
the tr.pr.m of the hidden Markov chain and the stepping matrices
induced by the elements of the observation space.

A difficulty  one needs to handle when analysing the filtering process  
of a HMM  with 
a {\em  denumerable and non-finite} state space is that the state space of 
the filtering process becomes   a  {\em nonlocally compact space}. To 
see this, note that 
in this case the state space of the filtering process is 
\begin{equation}\label{Kdefinition}
K= \{x=(x_1,x_2,...) \in {\tt R}^{\infty}:  x_i \geq 0, \sum_{i=1}^{\infty} x_i =1\}.
\end{equation}
If we let \({\overline B}(x_0, r_0)\) denote the closed ball 
under the \(l_1-topology\), with center \(x_0\) and radius \(r_0 >0\), 
it is easily proved and well-known that 
\({\overline B}(x_0, r_0)\) is not a compact set, from which 
follows that the set \(K\) is not  locally compact
under the topology induced by the \(l_1-norm\). This also implies
that the only real-valued continuous function on \(K\) with compact
support is the zero-function.

A nice property regarding probabilities on the set \(K,\) where 
thus \(K\) is defined by (\ref{Kdefinition}), 
is that the set of probability measures on \(K\) with {\em equal  barycenter} 
is a {\em tight} 
family of probability measures. Moreover, if one uses the 
{\em Kantorovich distance} (Vasershtein distance) to measure the distance
between the probabilities on the set \(K\), then the distance between 
the sets of probabilities with barycenter \(x\) and \(y\) respectively, is
 equal to \(||x-y||\).  
In \cite{Kai11}, these two facts, together 
with an equicontinuity property for the filtering process, 
made it possible to verify the convergence property, 
if also a certain  contraction condition, 
called Condition B, was satisfied.
 
A slight reformulation of Condition B reads as follows.
 Let \(\pi \) be a unique stationary probability vector
for the hidden Markov chain of the  HMM under consideration and
let {\bf P} denote the tr.pr.f of the filtering process.
(We call {\bf P} the filter kernel.)
\newline 
\({\bf Condition \;{\tilde B}}\):  
{\em To every \(\rho > 0\), there exists an integer
\(N\) and a constant \(\alpha > 0\), such that, if \(\{Z_{n,\mu},\;n=1,2,...\}\)
and 
\(\{Z'_{n,\nu},\;n=1,2,...\}\) are two independent Markov chains generated by the 
filter kernel {\bf P} and the initial distributions \(\mu\) and
\(\nu\) 
respectively, then
\[
Pr[\;||Z_{N,\mu} - Z'_{N,\nu}||< \rho ] \geq \alpha 
\]
if both \(\mu\) and \(\nu\) have barycenter \(\pi\).
}
We present the exact formulation of Condition B at  the end of Section 
\ref{sectionrandommapping}.
 
In \cite{Kai11}, a somewhat stronger condition called Condition B1  was
also 
introduced, a condition 
which is a more direct generalisation of the rank one condition 
of  Kochman and Reeds.  In brief, Condition B1 says essentially, that there shall exist
an infinite sequence \(a_1,a_2,...,a_n,...\) of observations such that
the normalised product of stepping matrices 
\[
\frac{\prod_1^N M(a_n)}
{||\prod_1^N M(a_n)||}\;\;
tends\;\; towards\;\;  a\;\; 
rank\;\; 1 \;\;matrix \;as \;N \rightarrow \infty.
\]

Also in \cite{Kai11}, a sufficient condition for when there
are more than one invariant probability measure was given.

In the paper 
\cite{CvH10} from 2010,  
P Chigansky and R van Handel  
prove the convergence property 
for HMMs with denumerable state space under a  contraction condition
which they 
call Condition C, a condition which they also prove is necessary.
 (For the formulation of Condition C, see \cite{CvH10} page 2325.)
In case  
the state space and the observation space are finite, 
they also verify that Condition C is equivalent
to 
Condition KR, thereby
proving that Condition KR is a {\em necessary} condition for convergence 
in distribution of filtering processes associated to a HMM with finite 
state space and finite observation space.
They also verify that both Condition B and Condition B1 of \cite{Kai11} imply
Condition C.

The work of Chigansky and van Handel in \cite{CvH10}
 has clear resemblance to the work 
of Kunita in \cite{Kun71}.
Just as in \cite{Kun71}, Chigansky and van
Handel 
considers two 
extremal invariant probability measures for the filtering process. 
By using Condition C and the  partial ordering for probabilities on the set of 
probabilities on the state space with the same barycenter induced by 
integration of convex functions, Chigansky and van Handel are 
able to prove that the two extremal invariant measures for the 
filtering process must be equal; at the same time they obtain that the 
convergence property holds.

In spite of the fact that Condition C has been proven to be both 
a necessary and sufficient condition for convergence in distribution 
of the filtering process
of a denumerable HMM for which the hidden Markov chain is strongly
ergodic,   
the theory regarding the  convergence property for HMMs with
denumerable 
state space is 
not quite complete, 
since, in some quite  concrete 
situations, it
is not clear how to verify any of Condition B, Condition B1 or 
Condition C. 

To illustrate the difficulty, consider a positively recurrent, aperiodic random
walk on the integers and suppose that our  observation system is such, 
that we only can tell whether the Markov chain is in an odd or even state.

For this example it is not clear how to verify for example 
Condition B1 introduced above.
One reason for this difficulty  is that in this case, the two 
stepping matrices that one obtains, 
will be {\em infinite} dimensional  matrices, and,
as far as we know,     
the generalisations to infinite dimensional  matrices 
of the  Perron-Frobenius theorem for finite dimensional  matrices that 
exist (see e.g. \cite{VeJ67}), do not seem to be sufficient for the
verification
of for example Condition B1.

In \cite{Kai11}, Condition A, which was originally formulated for 
a HMM with finite state space, was generalised to the case when a 
HMM has a denumerable state space.
However, in order to prove the convergence property, 
it was necessary to require,  that there exists a sequence of 
observations, such that the product of 
the corresponding stepping matrices is a subrectangular matrix  with  only
{\em finitely many} nonzero columns. Thereby, it was again possible 
to apply Perron's theorem for finite dimensional matrices in order
to verify Condition B1.

In this paper we shall  thus consider HMMs for which both the state 
space and observation space are  complete, separable, metric spaces.
An important decision we had to make was to decide which topology and which 
\(\sigma-algebra\) we should use for the  set of probabilities  on the state space of the HMM. 
For reasons described below, we decided to use 
the total variation distance as
metric for 
this set, 
and therefore 
it was natural to choose, as topology for 
this set, 
the topology determined by 
{\em the total variation distance} and as \(\sigma-algebra\) the Borel
field induced by this topology.

This choice of topology and \(\sigma-algebra\)
is in contrast to previous work on convergence in 
distribution for filtering processes associated to a 
HMM with nondenumerable state space.
As far as we know, 
in previous work  the topology on the set of probability measures
on the state space used, when proving the convergence property or proving the 
existence of a unique invariant probability measure, has {\em always}
been the weak topology and  the
\(\sigma-algebra\) has been  the Borel field induced by the weak topology.
(See e.g \cite{Kun71}, \cite{Ste89}, \cite{DS05}, \cite{vH09}, \cite{vH12}.)  
One  natural reason for this latter choice 
is that,  if   
the hidden Markov chain has an invariant
probability measure, then there also  exists at least one 
invariant probability  measure for the filtering process. 
(See e.g \cite{vH12}, Lemma A.5.) Unfortunately we have failed to prove 
a similar result when the topology is the stronger topology determined by
the total variation distance.

An important inequality, used  in \cite{Kai11}  as well as in 
\cite{Kai75}, is the inequality
\begin{equation}\label{gamma}
 \gamma({\bf T}u) \leq 3 \gamma(u)
\end{equation}
where thus \(\gamma(\cdot)\) is a generic symbol for 
the Lipschitz constant of a Lipschitz
 continuous function, and \({\bf T}\) denotes the
 transition operator associated
to the tr.pr.f of the filtering process.

When trying to  prove the inequality (\ref{gamma}) 
for  the case when the state space and the observation space of a HMM
are complete, separable, metric spaces,  it turned out, 
that it was necessary to assume  
\newline
1) that there exist a \(\sigma-finite\) measure \(\lambda\), say,
 on the state space
and a \(\sigma-finite \) measure \(\tau\), say, on the observation space, 
\newline
2) that the set of  initial distributions 
on the  state space are restricted  to the set   
of {\em absolutely continuous} probability distributions
with respect to \(\lambda\),  
\newline
3)  that 
the transition probability function \(M\)   of the HMM
has a {\em probability density kernel} with respect to the product measure
\(\lambda \otimes \tau\),  and 
\newline
4)  that the \(\sigma-algebra\)  for the set of probabilities on the
state space, is the Borel field  
generated by the metric defined by the {\em total  variation distance}.

A {\em regularity}  condition that we need  is - loosely speaking - that if two
observations are close, then the  
two conditional distributions that these two observations give rise to, 
shall also be close.
 This condition is thus a kind of {\em continuity 
condition}. 
(See Definition \ref{regulardefinition}.)  
If this condition holds and the tr.pr.f
 \(M\) has 
a density, then we call the HMM {\em regular}.

Now, if we consider a regular HMM,  let \(K\) denote
the set of probabilities on the state space 
which are absolutely continuous with respect to
the \(\sigma-finite\) measure \(\lambda\), 
and let \({\cal E}\) denote the \(\sigma-algebra\)
on \(K\) generated by the total variation distance, 
then, by using the tr.pr.f  \(M\),  we can define a tr.pr.f  on 
the measurable space \((K,{\cal E})\), a tr.pr.f which we call  
the {\em filter kernel} and usually denote 
by {\bf P}. (See Definition \ref{filterkerneldefinition}.) The filter kernel {\bf P}, together 
with  an initial distribution \(\mu\)  on \((K,{\cal E})\), 
generates a Markov chain 
on the space 
\((K, {\cal E})\)
which we call the {\em filtering process} generated by the  HMM 
and the initial distribution \(\mu\).

A complication when trying to extend the 
convergence result obtained in \cite{Kai11}, -
a complication we have not been able to overcome -, 
is due to the fact that the set of probability measures on 
\((K,{\cal E})\)
having  the same {\em barycenter} is  {\em not a tight set}.  
For this reason, in our main theorem (see
 Theorem \ref{maintheorem} below),
we partly have to be content with proving 
that the Kantorovich distance between the distributions of two 
filtering process with different initial distributions tends to zero.

The basic assumption we make about a HMM,  besides being regular,  
is that it shall be
 {\em strongly ergodic},
that is that there shall exist a unique invariant measure \(\pi\),
such that  for all starting points \(s\) in the state space
\[ \lim_{n\rightarrow \infty} ||P^n(s, \cdot) - \pi|| = 0,
\]
where thus \(P^n(s,\cdot)\) denotes the distribution of the hidden Markov chain
at time \(n\) when starting at \(s\) and \(||\cdot ||\) denotes the total variation distance.

The special assumption regarding a regular HMM  that we make in order
to 
be able to obtain limit 
results for the distributions of the filtering process  
is formulated as a {\em coupling condition}.
\newline 
{\em {\bf Condition E}: To every \(\rho >0\), there exist an integer \(N\) and a number \(\alpha\),
such that for any two probability measures \(\mu\) and \(\nu\) 
on  \((K,{\cal E})\)
with barycenter equal to the stationary measure \(\pi\),  
there exists a coupling
\({\tilde \mu}_N\), say,  of \(\mu {\bf P}^N\) and \(\nu {\bf P}^N\),
such that
\[{\tilde \mu}_N(\{(z_1,z_2) \in K\times K : ||z_1-z_2|| < \rho \}) \geq \alpha. \]
}

What  we state in our main theorem (Theorem \ref{maintheorem}) is that,
if    the  HMM  is {\em regular},
 the hidden Markov chain is {\em strongly ergodic}
 and Condition E is satisfied,  then 
the filter kernel  
is {\em weakly contracting}, that is, in the limit, the distribution 
of the filtering process is independent of the initial distribution; if 
moreover either  the hidden Markov chain
is {\em uniformly ergodic}, or the filter kernel has an invariant 
probability measure, or
there exists an element \(x_0\in K\) 
such that \(\{{\bf P}^n(x_0, \cdot), n=1,2,...\}\)
is a tight sequence, then the filter kernel is {\em weakly ergodic}, 
that is, the filtering process
converges in distribution to a unique limit measure independent 
of the initial distribution.

As pointed out above, 
in previous papers dealing with convergence in distribution or
the existence of invariant probability measures 
for a filtering process associated
to a HMM with a complete separable metric space as state space,  
the topology for the probabilities on the state space has been
the weak topology and the \(\sigma-algebra\) has been the Borel field induced
by the weak topology. 
Therefore previous results are not quite comparable to ours.
Let us though just mention, that it seems, as if in most  
papers where a correct 
proof of the convergence property has been given, 
an important assumption has been, that
the  probability density kernel \(m\), which determines   both the  HMM and the 
probability distribution for the  next observation, shall be {\em strictly positive}. Such 
an  assumption  is not necessary to make  in order to verify Condition E.

At this point we want  to mention a related problem, namely  
the problem to characterise, when the filtering process of a HMM  
has the {\em filter stability}
property. This property has to do with a computation problem regarding the 
filtering process; since one usually does not know the initial
distribution, it is
important to know, if, in the long run - with high probability 
(with probability one), the total 
variation distance between the distribution computed with the correct 
but unknown  initial distribution and
the distribution computed with the {\em guessed}  initial distribution 
tends to zero. 
With mathematical notations, if \(x\) and \(y\) are two initial distributions 
and \(a_1, a_2, a_3,...\) is a sequence of observations obtained when \(x\) 
is the initial distribution, does it hold that
\[
\lim_{n \rightarrow \infty}  ||h(x,(a_1,a_2,...,a_n))- h(y,(a_1,a_2,...,a_n))|| = 0,
\]
where thus
\(h(x,(a_1,a_2,...a_n))\)  
 denotes the ``true'' conditional distribution at time \(n\)  
of the   hidden Markov chain computed using \(x\) as the initial distribution, 
\(h(y,(a_1,a_2,...a_n))\)  
 denotes the ``guessed'' conditional distribution at time \(n\)  
of the hidden  Markov chain computed using \(y\) as  the 
initial distribution, 
 and \(||\cdot||\)  is e.g.  the \(l^1-norm\)?
 This problem has been much studied in the last two decades. (See e.g. 
 \cite{AZ97}, \cite{BCL04}, \cite{vH09}, \cite{Ata11}.) We will not 
discuss this problem further in this paper; we only want to mention
that  the inequality proved in Section    
\ref{sectioncontracting0} (see Theorem \ref{kerneltheorem0})
is similar to inequalities used in the literature, when proving
the filter stability property for filtering processes.

The plan of this paper is as follows. 
Recall, that \(K\) denotes the set of probability measures
on the state space of the HMM under consideration,  
which are absolutely continuous with respect
a given \(\sigma-finite\) measure \(\lambda\), and that \({\cal E}\) is the 
\(\sigma-algebra\) generated by the total variation distance.

In the next section, Section \ref{sectionbasics}, we introduce  some basic
definitions and  notations.
In Section \ref{sectionregularhmm} we make a precise definition of
 the concept {\em regular HMM}
and define the {\em filter kernel} of a regular HMM. We also 
 introduce the notion 
{\em compositions of HMMs} and the notion {\em iterations of a HMM} and state 
some simple facts regarding compositions of HMMs.

In Section \ref{sectionmaintheorem} we formulate the main theorem and in
  Sections \ref{sectionuniversalinequality} to 
\ref{sectionverifyingshrinkingproperty} we prove the main theorem.

 In Section \ref{sectionuniversalinequality} we prove that the filter kernel of {\em any} regular
HMM is {\em  Lipschitz equicontinuous} 
(see  Definition \ref{lipschitzdefinition1})
and in Section \ref{sectionkantorovich} we introduce  the Kantorovich
 distance for probability measures on the space \((K,{\cal E})\). 

 In Section 
\ref{sectionauxiliarytheorem}  we formulate and prove  an auxiliary theorem 
for Markov chains on a bounded, complete, separable,
 metric space.

 In Section \ref{sectionbarycenter} we prove
a simple result concerning the barycenters 
of the filtering process induced by a regular HMM and in
Section \ref{sectiondistanceproperty} we prove an inequality for  two 
different  probability measures on \((K,{\cal E})\)
with unequal  barycenters. 
Both these results are crucial to us, when  proving the main theorem.

  In section 
\ref{sectionverifyingshrinkingproperty} we conclude the proof of the main theorem by verifying
that the hypotheses of the auxiliary theorem are fulfilled.

In Section \ref{sectionrandommapping} we introduce the notion 
{\em random mapping associated to a regular HMM}. 
 The notion  random
mapping  is the same as the notion
  {\em random system with complete connections}
 (see e.g \cite{IT69}); other names for this concept is
 {\em learning model} (see e.g \cite{Nor72}) or  
{\em iterated function system with place-dependent probabilities}
(see e.g \cite{BDEG88}). 
That there is a strong connection between the theory of 
random systems with complete 
connections and the theory of HMMs (partially observed Markov chains), 
can be observed already in the paper \cite{Bla57} by Blackwell. 
(See also e.g 
\cite{IT69}, section 2.3.3.1.)

In Section \ref{sectionrandommapping} 
we also 
define the {\em Vasershtein coupling} of a random mapping
and introduce a  condition called
 Condition E1,  and by using
the Vasershtein coupling
we show  that
Condition E1 implies Condition E. 
 
At the end of Section \ref{sectionrandommapping} we
consider strongly ergodic HMMs with denumerable state space, finite or
infinite, and show, that 
the ``rank 1 condition'' introduced in the paper \cite{KR06} 
 and  Condition B introduced in  \cite{Kai11}, both
imply Condition E1.

In Section 
\ref{sectioncontracting0}, 
we prove  some inequalities for  iterations 
of positive,  integral kernels with {\em rectangular} support.
These results follow rather  easily from a theorem 
by E Hopf from 1963. (See \cite{Hop63}.) 
 In Section  \ref{sectioncontracting0}
we also
introduce yet another condition, which we call Condition P, and show that 
Condition P implies Condition E1.  
Condition P can be regarded as a generalisation of 
Condition A mentioned above.

Finally in Section \ref{sectionexamples}, 
we  present two examples.
In both examples
we start with a hidden Markov chain  
on a complete, separable, metric  state space such the tr.pr.f of 
the Markov chain
has  a probability density function \(p\) with respect to a 
\(\sigma-finite\) measure \(\lambda\).

  In the first example 
we assume, that  the state space  is  
partitioned into a denumerable
set  of subsets and  that  at each time epoch  it is only 
possible to determine in 
which subset
the hidden Markov chain is located. 
We prove that if 1)  the tr.pr.f of the
Markov chain has a probability density kernel  with respect to a \(\sigma-finite\) measure 
on the state space and 2) there exists
a subset belonging to the  partition such that  on this subset 
the probability  density kernel 
is bounded from above
and below by two positive constants, then Condition P is satisfied. 
Furthermore, by using a result in \cite{Sza06}, 
we prove that there exists  
an element \(x_0 \in K\), such that \(\{{\bf P}^n(x_0, \cdot),
n=1,2,...\}\) 
is a tight sequence, which together with the main theorem implies that 
the filter kernel is weakly ergodic.

In the other example we assume that 
the observation 
space is a complete,  separable, metric space, not necessarily denumerable, 
on which there is a \(\sigma-finite\) measure \(\tau\). 
We assume that the tr.pr.f \(M\), which determines the HMM,
has  a  probability density kernel  \(m\), which can be written as 
a product \(m = pq \) of two probability density kernels \(p\) and \(q\)
respectively, where thus \(p\) is the density kernel of the tr.pr.f of the
hidden Markov chain and \(q\) is the probability density kernel of a 
tr.pr.f \(Q\)
from the state space to the observation space.  
 
We prove that Condition P is satisfied,
if there exists a subset \(F_0\) of the state space and a 
subset \(B_0\) 
of the observation space,
such that, 
1) if the hidden Markov chain takes a value in \(F_0\), 
then the probability that the next
observation is in \(B_0\) is positive, and, 
2) if an observation  in \(B_0\) is obtained, 
then it follows that the position of the hidden Markov chain 
{\em must be} in the set \(F_0\).

We end this introductionary section with a few remarks. First, in  
 Section 11 of \cite{Kai11},  we gave an example of a  
HMM with finite state space and observation space such that the
filtering process  becomes a {\em  periodic} Markov chain,   
in spite of the fact that the hidden Markov chain is uniformly ergodic.
It is easy to generalise this example to a regular HMM for which  
the state space of the 
HMM is a finite interval, 
the observation space is finite and 
the hidden Markov chain is uniformly ergodic.  

Secondly, as pointed out above, in  \cite{CvH10} 
the authors proved, that the ``rank 1 condition'' of Kochman and Reeds,  introduced in \cite{KR06},  is also a necessary condition  for weak ergodicity 
of the filtering process associated to HMMs with finite state space
and observation space, when the hidden Markov chain is uniformly ergodic. 
We believe that similarly,  if we have a regular HMM with uniformly
ergodic hidden Markov chain, then 
Condition E is a necessary condition  for the converge property to hold. 

Thirdly, there are many other  open problems left.
 One important  problem is to generalise the
conclusions obtained in this paper to HMMs, which are not necessarily regular.
Another problem is to investigate whether, 
in the main theorem of this paper (Theorem \ref{maintheorem}), 
one can replace the conclusion  "weakly contracting" by the
conclusion  "weakly ergodic". This would follow 
if we could verify Condition \({\cal E}\) of \cite{Sza06}.

Still another problem we want to mention, is whether the technique invented
by Kunita and which was used by Chigansky and van Handel to prove the 
convergence property for HMMs with denumerable
state space, can be used also for HMMs for which the state space is a complete,
separable, metric space.

\section{Basic definitions and notations}\label{sectionbasics}

 In this section we introduce the basic concepts of 
the paper.

A hidden Markov model (HMM), as described in 
the classical paper \cite{RJ86}, consists 
of a finite state space \(S\), a finite observation space \(A\), 
a tr.pr.m \(P\) on \(S\), a
tr.pr.m \(R\) from  \(S\) to \(A\) and an initial distribution \(p_0\).
In the more modern literature, see e.g. \cite{CMR05},
 one allows both the
state space \(S\) and the observation space \(A\) to 
be measurable spaces, \((S,{\cal F})\) and \((A,{\cal A})\) say, and then, of
course, the
tr.pr.ms
\(P\) and \(R\) must be  replaced by tr.pr.fs.

 Our definition of a HMM is slightly more general than the one given in
 \cite{CMR05}, and will be based on a 
tr.pr.f
from the state space
 to the product
 of the state space
and the observation space.
 First though,  let us point out, that if  a measurable set 
\((X, {\cal X})\) and a  metric \(\phi\) on \(X\) are given,  
then we always assume implicitly, that 
there is a topology on \(X\) which is determined by the metric
\(\phi\),
 and that 
the \(\sigma-algebra \;{\cal X}\) is the Borel field induced 
 by this topology. We call such a space
 a {\em metric space} and  denote it   \((X, {\cal X}, \phi)\)
 or simply \((X,{\cal X})\).

\begin{definition}\label{HMMdefinition1}
Let
\((S,{\cal F})\) and \((A,{\cal A})\) be two measurable spaces, let
\newline
 \(M:S\times ({\cal F}\otimes {\cal A}) \rightarrow [0,1]\)
be a tr.pr.f from \((S,{\cal F})\) to \((S\times A, {\cal F}\otimes {\cal A})\)
and define the tr.pr.f  \(P:S\times {\cal F}\rightarrow [0,1]\)
by  
\(
P(s,F)= M(s, F\times A).\)
Then we call 
\begin{equation}\label{HMMequation}
{\cal H}=
\{(S,{\cal F}),P,(A,{\cal A}),M\}
\end{equation}
 a
{\em  Hidden Markov Model (HMM)}.
We call \((S,{\cal F})\) the {\em state space}, we
call \((A,{\cal A})\) the {\em observation space},
we call  
\(M\) the {\em Hidden Markov Model kernel } of \({\cal H}\) (the {\em  HMM-kernel})
and we call \(P\)  the {\em Markov kernel} of \({\cal H}\).

In case the state space is a complete, separable, metric space  
\((S,{\cal F}, \delta_0\}\), \(\lambda\)  is a positive \(\sigma-finite\)
measure on \((S, {\cal F})\), the observation space is 
a complete, separable, metric space  \( (A,{\cal A},\varrho)\), \(\tau\)
is a \(\sigma\)-finite positive measure on \((A,{\cal A})\) and 
\(m:S\times S\times A\rightarrow [0,\infty)\) is a
\({\cal F}\otimes {\cal F}\otimes {\cal A}- measurable\)  function such that
\[ 
M(s,F\times B)=\int_F\int_B m(s,t,a)\lambda(dt)\tau(da), 
\;\;\forall s\in S,\; \forall F \in{\cal F}, \;\forall B \in{\cal A},
\]
then we call \({\cal H}\) a {\em HMM with densities} and we call \(m\) the
{\em probability density kernel} of the HMM-kernel \(M\). 
We denote a HMM with
densities by
\begin{equation}\label{HMMdensity}
\{(S,{\cal F},\delta_0),(p,\lambda),(A,{\cal A}, \varrho),(m,\tau)\}
\end{equation}
where the function \(p:S\times S \rightarrow [0,\infty) \) is the function defined
by 
\[
p(s,t)=\int_Am(s,t,a)\tau(da).\]
 We call \(\lambda\) and \(\tau\) {\em base measures} and we call the tr.pr.f 
\(P:S\times {\cal F} \rightarrow [0,1]\), defined by 
\(P(s,F)=\int_Fp(s,t)\lambda(dt)\), the Markov kernel determined by \((p,\lambda)\).

If the state space \(S\) is denumerable we always assume that the associated
\(\sigma-algebra\) \({\cal F}\)  is the power set of \(S\),
 that \(\delta_0\) is the discrete metric and 
\(\lambda\)  is the  counting measure.

 Similarly,  if 
the observation  space \(A\) is denumerable, we always assume that the associated
\(\sigma-algebra\) \({\cal A}\) is the power set of \(A\), that \(\varrho\) is the discrete metric
 and that \(\tau\) is the counting measure. 
\(\; \Box\)
\end{definition}
{\bf Remark 1}.  Recall that if \((X_1, {\cal X}_1)\)
 and \((X_2, {\cal X}_2)\) are two complete, 
separable measurable spaces and \(\mu\) is a probability on
\((X_1\times X_2, {\cal X}_1\otimes{\cal X}_2)\)
 then \(\mu\) is determined by its values on 
rectangular sets \(B_1 \times B_2, \; B_1 \in {\cal X}_1, B_2 \in {\cal X}_2.\) \(\;\Box\)
\newline
{\bf Remark 2}.  Since the tr.pr.f \(P\) is determined by \(M\),
and the density kernel \(p\) is determined by \(m\),
we could have excluded \(P\) in the expression (\ref{HMMequation})
and \(p\) in the expression  
(\ref{HMMdensity}). We have included \(P\) and \(p\) for sake of
clarity.
\(\Box\)

  We shall next present  our choice of notations for some well-known
notions. Some of these notions will not be needed until later sections.

Let \((X,{\cal X},\phi)\)  be a metric space.
We let \({\cal P}(X,{\cal X})\)
denote the set of probabilities on \((X,{\cal X})\), 
we let \( {\cal Q}(X,{\cal X})\)
denote the set of finite, non-negative measures on \((X,{\cal X})\) and let 
\({\cal Q}^{\infty}(X,{\cal X})\) denote the set of 
\(\sigma-finite\), positive measures
on \((X,{\cal X})\).  If \(\mu, \nu \in {\cal Q}(X, {\cal X})\), we let 
\(\delta_{TV}(\mu,\nu) \) denote the total variation between \(\mu\) and \(\nu \)
defined by \[\delta_{TV}(\mu, \nu)=\sup\{\mu(F)- \nu(F):F \in {\cal X}\}+
\sup\{\nu(F)- \mu(F):F \in {\cal X}\}.\]
We shall also 
often
use the notation \( ||\mu -\nu ||\) instead  of \(\delta_{TV}(\mu,\nu)\).
If \(\nu \in {\cal Q}(X,{\cal X})\) we write \(||\nu||= \nu(X).\)
We always assume implicitly,
 that the topology on \({\cal Q}(X, {\cal X})\) is the toplogy
generated by the total variation metric \(\delta_{TV}\).

 We let 
\(B_u[X]\) denote the set of real, 
\({\cal X}-measurable\) functions on \(X\)  and 
let \(B[X]\) denote the set of  real, {\em bounded},  
\({\cal X}-measurable\) functions on \(X\). 
We may write \(B[X, {\cal X}]\) instead of \(B[X]\).
If \(u \in B[X]\), we set \(||u|| = \sup\{|u(x)|, x \in X\}\), we set 
\(osc(u) =\sup \{u(x)-u(y): x,y \in X\}\) and, if \(u \in B[X]\) and 
\(A \subset X\), we set \(osc_A(u) = \sup \{u(x)-u(y): x,y \in A\}\). 
If \(u \in B_u[X]\) and \(\nu \in {\cal P}(X,{\cal X})\) then, 
when convenient, we
write  \(\int_X u(x)\nu(dx) = \langle u,\nu\rangle\) if the integral exists.
If \(\lambda \in {\cal Q}^{\infty}(X,{\cal X})\) and 
 \(\nu \in {\cal Q}(X,{\cal X})\) are such that 
there exists a function \(f \in B_u[X]\)
such that 
\[\nu(F)=\int_Ff(x)\lambda(dx), \;\;\forall F \in {\cal X},\]
then we write \(\nu \in {\cal Q}_{\lambda}(X,{\cal X})\) and 
we call \(f\) a {\em representative} of \(\nu\).
If also \(\nu \in {\cal P}(X,{\cal X})\), 
we write \(\nu \in {\cal P}_{\lambda}(X,{\cal X})\).

We let \(C[X]\) denote the set of real,
bounded, continuous functions on \(X\).  
If \(u \in C[X]\), we define
\(\gamma(u)=\sup\{\frac{u(x_1)-u(x_2)}{\phi(x_1,x_2)}: x_1\not = x_2\}\), we 
define \(Lip[X]= \{u \in C[X]: \gamma(u) < \infty\}\) and we define 
\(Lip_1[X]=\{u\in Lip[X]: \gamma(u) \leq 1\}\).

 If
\(Q:X\times {\cal X} \rightarrow
[0,1]\) is a tr.pr.f on  \((X, {\cal X})\), then we define
\(Q^n:X\times {\cal X} \rightarrow
[0,1]\) recursively by \(Q^1=Q\) and 
\[
Q^{n+1}(x, F) = \int_X Q(x,dx')Q^n(x',F), \;\;
n=2,3,... \,.\]
We call the mapping  \(T:B[X] \rightarrow B[X]\)
defined by \(Tu(x)=\int_X u(y)Q(x,dy)\) the {\em transition operator} associated to
the tr.pr.f \(Q\). The tr.pr.f \(Q\) also induces a map 
\({\breve Q}: {\cal P}(X, {\cal X}) \rightarrow {\cal P}(X, {\cal X}) \)  
by \({\breve Q}(\mu)(F) = \int_XQ(x,F)\mu(dx)\). We shall usually 
write \({\breve Q}(\mu) = \mu Q\).
 As is well-known
\begin{equation}\label{PTequality} 
\langle u, \mu Q\rangle = \langle Tu, \mu \rangle .
\end{equation}
(See \cite{Rev75}, Section 1.2.)
Furthermore, if \(\mu,\nu \in {\cal P}(X, {\cal X})\) and \(u \in B[X]\), it is 
well-known that 
\begin{equation}\label{oscinequality}
|\int u(x)\mu(dx) - \int u(x)\nu(dx)| \leq osc(u)||\mu - \nu||/2,
\end{equation}
an inequality we shall have use of later.

If \((X_1,{\cal T}_1)\) and 
\((X_2,{\cal T}_2)\) are two topological  spaces, the topology
on \(X_1\times X_2\) will always be the product topology.

The  terminology below is not standard
and therefore we make a more formal definition.
\begin{definition}\label{lipschitzdefinition1}
Let \((X,{\cal X},\phi)\) be a metric space 
and
\(Q\)  a tr.pr.f on \((X,{\cal X})\).
\newline
I .  If the associated transition
operator \(T\) satisfies 
\[
 u \in Lip[X] \Rightarrow Tu \in Lip[X],\]
then we call \(Q\) {\bf Lipschitz-continuous}.
\newline
II. If \(Q\) is Lipschitz-continuous and also there exists
a constant \(C > 0\) such that the associated transition operator 
\(T\) satisfies 
\begin{equation}\label{Lipequicontinuity}
\gamma(T^nu) \leq C \gamma(u),\;\; n=1,2,..., \forall u \in Lip[X],
\end{equation}
then we call \(Q\) {\bf Lipschitz equicontinuous}. The smallest
constant C for which (\ref{Lipequicontinuity}) holds is called the 
{\bf bounding constant}.
\(\;\Box\) 
\end{definition}

We shall now introduce some terminology concerning  the limiting behaviour of 
the distributions of a Markov chain on a metric space.

\begin{definition}\label{strongergodicdefinition}
Let \((X,{\cal X},\phi)\) be a metric space and \(Q\) a tr.pr.f on 
\((X,{\cal X})\). 
\newline
1) If there exists a probability measure \(\pi \in {\cal P}(X,{\cal X})\) such  that 
\[
\lim_{n\rightarrow \infty} \delta_{TV}(Q^n(x,\cdot), \pi)=0, \forall x \in X,\]
then we call the tr.pr.f \(Q\) {\bf strongly ergodic}, and  we call 
\(\pi\) the limit measure.
\newline
2) 
 If furthermore
\[
\lim_{n\rightarrow \infty} \sup_{x \in X}\delta_{TV}(Q^n(x,\cdot), \pi)=0, \]
then we call the tr.pr.f \(Q\) {\bf uniformly ergodic}.
\end{definition}

\begin{definition}\label{strongergodicdefinitionHMM}
Let  \({\cal H}=\{(S,{\cal F}), P, (A,{\cal A}), M\}\) be 
a  HMM such that \((S,{\cal F})\) is a metric space.
 If the Markov kernel \(P\) is strongly ergodic 
(with limit measure \(\pi\)), then we also call  
\({\cal H}\) 
 strongly ergodic
 (with limit measure \(\pi\))  and, if furthermore
 the Markov kernel \(P\) is uniformly ergodic, we call 
\({\cal H}\)
uniformly ergodic.
\end{definition}
\begin{definition}\label{weakergodicdefinition}
Let   \((X,{\cal X},\phi)\) be a metric space  and let \(Q\) be  a tr.pr.f on 
 \((X,{\cal X})\). 
\newline
1)  If
\[
\lim_{n\rightarrow \infty} \sup\{\langle u, Q^n(x,\cdot) \rangle - 
\langle u, Q^n(y,\cdot)\rangle :
u\in Lip_1[K]\;\}= 0, 
\]
for all \(x,y \in X\), 
then we call the tr.pr.f  \(Q\) {\bf weakly contracting}.
\newline
2) If furthermore there exists a probability measure 
\(\pi \in {\cal P}(X,{\cal X})\),
such that 
\[
\lim_{n\rightarrow \infty} \langle u,Q^n(x,\cdot)\rangle =
\langle u, \pi \rangle, \;\forall u \in C[X], \forall x \in X,
\]
then we call the tr.pr.f \(Q\) {\bf weakly ergodic}
and we call \(\pi \) the limit measure.
\end{definition}

\section{Regular HMMs and the filter kernel}\label{sectionregularhmm}
In this section we shall introduce a more restricted  class of HMMs which we 
call {\em regular} HMMs. We shall define {\em iterations} of regular HMMs, we shall define 
the {\em filter kernel} of a regular HMM and 
shall state some simple facts regarding regular HMMs.   

We start with a HMM 
\(
{\cal H}=\{(S,{\cal F}, \delta_0),(p,\lambda),  (A,{\cal
  A},\varrho), 
(m, \tau)\}\)
with densities.  First, let us for  
each \(a \in A\) define a mapping 
\(M_a:{\cal Q}_{\lambda}(S,{\cal F})\rightarrow {\cal Q}_{\lambda}(S,{\cal F})\)
by
\begin{equation}\label{Mamap}
M_a(x)(F)=\int_{s \in S}\int_{t\in F}m(s,t,a)x(ds)\lambda(dt).
\end{equation}
We shall usually  write \(xM_a\) instead if \(M_a(x)\).

We also define a mapping 
\({\overline M}:{\cal Q}_{\lambda}(S,{\cal F})\times A 
\rightarrow {\cal Q}_{\lambda}(S,{\cal F})\) by
\begin{equation}\label{Mmap}
{\overline M}(x,a)=xM_a.
\end{equation}

In order to be able to verify that certain  sets are   measurable,  it 
has been necessary
for us  to introduce a more restricted class of HMMs with densities.

\begin{definition}\label{regulardefinition}
Let \({\cal H}=\{(S,{\cal F}, \delta_0),(p,\lambda),  (A,{\cal
  A},\varrho), 
(m, \tau)\}\) be a HMM with densities. 
If the function 
\({\overline M}:{\cal Q}_{\lambda}(S,{\cal F})\times A 
\rightarrow {\cal Q}_{\lambda}(S,{\cal F})\) 
 defined by 
 (\ref{Mmap}) and (\ref{Mamap}) is {\bf continuous} then 
we call \({\cal H}\) a {\bf regular} HMM.  \(\;\Box\)
\end{definition}

 A trivial    example of a regular HMM is a HMM with densities 
for which the observation space is denumerable, the metric \(\varrho\)
is the discrete metric and \(\tau\) is the counting measure,  since in this case
\[
||xM_a-yM_a|| \leq ||xP - yP|| \leq ||x-y||,\; \forall a\in A, \;\forall 
x,y \in Q_{\lambda}(S, {\cal F}).
\] 
For a less trivial  example see  Example
\ref{secondexample} in Section \ref{sectionexamples}.

Our next aim is to introduce a notion 
for regular  HMMs,
which we call the
{\bf filter kernel}.
Thus, let 
\(
{\cal H}=\{(S,{\cal F},\delta_0),(p,\lambda), (A,{\cal A},\varrho), (m,\tau)\}
\)
be a regular HMM.
In order to simplify the notations we shall let \(K\) be defined as
the set 
\[K ={\cal P}_{\lambda }(S, {\cal F}).\]
 Let
\(\delta_{TV}\) be the metric determined by the total variation on
 \(K\) 
and  let \({\cal E}\) be the \(\sigma-algebra\) on 
\(K\)
generated by \(\delta_{TV}\).
In agreement with our notations introduced above, 
we  let \({\cal P}(K,{\cal E})\) denote the set of probability measures 
on \((K,{\cal E})\).

We now define \(g:K\times A \rightarrow [0,\infty)\),
by
\begin{equation}\label{gdefinition00}
g(x,a)=||xM_a||,
\end{equation}
we define 
 \(G:K\times {\cal A} \rightarrow [0,1]\) by  
\begin{equation}\label{Gdefinition00}
G(x,B)=\int_B g(x,a)\tau(da), 
\end{equation}
and we define 
 \(h:K\times A \rightarrow K \) by
\begin{equation}\label{hdefinition00}
h(x,a)= xM_a/||xM_a|| \;\;if \;\; ||xM_a|| > 0
\end{equation}
\begin{equation}\label{hdefinition000}
h(x,a)=x \; \; if \; \; ||xM_a|| = 0.
\end{equation}
 
Since \({\cal H}\) is assumed to be regular,  it follows immediately that \(g\)
is continuous.
That \(G\) is  a tr.pr.f follows from the integral definition of 
\(G\) and the  fact that
\[
\int_A ||xM_a||\tau(da)= 1, \;  \forall x \in K.
\]
 That \(h\) is continuous on the set 
\(\{(x,a):||xM_a||> 0\}\)
follows as a simple consequence of the  following lemma.
\begin{lem}\label{norminequality}
Let \(x,y\) belong to a normed vector space and suppose that
\(||x||>0\) and \(||y||>0\). Then 
\[
||\frac{x}{||x||} - \frac{y}{||y||}|| \leq \frac{2||x-y||}{||x||}.
\;\Box
\] 
\end{lem}
The  inequality of Lemma \ref{norminequality}
is  easily proved by using the triangle inequality. 
We omit the details. (For details 
see e.g. \cite{Kai09}, section 3.)

Since \({\overline M}(x,a)\) is a continuous function it follows that  
the set \(\{(x,a): ||xM_a||=0\}\) is a closed set, and it is then
easily checked that
 \(\{(x,a): h(x,a) \in B\}\in {\cal K}\otimes {\cal A}\),
 if \(B\) is an open set in \({\cal E}\), 
from which follows that 
\(h:K\times A \rightarrow K\) is a measurable function. 

Next, for each \(x \in K\), we define \(A^+_x= \{a \in A: ||xM_a||>
0\}\), which is an open set for all \(x \in K\).
We now define
 the tr.pr.f \({\bf P}\) on \((K,{\cal E})\)
by 
\begin{equation}\label{filter}
{\bf P}(x,E)= \int_
{ A_x^+}
I_E(\frac{xM_a}{||xM_a||})||xM_a||\tau(da)
\end{equation}
and we define 
\({\bf T}:B[K]\rightarrow B[K] \) by
\begin{equation}\label{operator}
{\bf T}u(x)=\int_
{ A_x^+}
u(\frac{xM_a}{||xM_a||})||xM_a||\tau(da).
\end{equation}
That \({\bf P}(x, \cdot)\) is a probability measure in \({\cal
  P}(K,{\cal E})\)
for every \(x \in K\) follows from the  integral definition of \({\bf P}\).

To verify that \({\bf P}(\cdot, E)\) is measurable for each \(E \in {\cal E}\) we argue as follows.
Define \({\bf P}' : K\times {\cal E} \rightarrow [0,1]\) by 
\[
{\bf P}'(x,E) = G(x, B(x,E))
\]
where \(G\) is defined by (\ref{Gdefinition00}) and \(B(x,E)=\{a \in A: h(x,a) \in E\}\).
Clearly 
\[
G(x,B(x,E))= \int_{B(x,E)} ||xM_a||\tau(da) = \int_{A^+_x}I_E(\frac{xM_a}{||xM_a||})
 ||xM_a||\tau(da) = {\bf P}(x,E).\]
Since \(G:K\times {\cal A}\rightarrow [0,1]\) is a tr.pr.f and 
\(h:K\times A \rightarrow A\) is measurable,
it follows that  
\(
{\bf P}': K\times {\cal E} \rightarrow [0,1]\) is a tr.pr.f on \((K,{\cal E})\)
 (see e.g \cite{Kal02}, Lemma 1.41) and since \({\bf P}' = {\bf P}\) 
we can conclude that \({\bf P}:K\times {\cal E} \rightarrow [0,1]\)
is a tr.pr.f.

That \({\bf T}\) is the transition operator associated to {\bf P},
is evident from (\ref{filter}) and (\ref{operator}).
\begin{definition}\label{filterkerneldefinition}
 We call \({\bf P}:K\times {\cal E} \rightarrow [0, \infty)\), 
defined by (\ref{filter}), the {\bf filter kernel} induced by 
the regular  HMM 
\({\cal H}=
\{(S,{\cal F}, \delta_0),(p,\lambda),  (A,{\cal A},\varrho), (m, \tau)\}\). If
 \(\{Z_{n, \mu}, n=0,1,2,...\}\) 
denotes the Markov chain generated by \(\mu\in {\cal P}(K, {\cal E})\) and the 
filter kernel \({\bf  P}\), we call 
 \(\{Z_{n, \mu}, n=0,1,2,...\}\) the {\bf filtering process } induced by 
\({\cal H}\) and the initial distribution \(\mu\).
\end{definition}

We shall also need the map  \({\breve {\bf P}}:{\cal P}(K,{\cal E}) \rightarrow
{\cal P}(K,{\cal E})\)  defined by
\[
{\breve{\bf P}}(\mu)(E) = \int_K{\bf P}(x,E)\mu(dx), \;\forall E \in {\cal E}.
\]
We usually write \(\mu{\bf P}\) instead of 
\({\breve {\bf P}}(\mu)\). From (\ref{PTequality}) follows  that
\begin{equation}\label{TPintegralequality}
\langle {\bf T} u, \mu\rangle = \langle u, \mu {\bf P}\rangle.
\end{equation}

We shall next introduce a notion we call  compositions of HMMs.
Let 
\({\cal H}_1=
\{(S,{\cal F}),P_1,(A_1,{\cal A}_1),M_1\}
\)
and 
\({\cal H}_2=
\{(S,{\cal F}),P_2,(A_2,{\cal A}_2),M_2\}
\)
be two HMMs with the same state space. 
 Define
\(
A^{1,2}=A_1\times A_2,\;\; {\cal A}^{1,2}={\cal A}_1 \otimes {\cal A}_2,\) 
define
\(M^{(1,2)}:S\times {\cal F}\times {\cal A}^{1,2} \rightarrow [0,1]\)
by
\[
M^{(1,2)}(s, F\times B_1\times B_2)= \int_SM_1(s,dt,B_1)M_2(t,F,B_2),
\]
define
\(P^{(1,2)}:S\times {\cal F}\times {\cal A}^{1,2} \rightarrow [0,1]\)
by
\[
P^{(1,2)}(s, F)=
M^{(1,2)}(s, F\times A_1\times A_2)
\]
and define
\[
{\cal H}^{1,2}=\{(S,{\cal F}),P^{(1,2)},(A^{1,2},{\cal A}^{1,2}),M^{(1,2)}\}.
\]
Obviously \({\cal H}^{1,2}\) is also a HMM; we call 
 \({\cal H}^{1,2}\)  the {\em composition} of \({\cal H}_1\) 
and \({\cal H}_2\).  For simplicity we write
\[
{\cal H}^{1,2} ={\cal H}_1 * {\cal H}_2.
\]
By Fubini's theorem follows that if \({\cal H}_1\), \({\cal H}_2\) and \({\cal H}_3\) are
three HMMs with the same state space, then
\[
({\cal H}_1*{\cal H}_2)*{\cal H}_3 =  
{\cal H}_1*({\cal H}_2*{\cal H}_3). 
\]

If \({\cal H}\)  is a HMM and \({\cal H}_n={\cal H},\; n=1,2,...,N\), where \(N\geq 2\),
we set
\[
{\cal H}^N = {\cal H}_1*{\cal H}_2*...*{\cal H}_N.
\]
We  call \({\cal H}^N\) the  \(Nth\;\) 
 {\em  iterate} or the  \(Nth\)
{\em iteration} of \({\cal H}\).
Loosely  speaking, the \(Nth\)   iteration \({\cal H}^N\) of a HMM 
\({\cal H}\) is the HMM obtained from \({\cal H}\), when one collects  the observations  
in groups of \(N\)  instead of collecting them one by one.  

Next some simple facts regarding  HMMs with densities and regular HMMs.
Thus, 
let \({\cal H}_1=
\{(S,{\cal F},\delta_0),(p_1,\lambda), (A_1,{\cal A}_1,\varrho_1), (m_1,\tau_1)\}
\)
and 
\({\cal H}_2=
\newline
\{(S,{\cal F},\delta_0),(p_2,\lambda), (A_2,{\cal A}_2,\varrho_2), (m_2,\tau_2)\}\)
be two HMMs with densities and with the same state space.  We define 
 \(m^{(1,2)}:S\times S \times A_1\times A_2 \rightarrow [0,\infty)\) by
\[
m^{(1,2)}(s,t,a_1,a_2) = \int_S m_1(s,s', a_1)m_2(s',t,a_2)\lambda(ds').
\]
Again  by using  Fubini's theorem,  it follows that   
\({\cal H}_1*{\cal H}_2\) is a HMM with densities 
such that the HMM-kernel  \(M^{(1,2)}:S\times {\cal F}\otimes {\cal A}_1\otimes {\cal A}_2 \rightarrow [0,1]\)
satisfies  
\[
M^{(1,2)}(s,F\times B_1\times B_2)= \int_F\int_{B_1}\int_{B_2} m^{(1,2)}(s,t,a_1,a_2)
\lambda(dt)\tau_1(da_1)\tau_2(da_2).
\]
Furthermore, if both \({\cal H}_1\) and \({\cal H}_2\) are {\em regular},
then it is elementary  to prove that also \({\cal H}_1*{\cal H}_2\) is regular.

Next let us note that 
the following  ``scaling property'' holds:
\begin{equation}\label{scalingproperty1}
 \frac{xM^{(1,2)}_{(a_1,a_2)}}{||xM^{(1,2)}_{(a_1,a_2)}||} = 
 \frac{\frac{ xM^1_{a_1}}{||xM^1_{a_1}||}M^2_{a_2}}
 {||\frac{ xM^1_{a_1}}{||xM^1_{a_1}||}M^2_{a_2}||}, \;\;
\;\;if \; ||xM^{(1,2)}_{(a_1,a_2)}|| > 0.
\end{equation}
If we let
\({\bf P}_1\) and \({\bf P}_2\) denote the induced  filter kernels,  let 
\({\bf T}_1\) and \({\bf T}_2\) denote the associated transition operators,
let
\({\bf P}^{(1,2)}\) denote 
the filter kernel of  \({\cal H}_1 *{\cal H}_2\)  and let 
\({\bf T}^{(1,2)}\) denote the associated transition operator,
then, by using the scaling property (\ref{scalingproperty1}),
it is not difficult to prove that 
\begin{equation}\label{cocycleT}
{\bf T}_1{\bf T}_2 = {\bf T}^{(1,2)}
\end{equation}
and that
\begin{equation}\label{cocycleP} 
{\bf P}_1{\bf P}_2 = {\bf P}^{(1,2)}.
\end{equation}
Since these relations are of importance for our proof  of the main theorem (Theorem \ref{maintheorem}), 
we prove (\ref{cocycleT}) and (\ref{cocycleP}).

The equality (\ref{cocycleP})  follows from the equality in (\ref{cocycleT})  if one uses 
the identity (\ref{TPintegralequality}).
To prove (\ref{cocycleT}),  let \(u \in B[K]\)  and set \(u_2= {\bf T}_2 u\).
From (\ref{operator}) we find that
\[
u_2(x)=  \int_{ A_{2,x}^+} u(\frac{xM_{a_2}}{||xM_{a_2}||})||xM_{a_2}||\tau_2(da_2).
\]
Hence 
\[
{\bf T}_{1}{\bf T}_{2}u(x)=
\int_{ A_{1,x}^+}u_2(\frac{xM_{a_1}}{||xM_{a_1}||})||xM_{a_1}||\tau_1(da_1) =
\]
\[
\int_{ A_{1,x}^+}\int_{A_{2, x(a_1)}^+} 
u(\frac{(\frac{xM_{a_1}}{||xM_{a_1}||}M_{a_2})}{
||\frac{xM_{a_1}}{||xM_{a_1}||}M_{a_2}||})
||\frac{xM_{a_1}}{||xM_{a_1}||}M_{a_2}||\tau(da_2)||xM_{a_1}||\tau(da_1)
=\]
\[
\int_{ A_{1,x}^+}\int_{A_{2, x(a_1)}^+} u(\frac{xM_{a_1}M_{a_2}}
{||xM_{a_1}M_{a_2}||})
||xM_{a_1}M_{a_2}||\tau(da_2)\tau(da_1)
\]
where thus
\(x(a_1)\) and \(A^+_{2,x(a_1)}\) are defined by
\[
x(a_1)= xM_{a_1}/||xM_{a_1}||,\; \; a_1 \in A_{1,x}^+  \]
and
\[
A_{2,x(a_1)}^+ = \{a_2 \in A_2: ||x(a_1)M_{a_2}||>0\}
\]
respectively.

It is easily checked that the set 
\[
B(x)=\{(a_1,a_2) \in A_1 \times A_2: || xM_{a_1}M_{a_2}|| > 0\}\]
satisfies 
\[B(x)= \{(a_1,a_2) \}\in A_1 \times A_2 : a_1 \in A_{1,x} \;and \;
a_2 \in 
A_{2,x(a_1)}\}.\]
Hence 
\[
{\bf T}_{1} {\bf T}_{2}u(x)=
\int_{B(x)}u(\frac{xM_{a_1}M_{a_2}}{||M_{a_1}M_{a_2}||})||xM_{a_1}M_{a_2}||
\tau^2(da_1,da_2)=
 {\bf T}^{(1,2)}u(x)
\]
and hence (\ref{cocycleT}) holds.

By induction  follows  that if 
\({\cal H}\)  is a regular  HMM and \({\cal H}_n={\cal H},
n=1,2,...,N\), 
where \(N\geq 2\),
then \({\cal H}^N\)  is also regular, 
and if we let \({\bf P}^{(N)}\) denote the filter kernel induced by  
\({\cal H}^N\) 
and let \({\bf T}^{(N)}\) denote the transition operator 
associated to \({\bf P}^{(N)}\),
then it follows from  (\ref{cocycleP}) and (\ref{cocycleT}) that
\begin{equation}\label{PTrelations}
{\bf P}^N = {\bf P}^{(N)} \;\; and \;\; {\bf T}^N= {\bf T}^{(N)}.
\end{equation}
The second of these equalities is used in order to prove that the filter kernel 
of a regular HMM is Lipschitz equicontinuous
 and not only Lipschitz continuous, a fact
which is crucial to us, when proving the main theorem.

We end this section emphasizing 
that whenever we introduce
a HMM \({\cal H}=\{(S,{\cal F},\delta_0),(p,\lambda), (A,{\cal A},\varrho), (m,\tau)\}\)
 with densities,
then \(K\) will denote the set \({\cal P}_{\lambda}(S,{\cal F})\) and 
\({\cal E}\) will denote the \(\sigma-algebra\) on \(K\) 
generated by the total variation metric.

\section{The main theorem}\label{sectionmaintheorem}
 In this section we shall formulate the main theorem.
We shall first recall the well-known concept {\em barycenter}.

Let 
 \({\cal H}=\{(S,{\cal F}, \delta_0),(p,\lambda) ,
 (A,{\cal A},\varrho), (m,\tau)\}\) be a
HMM with densities
  and let 
\(\mu \in {\cal P}(K,{\cal E})\). The {\em barycenter} of \(\mu\),
 which we denote by
\({\overline b}(\mu)\),  is a 
probability measure  in \(K\) defined by
\[
{\overline b}(\mu)(F)= \int_K \int_Fx(ds)\mu(dx), \;\;F \in {\cal F}.
\]
That the function \({\overline b}(\mu):{\cal F} \rightarrow [0,1]\) 
is a probability in
\(K\)
is easily verified.

We let \({\cal P}(K|x)\) denote the set of probability measures in 
\({\cal P}(K,{\cal E})\) for which the barycenter is equal to \(x\).

We shall next recall the concept {\em coupling}. Let  \(\mu \) be a probability measure on 
the measurable space \((X_1, {\cal X}_1)\)
and let \(\nu\)  be a probability measure on the measurable space 
\((X_2, {\cal X}_2)\).
If \({\tilde \mu}\) is a probability measure on the product space
\((X_1\times X_2,{\cal X}_1\otimes{\cal X}_2)\) such that 
\[
{\tilde \mu}(F \times X_2)= \mu(F), \; \forall \; F \in X_1
\]
and
\[
{\tilde \mu}(X_1 \times F)= \nu(F), \; \forall F \in X_2
\]
 then we 
call \({\tilde \mu}\) a {\em coupling} of \(\mu\) and \(\nu\).

\begin{definition}\label{conditionE}
 Let \({\cal H}=\{(S,{\cal F}, \delta_0),(p,\lambda),   (A,{\cal A},\varrho), (m,\tau)\}\) be a
strongly regular HMM with limit measure \(\pi\) and let {\bf P} be the induced
filter kernel. 
We define   
{\bf Condition E} as follows :
\newline
To every \(\rho >0\), there exist an integer \(N\) and a number \(\alpha\)
such that, for any two measures \(\mu\) and \(\nu\) in \({\cal P}(K|\pi)\),
there exists a coupling
\({\tilde \mu}_N\), say,  of \(\mu {\bf P}^N\) and \(\nu {\bf P}^N\),
such that, if we set \(D_{\rho}=\{(x,y)\in K\times K: \delta_{TV}(x,y) 
<\rho\}\),
then  
\[{\tilde \mu}_N(D_{\rho}) \geq \alpha. \;\;\Box\]
\end{definition}

\begin{thm}\label{maintheorem}
Let \({\cal H}=\{(S,{\cal F},\delta_0),
(p,\lambda),(A,{\cal A},\varrho), (m,\tau)\}\)  
be a  strongly ergodic,  regular HMM with limit measure \(\pi\) and let
\({\bf P}\) be the induced filter kernel.
 Suppose also that \(\; {\cal H}\) fulfills  Condition E. 
Then \({\bf P}\)
is weakly contracting.

If furthermore, either \newline
1) there exists a measure \(\mu \in {\cal P}(K,{\cal E})\) which is
 invariant with respect to \({\bf P}\) 
or \newline
2) there exists an element 
\(x_0 \in {\cal P}_{\lambda}(S,{\cal F})\) such that 
the sequence \[\{{\bf P}^n(x_0, \cdot), n=1,2,...\}\] is a  tight
sequence, or 
\newline
3)
 \({\cal H}\)
is also  uniformly ergodic
 \newline
- then 
 \({\bf P}\)
 is  weakly ergodic. \(\;\Box\)
\end{thm}
{\bf Remark}.
In the paper \cite{Kai11} it was  proved that, 
if the state space  of a strongly ergodic, regular HMM with limit measure \(\pi\) is
denumerable, then 
\(\{{\bf P}^n(\pi, \cdot), n=1,2,...\}\)
is a tight sequence.
 We  believe the same is true, if  the state space is a complete, 
separable, metric space.
Therefore, we believe that the second part of the theorem
could be omitted and that the conclusion in the first part of the theorem 
ought to be that 
 the filter kernel \({\bf P}\)
is  "weakly ergodic"   instead of just  "weakly
contracting".
\(\;\Box\)

\section{A universal inequality}\label{sectionuniversalinequality}

\begin{lem}\label{equicontinuitylemma1}
Let
\({\cal H}=\{(S,{\cal F},\delta_0),(p,\lambda), (A,{\cal A},\varrho), (m,\tau)\}\)  
be a regular  HMM, let \({\bf P}\) be the induced  filter kernel.
Then {\bf P} is Lipschitz equicontinuous with bounding constant \(\leq 3\).
\end{lem}
{\bf Proof}. 
We shall first show that for all \(x,y \in K\) and all \(u \in Lip[K]\) 
\begin{equation}\label{Lipschitzestimate}
|{\bf T}u(x)-{\bf T}u(y)|\leq (||u||+2\gamma(u))\delta_{TV}(x,y),
\end{equation}
where thus {\bf T} is the transition operator associated to {\bf P}.

Recall, that for \(x \in K\), the set \({A^+_x}\) is defined as 
\({A^+_x}=\{a:||xM_a||>0 \}\),
where thus,  for each \(a \in A\), the mapping  
\(M_a:S\times {\cal F} \rightarrow [0,\infty)\) is 
 defined by
\(
M_a(s,F)=\int_Fm(s,t,a)\lambda(dt)
\)
 and  \(||xM_a||\) is defined by
\(
||xM_a|| = \int_S M_a(s,S)x(ds).
\)
Recall also, that \({\bf T}:B[K] \rightarrow B[K]\) is defined by
\(
{\bf T}u(x) = \int_{A^+_x} u(\frac{xM_a}{||xM_a||})||xM_a||\tau(da). 
\)

Next, let  us note, that if \(x, y \in K\) and \(a \in A\), then 
\begin{equation}\label{xyinequality}
|(||xM_a||-||yM_a||)| \leq ||xM_a- yM_a|| = ||(x-y)M_a||. 
\end{equation}
Furthermore, if \(x\) and \(y\)  in \( K\), 
and \(f\) and \( g\) in \( B_u[S]\) are representatives 
of \(x\) and \(y\) respectively, we find
that 
\[\int_A||(x-y)M_a||\tau(da) = 
\int_A\int_S|f(s)-g(s)|m(s,S,a )\lambda(ds)\tau(da) =\]
\begin{equation}\label{normequality}
\int_S|f(s)-g(s)|p(s,S)\lambda(ds) = \int_S|f(s)-g(s)|\lambda(ds) = ||x-y||.
\end{equation}

We  shall below  have use of the following proposition, the proof of which is an
immediate consequence of Lemma \ref{norminequality}.  
\begin{prop}\label{matrixestimate} Let \(x,y \in K\) and \(a\in A\) be such that 
\(||xM_a|| >0\) and \(||yM_a|| > 0\). Then
\[
\|\frac{xM_a}{\|xM_a\|} - \frac{yM_a}{\|yM_a\|} \|
\leq \frac{2\|xM_a - yM_a\|}{\|xM_a\|}. \;\;\Box
\]
\end{prop}

Now let \(x, y \in K\),
and define \(B \subset A\) by
\(
B=\{a \in A: ||xM_a||>0, ||yM_a||>0 \}.
\)
Clearly \(B\) is an open set, 
since we have assumed that \({\cal H}\) is regular.
Define \(B_1\) and \(B_2\) by
\(
B_1 = {A^+_x}\setminus B, \;\; B_2=  {A^+_y}\setminus B
\). Obviously \(B, B_1, B_2\) are disjoint, measurable sets. 
For \(u \in Lip[K]\) we now find that 
\[
|{\bf T}u(x)-{\bf T}u(y)|\leq  
\]
\[|\int_{B} (u(\frac{xM_a}{||xM_a||})||xM_a||- 
u(\frac{yM_a}{||yM_a||})||yM_a||)\tau(da)|  +
\]
\[
||u||\int_{B_1} ||xM_a||\tau(da) +||u||\int_{B_2} ||yM_a||\tau(da)\]
and by using
 Proposition \ref{matrixestimate}, (\ref{xyinequality}) 
and (\ref{normequality}),  it follows that
\[
|{\bf T}u(x)-{\bf T}u(y)|\leq\]  
\[
2\gamma(u)\int_{B}\|xM_a-yM_a\|\tau(da) + 
||u||\int_{B}||xM_a -yM_a||\tau(da) +\]
\[
||u||\int_{B_1\cup B_2}||xM_a -yM_a||\tau(da)
\leq 
 (2\gamma(u) + ||u||)||x-y|| 
\]
and thereby the inequality (\ref{Lipschitzestimate}) is proved.

From (\ref{Lipschitzestimate}) it immediately follows that 
\(
\gamma({\bf T}u) \leq 2\gamma(u) + ||u||
\)
for all \(u\) in \(Lip[K],\)
from which follows that 
\begin{equation}\label{gammainequality1}
\gamma({\bf T}u)\leq 3, \forall \,u \in Lip_1[K], 
\end{equation}
since  \(\sup\{||x-y||: x, y \in K\} =2\), and from  (\ref{gammainequality1})
then follows that 
\begin{equation}\label{gammainequality}
\gamma({\bf T}u)\leq 3\gamma(u),\; \forall u \in Lip[K].
\end{equation}

That the inequality \(\gamma({\bf T}^nu)\leq 3\gamma(u) \) also
holds for \(n\geq 2\) and all \(u \in Lip[K]\)
is an immediate consequence of equality (\ref{PTrelations})
and  the fact that 
the inequality (\ref{gammainequality}) holds for {\em all} regular HMM.
Hence the filter kernel 
\({\bf P}\)  is 
Lipschitz equicontinuous with bounding constant \(\leq 3. \)\(\;\Box\)
\newline
{\bf Remark}. It is easy to construct an example which shows that 
the   bounding constant can not be less than 2. (See \cite{Kai09}.) 
We believe the bounding constant is in fact exactly 2.

\section{ The Kantorovich distance  on the space 
\({\cal P}(K,{\cal E})\)}\label{sectionkantorovich}

Let \({\cal H}=\{(S,{\cal F},\delta_0),(p,\lambda), 
(A,{\cal A},\varrho), (m,\tau)\}\)  
be a HMM with densities and, as usual, let
 \(K={\cal P}_{\lambda}(S,{\cal F})\), and
let \({\cal E}\) be the \(\sigma-algebra\) on \(K\)
generated by \(\delta_{TV}\). 

Now, let \(K^2 = K \times K\) and \({\cal E}^2 ={\cal E} \otimes {\cal E}\).
If \(\mu\) and \(\nu\) belong to \({\cal P}(K,{\cal E})\), 
we let \({\cal P}(K^2; \mu, \nu)\) denote the subset  of 
\({\cal P}(K^2,{\cal E}^2)\) defined by
\[
{\cal P}(K^2; \mu, \nu) = \{
{\tilde \mu} \in {\cal P}(K^2,{\cal E}^2)):{\tilde \mu}(E\times K)=\mu(E),\;\;   
{\tilde \mu}(K \times E)=\nu(E),\; \forall E \in {\cal E}\}.
\]

The {\em Kantorovich distance} 
\(d_{K}(\mu,\nu)\),
for \(\mu, \nu \in {\cal P}(K, {\cal E})\), is defined as 
\begin{equation}\label{kantorovichdistance}
d_{K}(\mu,\nu) =\inf\{ \int_{K^2} \delta_{TV}(x,y) 
{\tilde \mu}(dx,dy): {\tilde \mu} \in {\cal P}(K^2;\mu,\nu)\}.
\end{equation}

Since \(0\leq \delta_{TV}(x,y) \leq 2\) for \(x,y \in K\),
it is clear that \(d_K(\mu,\nu)\) is well-defined.

From the Kantorovich-Rubenstein theorem (see \cite{Dud02},
 Theorem 11.8.2, see also \cite{Kan42}), it follows that the Kantorovich distance 
\(d_{K}\) can also be defined by
\begin{equation}\label{dualdefinition}
d_{K}(\mu,\nu) = \sup\{\int_K u(x)\mu(dx) - 
\int_K u(x)\nu(dx):
u\in Lip_1[K]\;\}.
\end{equation}
That \(d_{K}\) is a metric on \({\cal P}(K,{\cal E})\)
follows 
from  (\ref{dualdefinition}).

Since \((K,{\cal E})\) is a bounded space, it is clear that 
the metric \(d_K\) is equivalent to the metric \(\beta\) defined by
\[
\beta(\mu,\nu) =  \sup\{\int_K u(x)\mu(dx) - 
\int_K u(x)\nu(dx): \gamma(u)+||u||\leq 1\}
\]
and, as is well-known and shown in 
e.g Chapter 11 of \cite{Dud02}, the topology induced by \(\beta\)
is equivalent to the weak topology; hence the topology induced
by \(d_K\) is also equivalent to the weak topology.

Let us also note,  that if \({\bf P}\) 
denotes the filter kernel of a regular HMM, then 
it follows from  (\ref{dualdefinition}) and Definition 
\ref{weakergodicdefinition}, that
the statement "{\em \({\bf P}\) is  weakly contracting}"  is
equivalent to the statement 
\[
\lim_{n\rightarrow \infty}
d_K({\bf P}^n(x, \cdot), {\bf P}^n(y, \cdot)) = 0, 
\;\; \forall \;x,y \in K.
\]

\section{An auxiliary theorem
}\label{sectionauxiliarytheorem}
In this section 
\((K,{\cal E})\) will denote an arbitrary, bounded, complete,
separable, metric space, with metric \(\delta\) and where
 \({\cal E}\)  is the  Borel field associated to 
the topology generated  by \(\delta\).
The purpose of this section is to  state and prove  a  
limit theorem for Markov chains in
bounded, complete, separable, metric spaces.

Let  \(Q: K\times {\cal E}\rightarrow [0,1]\) 
be a tr.pr.f
on  \((K, {\cal E})\)
and  let \( T: B[K] \rightarrow B[K] \)  denote the {\em transition  operator} 
associated  to \(Q\).
We define \(T^0u(x)=u(x)\). 
Recall  that 
\begin{equation}\label{S6eq0a}
osc(T^{n+1}u) \leq 
osc(T^{n}u), \; \; \; n=0,1,2,...,  \; \; \; u\in B[K],
\end{equation}
since \(T\) is an "averaging" operator.

We shall next  define two properties  that will be part of the hypotheses
of Theorem \ref{auxiliarytheorem}.
\begin{definition}
Let \(Q\) be a tr.pr.f 
on  \((K, {\cal E})\), and let \(T\) be the
associated transition operator.
\newline
A.
If for every \(\rho > 0\), there exists a number 
\(\alpha \), \(\;0 < \alpha< 1\), and an integer \(N\)
such that,
if the integer \(n\geq N\), 
then for all \(u \in Lip[K]\)
\[
osc(T^{n}u)\leq
\alpha \rho \gamma(u) + 
(1-\alpha)
osc(T^{n-N}u),
\]
then we say that \(Q\) has the  {\bf strong shrinking property}.
We call \(\alpha \) a
{\bf shrinking number} associated to \(\rho\).
\newline
B.
If 
for every \(\rho > 0\), there exists a number 
\(\alpha \), \(\;0 < \alpha< 1\),  such that 
for every nonempty, compact set \(E \subset K\),
every \(\eta > 0\) and every \(\kappa > 0\), 
there exist an integer \(N\) and another 
nonempty, compact set \(F \subset K\) such that, 
if the integer \(n\geq N\),
then for all \(u \in Lip[K]\)
\[
osc_E(T^{n}u)\leq
\eta \gamma(u) + \kappa osc(u)+ 
\alpha \rho \gamma(u) + 
(1-\alpha)
osc_F(T^{n-N}u),
\]
then we say that \(Q\) has the  {\bf shrinking property}.
We call \(\alpha \) a
{\bf shrinking number} associated to \(\rho\).\(\,\;\Box \)
\end{definition}

\begin{thm}\label{auxiliarytheorem}
Let
\((K,{\cal E})\) be a complete,
separable,  bounded, metric space with metric \(\delta\),
let \(Q\) be a tr.p.f on \((K,{\cal E})\) and suppose that 
\(Q\) is  Lipschitz equicontinuous. 
\newline
A. Suppose also that \(Q\) has the  shrinking property. Then 
\(Q\) is  weakly contracting
(as defined in Definition \ref{weakergodicdefinition}, part 1).
\newline
B. 
Suppose furthermore that either
\newline
a)   there exists an invariant probability measure with respect to Q or
\newline
b)
there exists \(x^* \in K\) such that 
\(\{Q^n(x^*,\cdot),\; n=1,2....\}\) is a tight sequence or 
\newline
c) 
\(Q\) has the strong shrinking property,
\newline
then 
\(Q\) is  weakly ergodic (as defined in Definition \ref{weakergodicdefinition}, part 2).
\end{thm}
{\bf Proof}. 
Set 
\( D =\sup\{\delta(x,y): x,y \in K\}.
\)
Since \(K\) is assumed to be bounded we have \(D < \infty\);
it is clearly no loss of generality to assume  that \(D=2\), since the shrinking properties also hold 
if we replace the given metric \(\delta\)  by \(2\delta/D \).

In order to prove that \(Q\) is weakly contracting, we need to show,
that for all \(x,y\in K\) 
\begin{equation}\label{shrinkinglimit}
\lim_{n\rightarrow \infty} \sup 
\{|\int_{K}u(z) Q^n(x,dz) - \int_{K}u(z) Q^n(y,dz)|: u \in Lip_1[K]\} = 0.
\end{equation}

Let \(\epsilon > 0\), \(x,y \in K\) and \(u \in Lip_1[K]\) be given. 
In order to prove 
(\ref{shrinkinglimit}) we shall show, that  
 we can find an integer \(N\), which may depend  on \(x\) and
\(y\), but which does not depend on \(u\),
 such that
\begin{equation}\label{shrinkingestimate1}
\{|\int_{K}u(z) Q^n(x,dz) - \int_{K}u(z) Q^n(y,dz)| < 6\epsilon,
\;\;\;\forall \;n \geq N.
\end{equation}

This is not difficult to do, if one uses the
shrinking property. We first choose the number \(\rho\) sufficiently 
small, more precisely we set 
\(
\rho= \epsilon.
\)
Next, let \(\alpha\) be a shrinking number associated to \(\rho\).
Since \(\{x,y\}\) is a compact set, it follows from the shrinking property, 
that if we define \(\eta=\eta_1 = \epsilon/2\) and 
\(\kappa=\kappa_1=\epsilon/2\), then  we can find an integer \(N_1\)
and a compact set \(E_1\) such that, if \(n\geq N_1,\) then
\[
|\langle u,Q^n(x,\cdot)\rangle - \langle u,Q^n(y,\cdot)\rangle| =
|T^nu(x)-T^nu(y)|
\leq \eta_1+2\kappa_1 + \alpha \epsilon +
\]
\[
(1-\alpha) \sup_{z_1, z_2 \in E_1} 
|T^{n-N_1}u(z_1) - T^{n-N_1}u(z_2)|,
\] 
where we have used the fact that 
\(\gamma(u)\leq
1\), 
 \(osc(u) \leq 2\) 
and 
\(\rho = \epsilon\).

We now choose \(M=\min\{m:(1-\alpha)^m < \epsilon\}\).
For \(i =2,3,...,M,\) we  define 
the numbers \(\eta_i\)
 by
\(\eta_i=\epsilon/2^i,\)
the numbers \(\kappa_i\) by
\( \kappa_i=\epsilon/2^i\)
and having defined the compact sets \(E_i,\;\) for \(i=1,2,...,j-1,\) 
and the integers \(N_i,\) for \(i=1,2,...,j-1\),  
it follows from the shrinking property, that  we can find a compact set
\(E_j\) and an integer \(N_{j}\), such that 
\[
\sup_{z_1,z_2 \in E_{j-1}}\;|T^{n}u(z_1)-T^{n}u(z_2)|\leq
\eta_j + 2\kappa_j +\alpha \rho  +\]
\begin{equation}\label{Tuestimate} 
(1-\alpha)\sup_{z_1,z_2 \in E_j}\;|T^{n-N_j}u(z_1)-T^{n-N_j}u(z_2)|,
\end{equation}
if 
\(n \ge N_j.\;\;\)
By using (\ref{Tuestimate}) repeatedly it
 follows, that if the integer \(n\) satisfies
\(n\geq N_1+N_2+...+N_j\), then
\[
|T^{n}u(x)-T^{n}u(y)|\leq
\epsilon/2 + 2\epsilon/2 + \alpha \epsilon +
\]
\[ 
(1-\alpha) \sup_{z_1,z_2 \in E_1}\;|T^{n-N_1}u(z_1)-T^{m-N_1}u(z_2)| \leq\]
\[
\sum_{i=1}^j \epsilon/2^i + 
2\sum_{i=1}^j \epsilon/2^i + 
\epsilon \alpha(1+(1-\alpha)+(1-\alpha)^2+...
+(1-\alpha)^{j-1}) + 
\]
\[
(1-\alpha)^j\sup_{z_1,z_2\in E_j}|T^{n-(N_1+N_2+...+N_j)}u(z_1)
-T^{n-(N_1+N_2+...+N_j)}u(z_2)|.
\]
In particular, if \(j=M\) and the integer \(n\) satisfies
\(n\geq N_1+N_2+...+N_M,\) then 
\[
|T^{n}u(x)-T^{n}u(y)|\leq
\]
\[
\sum_{i=1}^M \epsilon/2^i + 
2\sum_{i=1}^M \epsilon/2^i + 
\epsilon \alpha(1+(1-\alpha)+(1-\alpha)^2+...
+(1-\alpha)^{M-1}) + 
\]
\[
(1-\alpha)^M\sup_{z_1,z_2\in E_M}|T^{n-N}u(z_1)-T^{n-N}u(z_2)|,
\]
where \(N=N_1+N_2+...+N_M\), 
and by using the  fact 
that \(osc(Tu) \leq osc(u)\), the fact that \(osc(u) \leq 2\)
 and the fact that 
\[
\epsilon \alpha(1+(1-\alpha)+(1-\alpha)^2+...(1-\alpha)^M) < \epsilon,
\]
we find that, if \(n\geq N\), then 
\[
 |T^{n}u(x)-T^{n}u(y)|<
\epsilon + 2\epsilon + \epsilon +
2(1-\alpha)^M \leq 4 \epsilon +2(1-\alpha)^M 
\]
and,  since \(M\) is defined  in such a way that 
\(
(1-\alpha)^M  < \epsilon,
\)
it follows that
\[
|T^{n}u(x)-T^{n}u(y)|=|\int_K u(z)Q^n(x,dz) -
\int_K u(z)Q^n(y,dz) |
<6\epsilon, 
\]
if 
\(n\geq N\). Hence (\ref{shrinkingestimate1}) holds 
from which follows that (\ref{shrinkinglimit})
is satisfied.
Hence \(Q\) is {\em weakly contracting}.
Thereby, the first part of Theorem \ref{auxiliarytheorem} is proved.

It remains to prove Part B of  Theorem \ref{auxiliarytheorem}.  
In order to do this 
we first prove the following lemma which is easily proved by 
using (\ref{shrinkinglimit}),
 compactness and
the Lipschitz equicontinuity property. 
We include a proof for sake of completeness.

\begin{lem}\label{uniformlimit}
 As before, let
\((K,{\cal E})\) be a complete,
separable,  bounded, metric space with metric \(\delta\),
let \(Q\) be a tr.p.f on \((K,{\cal E})\) and suppose that 
\(Q\) has the  shrinking property. Suppose also that 
\(Q\) is  Lipschitz equicontinuous. Then,
to every nonempty, compact set \(E \in {\cal E}\)
and every  \(\epsilon > 0\), we can find an integer \(N\),
such that, for any function \(u \in Lip_1[K]\),
\begin{equation}\label{shrinkingestimateagain}
\sup_{x,y \in E}|\int_{K}u(z) Q^n(x,dz) - \int_{K}u(z) Q^n(y,dz)| \leq \epsilon,
\end{equation}
for all \(n \geq N\).
\end{lem}
{\bf Proof of Lemma \ref{uniformlimit}}. 
Let \(E\in {\cal E}\) and \(\epsilon>0\) be given, where \(E\) is a 
nonempty, compact set.
Since we have assumed that \(Q\) has the Lipschitz 
equicontinuity property, there exists a constant \(C\), 
such that for all \(n \geq 1\)
\begin{equation}\label{equiproperty}
|\int_{K}u(z) Q^n(x,dz) - \int_{K}u(z) Q^n(y,dz)| \leq C\delta(x,y)\gamma(u)
\end{equation}
for all \(x,y \in K\) and all \(u \in Lip[K]\).

Next,  set \(\epsilon_1 = \epsilon/3C\). 
Since \(E\) is compact we can find  
a finite set 
\(\Psi=\{x_i, i=1,2,..., M\}\) 
consisting of \(M\) elements 
such that, for every \(x\in E\),
\(
\inf \{\delta(x,x_i): x_i      \in \Psi\} < \epsilon_1.
\) 

Further,  let \(x_i,x_j \in \Psi\) be two arbitrary elements. 
From
(\ref{shrinkinglimit})  follows that
for every pair \(x_i,\, x_j\) in \(\Psi\), we can find an integer \(N_{x_i,x_j}\)
such that, if \(n \geq N_{x_i,x_j},\) then 
\[
|\int_{K}u(z) Q^n(x_i,dz) - \int_{K}u(z) Q^n(x_j,dz)| < \epsilon/3,
\]
for all \(u \in Lip_1[K]\).
Therefore,  if we define 
\(N = \max \{N_{x_i,x_j}: (x_i,x_j) \in 
\Psi\times \Psi, x_i \not = x_j\}\),
it follows that 
\[
|\int_{K}u(z) Q^n(x_i,dz) - \int_{K}u(z) Q^n(x_j,dz)| < \epsilon/3,
\]
if \(n \geq N\),\(\;x_i,x_j \in \Psi\) and \(u\in Lip_1[K]\).

Now, let \(x,y\in E\) be chosen arbitrarily,
choose 
\(x_i\in\Psi\) such that 
\(\delta(x,x_i) < \epsilon_1\) and choose \(x_j\) such that 
\(\delta(y,x_j) < \epsilon_1\). Let
 \(u \in Lip_1[K]\).
Using the triangle inequality, (\ref{equiproperty}) and that 
\(\epsilon_1 = \epsilon/3C\), 
 we now find
that, if \(n \geq N\), then 
\[
|\int_{K}u(z) Q^n(x,dz) - \int_{K}u(z) Q^n(y,dz)| =\]
\[ 
|\int_{K}u(z) Q^n(x,dz) - \int_{K}u(z) Q^n(x_i,dz)| +
\]
\[
|\int_{K}u(z) Q^n(x_i,dz) - \int_{K}u(z) Q^n(x_j,dz)| +
\]
\[
|\int_{K}u(z) Q^n(x_j,dz) - \int_{K}u(z) Q^n(y,dz)| <\]
\[
C\delta(x,x_i) + \epsilon/3 + C\delta(x_j,y) =
C\epsilon_1 + \epsilon/3 + C\epsilon_1 = \epsilon/3 + \epsilon/3 +
\epsilon/3 = \epsilon.
 \]
Hence, 
\[
\sup \{|\int_{K}u(z) Q^n(x,dz) - \int_{K}u(z) Q^n(y,dz)|:
x,y \in E, u\in Lip_1[K]\}
\leq \epsilon, \]
if \(n\geq N\), and thereby the lemma is proved. \( \;\Box\).

We shall now complete the proof of Theorem \ref{auxiliarytheorem}
by proving Part B of the theorem.

Let us first consider the case when the tr.pr.f \(Q\) has at least one  
invariant probability  measure. 
That \(Q\) then must have precisely  one invariant probability measure
is then easily proved by a contradiction argument, if one uses 
Lemma \ref{uniformlimit} and the fact that \(Lip[K]\) is measure determining.
 We omit the details. 

Now let \(\nu\) denote the unique  invariant measure of \(Q\).
In order to prove that 
\begin{equation}\label{weaklimit}
\lim_{n \rightarrow \infty}\int_K u(y)Q^n(x,dy)=
\int_Ku(y)\nu(dy)
\end{equation}
for all \(x \in K\) and all \(u\in Lip_1[K]\),
we argue as follows. 
Let \(x \in K\) and \(u\in Lip_1[K]\) be given. Since \(\nu\) is invariant, we find
\begin{equation}\label{Testimate}
|\int_K u(y)Q^n(x,dy)-\int_Ku(y)\nu(dy)|
\leq 
\int_K |T^nu(x)-T^nu(y)|\nu(dy).
\end{equation}
Now, let \(\epsilon > 0\) be given and  choose the compact set \(C\) so
 large that \(x\in C\) and \(\nu(C) > 1- \epsilon\). This we can do
 since \((K, {\cal E})\) is a complete, separable, metric space 
and therefore every probability
measure is tight. (See \cite{Bil68}, Theorem 1.4.) From Lemma 
\ref{uniformlimit}
it follows that we can choose an integer \(N\), independent of \(u\),
so large that 
\begin{equation}\label{supestimate}
\sup_{z_1,z_2 \in C}|T^nu(z_1)- T^nu(z_2)| < \epsilon, \;\forall n \geq N.
\end{equation}
By using the inequalities (\ref{Testimate}) and (\ref{supestimate})
it now follows, that if \(n\geq N\), then
\[
|\int_K u(y)Q^n(x,dy)-\int_Ku(y)\nu(dy)|< 
\epsilon(1-\epsilon) + \epsilon osc(u) 
\leq \epsilon + 2\epsilon =3\epsilon,
\]
from which follows that
(\ref{weaklimit}) holds for all \(u \in Lip_1[K]\) 
and all \( x \in
K\). That (\ref{weaklimit}) holds for all \(u\in Lip[K]\) and all \(x\in K\),
then follows from the fact that if \( u \in Lip[K] \) and \(\gamma(u) > 0 \) 
then \(v = u/\gamma(u) \in Lip_1[K]\). 
Then, by using the same argument as used in \cite{Bil68} when proving
that (ii) of Theorem 2.1 in \cite{Bil68} implies (iii) of Theorem
2.1, it follows that 
\(\limsup_{n \rightarrow \infty} Q^n(x,F) \leq \nu(F)\) for all closed sets 
\(F \in {\cal E}\). Now, by referring to 
Theorem 2.1 of \cite{Bil68}, we find  that 
 (\ref{weaklimit}) holds for all \(u \in C[K]\) and all \(x \in K\).   
Hence \(Q\) is {\em weakly ergodic} with limit measure \(\nu\).

Next, let us assume that there exists
an element \(x^* \in K\) such that 
 \(\{Q^n(x^*,\cdot),\)
\newline 
\(n=1,2....\}\) is a tight sequence. 
To prove weak ergodicity under this assumption,
  it follows from the preceding result, that it  suffices 
to prove that
there exists an invariant 
probability measure. 
To do this we shall  use well-known arguments together with the fact that 
\(Q\) is Lipschitz-continuous.

As usual, let \(T\) denote the transition operator associated to \(Q\).
For \(n=1,2,...\) we define
 \(T^{(n)}\)
by 
\( T^{(n)} = (1/n)\sum_{k=1}^nT^k
\)
and we define \(Q^{(n)}
\) by 
\( Q^{(n)} = (1/n)\sum_{k=1}^nQ^k.
\)
Now, since 
 \(\{Q^n(x^*,\cdot), n=1,2....\}\) is a tight sequence, 
it follows immediately that also
 \(\{Q^{(n)}(x^*,\cdot), n=1,2....\}\) is a tight sequence.
Therefore  we can extract 
a subsequence \(n_j, j=1,2,...\) such that
\(\{Q^{(n_j)}(x^*,\cdot),, j=1,2,...\}\) converges weakly towards a
probability measure \(\nu\), say. 
Hence 
\begin{equation}\label{tightlimit}
\lim_{j \rightarrow \infty} T^{(n_j)}u(x^*) =\langle u, \nu \rangle 
\end{equation}
for all \(u \in C[K]\).

Now assume that \(u \in Lip[K]\).
By considering the sequence \(\{T^{(n_j+1)}u(x^*), j=1,2,...\}\)
it is easily proved that on the one hand 
\[
\lim_{j \rightarrow \infty}T^{(n_j+1)}u(x^*)=
\lim_{j \rightarrow \infty}T^{(n_j)}u(x^*)=
\langle u, \nu \rangle \]
and on the other hand 
\[
\lim_{j \rightarrow \infty}T^{(n_j+1)}u(x^*)=
\lim_{j \rightarrow \infty}T^{(n_j)}Tu(x^*)=
\langle Tu, \nu \rangle = \langle u, \nu Q\rangle,\]
where we thus have used the fact that \(Tu \in Lip[K]\)  if \(u \in Lip[K]\). 

Hence, if \(u \in Lip[K]\),
then 
\begin{equation}\label{invarianceequation}
\langle u,\nu Q \rangle  =\langle u,\nu \rangle 
\end{equation}
 holds,
and since the set of Lipschitz continuous functions is measure 
determining, it follows, that 
(\ref{invarianceequation}) holds for \(u \in C[K]\), which 
was what we wanted to prove.

To complete the proof  of Theorem \ref{auxiliarytheorem} it remains 
to prove that \(Q\) is weakly ergodic, if \(Q\) has the strong
shrinking  property.
We shall first prove that
\begin{equation}\label{oscillationlimit}
\lim_{n \rightarrow \infty} 
\sup\{ osc(T^nu): u \in Lip_1[K]\} = 0.
\end{equation}

Let \(\epsilon > 0 \) be given.   Choose \(\rho = \epsilon\). From the 
strong shrinking property follows, that we can find a 
number \(\alpha > 0\)  and an integer
\(N\), such that, if \(u \in Lip_1[K]\) and \(n> N\), then 
\begin{equation}\label{oscestimate}
osc(T^nu) \leq \epsilon \alpha + (1-\alpha)osc(T^{n-N}u).
\end{equation}
Now define 
\(
M = \min \{m: (1-\alpha)^m < \epsilon/2\}.
\)
Then, if \(n > NM\), it follows from (\ref{oscestimate}) and the fact that \(osc(u)\leq 2\)
if \(u\in Lip_1[K]\),  that
\[
osc(T^nu) 
\leq
 \epsilon \alpha + (1-\alpha)osc(T^{n-N}u) \leq\]
\[
\epsilon \alpha + (1-\alpha)
 (\epsilon \alpha + (1-\alpha)osc(T^{n-2N}u)) \leq ... <
\]
\[\epsilon \alpha (1/(1-(1-\alpha))) +2 (1-\alpha)^M 
< 2 \epsilon,
\]
if \(u \in Lip_1[K]\),  and  since  \(\epsilon \) is arbitrarily chosen,
(\ref{oscillationlimit}) follows.

Next, let \(x_0\in K\) be given. We shall now prove, that to every 
\(\epsilon>0\) we can find an integer \(N\) such that, for 
every integer  \(m\geq 1\), and every integer \(n \geq N\), 
\begin{equation}\label{cauchyestimate}
\sup \{|\int_Ku(y)Q^{n}(x_0,dy)-\int_Ku(y)Q^{n+m}(x_0,dy)|: 
u \in Lip_1[K]\} < \epsilon.
\end{equation}
Thus, let \(\epsilon > 0\)
and the integer \(m \geq 1 \)
be given. Set \(\nu_{x_0} = \delta_{x_0}Q^m\).
Then, if \(u \in Lip_1[K]\),   we find, for \(n=1,2,...\), that 
\[
|\int_Ku(y)Q^{n}(x_0,dy)-\int_Ku(y)Q^{n+m}(x_0,dy)| \leq
\int_K|T^nu(x_0)- T^nu(y)|\nu_{x_0}(dy). \]
From the limit relation (\ref{oscillationlimit}) it follows, that
we can find an integer \(N\), which is  independent of the integer \(m\), 
 such that for any  \(u \in Lip_1[K]\) and all \(y \in K\)
\( |T^nu(x_0)-T^nu(y)| < \epsilon,\) if \(n \geq N\),
which implies 
that  (\ref{cauchyestimate}) holds for all \(n\geq N\).
From the definition of the Kantorovich distance it follows that 
\[
d_K(Q^n(x_0,\cdot), Q^{m}(x_0,\cdot)) < \epsilon 
\]
if \(n,m \geq N\).
This shows that \(\{d_K(Q^n(x_0,\cdot), Q^{m}(x_0,\cdot))\}\)
is a Cauchy sequence.

Since we have assumed that \((K, {\cal E})\) is a complete, 
separable, metric space 
 it follows that
\(({\cal P}(K,{\cal E}), {\cal T}, d_K)\)  is also a 
complete, separable, metric space,
if we let \({\cal T}\) denote the Borel field generated by the 
Kantorovich metric \(d_K\).
(See e.g \cite{Dud02}, Corollary 11.5.5 and Theorem 11.8.2.)
Therefore  it follows that there exists 
a probability measure \(\mu\), say, in \({\cal P}(K,{\cal E})\),  such that
\(
\lim_{n\rightarrow \infty} 
d_K(Q^{n}(x_0,\cdot),\mu)=0.
\)
But since  
\(\lim_{n \rightarrow \infty}\sup \{ osc(T^nu): u \in Lip_1[K]\} = 0\)
because of (\ref{oscillationlimit}),
it now also follows that 
\[
\lim_{n\rightarrow \infty} 
d_K(Q^{n}(x,\cdot),\mu)=0, \;\;\forall \; x \in K,
\]
which implies that for all \(u \in Lip_1[K]\)
\begin{equation}\label{strongshrinkingresult}
\lim_{n\rightarrow \infty} \int_K 
u(y)Q^{n}(x,dy) -\int_K u(y)\mu(dy)=0, \;\;\forall \; x \in K.
\end{equation}
But if (\ref{strongshrinkingresult}) holds for
all \(u \in Lip_1[K]\), as was shown above,  it also holds for all  
\(u \in Lip[K]\).
Again referring to the proof of Theorem 2.1 in \cite{Bil68}, we can
conclude 
that \(\limsup_{n \rightarrow \infty} Q^n(x,F) \leq \nu(F)\) for all closed sets 
\(F \in {\cal E}\) and then referring to 
Theorem 2.1 of \cite{Bil68}, it follows again
that  (\ref{weaklimit}) holds for all \(u \in C[K]\) and all \(x \in K\). 
Hence \(Q\) is {\em weakly ergodic} with limit measure \(\nu\).
Thereby Theorem \ref{auxiliarytheorem} is proved.\( \;\Box.\)

\section{The barycenter of the filtering process}\label{sectionbarycenter}

From the auxiliary theorem of the previous section and Lemma 
\ref{equicontinuitylemma1}
it follows, that in order to 
prove Theorem \ref{maintheorem},  it remains to verify, that, 
if the Markov kernel 
\(P\) of the HMM under consideration is strongly ergodic, then the shrinking 
property  is satisfied, and, if \(P\)  is uniformly ergodic, 
then the {\em strong} shrinking property
is satisfied.  
In order to accomplish this we shall need two  results on barycenters
both of which are  of some independent interest. 
The first of these we shall state and prove in this section.

\begin{thm}\label{kunita}
Let \({\cal H}=\{(S,{\cal F},\delta_0),(p,\lambda),(A,{\cal A},\varrho), 
(m,\tau)\}\)  
be a regular  HMM 
 and let \({\bf P}\) be the filter kernel.
Let  
\(P\) be the 
Markov kernel
of \({\cal H}\).
Then for all \(x \in K\)  
\[
{\overline b}({\bf P}^n(x,\cdot)) = xP^n, \; n=1,2, ...\, .
\]
\end{thm}
{\bf Remark}. 
The  theorem is essentially due to Kunita. (See \cite{Kun71}. 
See also \cite{CvH10}, Lemma A.5.)\(\;\Box\)
\newline
{\bf Proof}.  
 Let \(F \in {\cal F}\) and 
\(I_F:S \rightarrow \{0,1\} \) denote the indicator function 
of \(F\). 
From the definition of the barycenter we find 
\[
{\overline b}(\delta_x{\bf P})(F)=
 \int_{A_x^+} \langle I_F,\frac{xM_a}{||xM_a||}\rangle||xM_a|| \tau(da) = \]
\[
 \int_{A_x^+} \int_F xM_a(dt) \tau(da) =
 \int_{A_x^+} \int_F\int_S m(s,t,a)x(ds)\lambda(dt) \tau(da) =\]
\[
\int_F\int_S  p(s,t)x(ds)\lambda(dt) =
 \int_F( xP)(dt) = xP(F) 
\]
from which follows that 
\(
{\overline b}(\delta_x{\bf P})=xP\).
That \(
{\overline b}(\delta_x{\bf P}^n)=xP^n\),  for \(n\geq 2\), then  follows from
the relation (\ref{PTrelations}).  
\( \;\Box\)

The following lemma is not needed in the proof of the main theorem,
but will be needed later, when we want to verify that Condition E holds.
 We present it here, since it  gives some insight into the
sets of probability measures on \((K,{\cal E})\)  with equal barycenter.

\begin{lem}\label{barycenterlemma}
Let \((S,{\cal F}, \delta_0)\) be a complete, separable metric space, let \(\lambda\) be
a \(\sigma-finite\) measure on \((S,{\cal F})\), let \(K = {\cal P}_{\lambda}(S, {\cal F})\),
 let \({\cal E}\) denote the \(\sigma-algebra\) generated by the total variation
metric and let \(\pi \in K\). 
 For \(F \in {\cal F}\) define
\(
E(F)=\{x \in K : x(F)  \geq \pi(F)/2\}.
\)
Then, for all \(\mu \in {\cal P}(K|\pi)\) and all \(F \in {\cal F}\),
\begin{equation}\label{usefulestimate}
 \mu(E(F)) \geq \pi(F)/2. 
\end{equation}
\end{lem}
{\bf Proof}. The inequality (\ref{usefulestimate}) holds trivially if 
\(\pi(F)=0\). Thus assume \(F \in {\cal F}\) is such  that 
\(\pi(F)>0\). Clearly \(E(F) \in {\cal E}\). Set \(E(F)=E.\) 
Since \(\mu \in {\cal P}(K|\pi)\) we have
\(\int_K \langle I_F, x\rangle \mu(dx) = \pi(F).\)
Hence
\[
 \pi(F)=
\int_{E}  \langle I_F, x\rangle \mu(dx)+
\int_{K\setminus E}  \langle I_F, x\rangle \mu(dx) =
\int_E \int_Fx(ds)\mu(dx) +\]
\[ \int_{K\setminus E} \int_Fx(ds)\mu(dx) \leq
\ \mu(E) + (1-\mu(E))\pi(F)/2.
\]
Hence
\(\mu(E)( 1- \pi(F)/2) \geq \pi(F)/2
\)
and hence 
\(
\mu(E(F)) > \pi(F)/2
\)
which is  more than  we needed to prove.
\(\Box\)

\section{On the Kantorovich distance between sets 
with different barycenters}\label{sectiondistanceproperty}

Let \((S,{\cal F})\)
 be a complete, separable, measurable space 
with metric \(\delta_0\), let \(\lambda\) denote a \(\sigma\)-finite, 
nonnegative measure on 
\((S,{\cal F})\) and set \(K= {\cal P}_{\lambda}(S,{\cal F})\). 
As before, let 
\(\delta_{TV}\)
denote the metric on \(K\) induced by the total variation and let \({\cal E}\)
denote the \(\sigma-algebra\) generated by \(\delta_{TV}\).
Instead of writing \(\delta_{TV}(x,y) \)  we shall in this section  usually 
write \(||x-y||\).
 Let 
\({\cal P}(K,{\cal E})\) denote the set of probability measures
on \((K, {\cal E})\),  let
\({\cal Q}(K,{\cal E})\) denote the set of positive and finite 
measures on \((K,{\cal E})\) and for \(r > 0\) let \({\cal Q}^r(K,{\cal E})\)  
denote the set of positive, finite measures
on \((K,{\cal E})\) with total mass equal to \(r\).

Let \(d_K:{\cal P}(K,{\cal E})\times{\cal P}(K,{\cal E}) \rightarrow [0,2]\) denote the
Kantorovich distance on \({\cal P}(K,{\cal E})\) (see Section \ref{sectionkantorovich}). 
 Recall that
the Kantorovich distance on \({\cal P}(K,{\cal E})\) has two
equivalent definitions namely either by the formula (\ref{kantorovichdistance})
or by the formula (\ref{dualdefinition}).

For the set \({\cal Q}^r(K,{\cal E})\) we also define a metric, 
which we also denote by \(d_K\),
simply by
\[
d_K(\mu,\nu) = rd_K(\mu/r, \nu/r), \; \mu, \nu \; \in {\cal Q}^r(K,{\cal E}).
\]
Also in this case we call \(d_K\) the Kantorovich distance.

As  in Section  \ref{sectionmaintheorem}, we let  
\({\cal P}(K|x)\) denote the set of probability measures
on \((K, {\cal E})\) for which the barycenter is equal to \(x\).
For  \(\mu \in {\cal Q}^r(K,{\cal
  E})\) 
we also define
a barycenter \({\overline b}(\mu)\) simply by 
\[
{\overline b}(\mu) = r {\overline b}(\mu/r).
\]
Thus, if \(\mu \in  {\cal Q}^r(K,{\cal E})\) then 
\({\overline b}(\mu) \in {\cal Q}_{\lambda}(S,{\cal F})\) and 
\(|| {\overline b}(\mu)|| = r\).
For \(x \in K\) and \(r > 0\), we let \({\cal Q}^r(K|x)\) denote the set of measures in  
\( {\cal Q}^r(K,{\cal E})\) which have barycenter equal to \(rx\).

The  purpose of this section is to prove the following result:
\begin{thm}\label{barycentertheorem}
Let \(r>0\), let \(x,y \in K\)
 and  let \(\mu \in {\cal Q}^r(K|x)\). Then
\[
\inf \{d_K(\mu,\nu):\nu \in {\cal Q}^r(K|y)\} = r||x-y||.
\]
\end{thm}
{\bf Proof}.
Let us first note 
that  if \(x,y \in K\),  then
\(
d_K(\delta_x,\delta_y) = 
||x-y||,
\)
where thus \(\delta_x\) and \(\delta_y\) denote the Dirac measures
 at \(x\) and \(y\) respectively. This follows  from 
(\ref{kantorovichdistance}).

The following lemma gives a  lower bound for the Kantorovich distance 
between two measures in \({\cal Q}^r(K,{\cal E})\) in terms of their
 barycenters.
\begin{lem}\label{lowerbound}  
Let \( r > 0\) and  let \(\mu, \nu  \in {\cal Q}^r(K,{\cal E})\).
Then 
\(
d_K(\mu,\nu) \geq
|| \overline{b}(\mu)- \overline{b}(\nu)||. 
\)
\end{lem}
{\bf Proof}.  
The conclusion of the lemma is trivially true if 
\(\overline{b}(\mu) = \overline{b}(\nu).\)
We thus assume that \(\overline{b}(\mu) \not = \overline{b}(\nu)\).
From the definition of the Kantorovich distance in 
\({\cal Q}^r(K,{\cal E})\) and the definition of
the barycenter of a measure in \({\cal Q}^r(K,{\cal E})\), 
it follows that it suffices to prove the inequality if \(r=1\), that is when
 \(\mu, \nu \in {\cal P}(K,{\cal E})\). 

Thus, let \(\mu,\nu \in {\cal P}(K,{\cal E})\) and set 
\(x= \overline{b}(\mu)\) and \(y=\overline{b}(\nu).
\)
Let \(F_1, F_2 \in {\cal F}\) be such that \(F_2 =S\setminus F_1\) and such that 
\(x(F\cap F_1) \geq y(F\cap F_1),\;\; \forall F \in {\cal F}\) such that \(F\subset F_1\)
and \(x( F \cap F_2) < y(F\cap F_2),\;\; \forall F \in {\cal F}\) such that 
\(F\subset F_2\).
Define the function \(J:S\rightarrow [-1,1]\)
by
\begin{equation}\label{indicatorrelation}
 J(s)=I_{F_1}(s)- I_{F_2}(s) ,
\end{equation}
where thus \(I_{F_1}\) and \(I_{F_2}\)  denote the indicator functions of the sets \(F_1\)
and \(F_2\).
 
 Next, define \(v \in 
B[K]\)
 by  
\(
v(z)=\langle J, z \rangle.
\)
Since \(osc(J)\leq 2\), it follows from (\ref{oscinequality}),
that 
\[
|v(z_1)-v(z_2)| =
|\langle J, z_1 \rangle-\langle J, z_2 \rangle|
 \leq osc(J)||z_1-z_2 ||/2 \leq ||z_1-z_2 ||\]
and hence \(v \in
 Lip_1[K]\).
 From the definition of
the Kantorovich distance it then follows that
\begin{equation}\label{kantestimatebelow} 
d_K(\mu,\nu) \geq |\int_K v(z)\mu(dz)- \int_K v(z) \nu(dz)|
\end{equation}
 and from the definition of the 
barycenter and (\ref{indicatorrelation}),
it  follows that  
\[
|\int_K v(z)\mu(dz)- \int_K v(z) \nu(dz)|=
|\int_K \langle J, z \rangle\mu(dz)- 
\int_K \langle J, z \rangle \nu(dz)| =\]
\[
| \langle I_{F_1}, \overline{b}(\mu) \rangle - 
\langle I_{S \setminus F_1}, \overline{b}(\mu) \rangle - 
\langle I_{F_1}, \overline{b}(\nu) \rangle
+\langle I_{S\setminus F_1}, \overline{b}(\nu) \rangle  
| =
\]
\[
|x(F_1)-y(F_1)+y(S\setminus F_1) -x(S\setminus F_1)| =||x - y|| =
||\overline{b}(\mu) - \overline{b}(\nu)||,\]
which together with (\ref{kantestimatebelow})
implies that 
\(
d_K(\mu,\nu)\geq  
||\overline{b}(\mu) - \overline{b}(\nu)||.\)
\(\;\Box \)

We now continue our proof of Theorem \ref{barycentertheorem}
by proving that,
if the measure \(\mu \in {\cal Q}(K,{\cal E})\) 
is a weighted finite sum of Dirac measures,
then for every \(y \in K\) we can find a measure \(\nu \in 
{\cal Q}(K,{\cal E})\), such that \(\mu(K)=\nu(K)\),
\({\overline b}(\nu)=y\mu(K)\) and 
\(
d_K(\mu,\nu)= ||{\overline b}(\mu)-{\overline b}(\nu)||.
\)
As usual, if \(\xi\) denotes an arbitrary element in \(K\), 
we let \(\delta_{\xi}\) denote the 
Dirac measure at \(\xi\).
 
\begin{lem}\label{finitecase} 
Let \(N\) be a positive integer 
and let \(\xi_k, \;k=1,2,...,N,\) be elements  in \(K\).
Let \(\beta_k>0, \; k=1,2,...,N\), let the measure 
\(\varphi \in {\cal Q}(K,{\cal E})\) be defined by 
\(
\varphi = \sum_{k=1}^N\beta_k\delta_{\xi_k}
\)
and define the element \(a\in 
{\cal Q}_{\lambda}(S,{\cal F})\) 
by
\( a = \sum_{k=1}^N\beta_k\xi_k.
\)
Let \(b \in 
{\cal Q}_{\lambda}(S,{\cal F})\) 
be an element 
satisfying 
\(
||b|| =||a||.
\)

Then, there exist elements \(\zeta_k,\;\; k=1,2,...,N,\;\) in \(\;K\),
such that 
\(
b= \sum_{k=1}^N\beta_k\zeta_k,\)
and such that, if we define 
\(
\Psi =  \sum_{k=1}^N \beta_k \delta_{\zeta_k},
\)
then
\[
d_{K}(\varphi, \Psi) = ||a-b||.
\]
\end{lem}
{\bf Proof}.
First let us observe  that, if \(\psi \in  
{\cal Q}(K,{\cal E})\) is defined by 
\(
\psi =\sum_{k=1}^N \beta_k \delta_{\zeta_k},
\)
where \(\beta_k\), for \(\; k=1,2,...,N,\) is a positive number,
and  \(\zeta_k\), for \(\;\; k=1,2,...,N,\;\) belongs to \(\;K\), then 
\begin{equation}\label{discretebarycenter}
{\overline b}(\psi) = \sum \beta_k \zeta_k .
\end{equation}
This follows  from the fact that, if 
\(\mu \in {\cal Q}(K,{\cal E})\) is defined by \(\mu = \delta_{z_0}\)
and 
\(F \in {\cal F}\),  then 
\(
\int_K \langle I_F,z \rangle \mu(dz) = 
\langle I_F,z_0 \rangle = z_0(F).
\)

Next, let \(\zeta_1, \zeta_2,...,\zeta_N\) denote an arbitrary set of \(N\) elements in \(K\) 
and define 
\(\theta \in {\cal Q}(K,{\cal E})\)  by
\(
\theta = \sum_{k=1}^N \beta_k \delta_{\zeta_k}.
\)
Clearly \(\theta(K)=\sum_{k=1}^N \beta_k\) and hence \(\theta(K) = \varphi(K)=||a||.\)
We now define the measure \({\tilde \varphi}\) on \((K^2,{\cal E}^2)\) by 
\(
{\tilde \varphi}(\{(\xi_k,\zeta_k)\})=\beta_k,\; k=1,2,...,N .
\)  
Then clearly 
\(
{\tilde \varphi}(A\times K) = \varphi(A),\; \forall A \in {\cal E}, \;\;\)
and
\(\;\;
{\tilde \varphi}(K\times A) = \theta(A), \;\forall A \in {\cal E},
\)
from which follows that the Kantorovich distance \(d_K(\varphi, \theta)\)
satisfies 
\begin{equation}\label{Kantinequality}
 d_K(\varphi, \theta) \leq \sum_{k=1}^N \beta_k ||\xi_k - \zeta_k||,
\end{equation}
since
\[
d_K(\varphi, \theta) \leq
\int_{K \times  K}||x-y||{\tilde \varphi}(dx,dy) =
\sum_{k=1}^N \beta_k ||\xi_k - \zeta_k|| .
\]
By combining
 (\ref{Kantinequality}) and 
(\ref{discretebarycenter})  with  Lemma \ref{lowerbound},  
it follows, that in order to prove 
Lemma \ref{finitecase},  it suffices to find probability measures  
 \(\zeta_k,\;\; k=1,2,...,N,\) belonging to \( K\), such that  
\begin{equation}\label{beta}
b= \sum_{k=1}^N\beta_k\zeta_k 
\end{equation}
and also
\begin{equation}\label{xi}
\sum_{k=1}^{N} \beta_k ||\xi_k - \zeta_k|| = ||a-b||.
\end{equation}

That we can do this when  \( N=1\), 
that is, when \(\varphi = \beta_1 \delta_{\xi_1}\),
is trivial.  Simply define \(\zeta_1 = b/\beta_1\); then \(\beta_1||\xi_1-\zeta_1|| =
||a-b||\), as we want it  to be.
The case when \(b=a\) is also trivial. Just take
\(\zeta_k = \xi_k,\;\; k=1,2,..., N\).
In the remaining part of the proof we therefore assume that \(a\neq b\).

We shall now prove - by induction -,  that we can find 
probability measures 
\(\zeta_k \in K,\;\; k=1,2,...,N,\) such that (\ref{beta})
and (\ref{xi}) hold.
Thus, let us  assume, that, if 
\(N=M-1\), where \(M\geq 2\),  if 
\(a= \sum_{k=1}^N \beta_k \xi_k\) where 
\(\beta_k > 0,
\;k=1,2,...N,\)
and  \(\xi_k \in K, k=1,2,...,N\), if  \(b \in {\cal Q}_{\lambda}(S,{\cal F})\) 
and also \(||a|| = ||b||\), then  we can find 
\(\zeta_k, k=1,2,...,N\) in \(K\), such that (\ref{beta})
and (\ref{xi}) hold.

Now, let \(N=M\), let   
\(\beta_k > 0,
\;k=1,2,...,M,\) let 
\(\xi_k \in K, k=1,2,...,M,\)  set 
 \(a= \sum_{k=1}^M \beta_k \xi_k\) and suppose 
that \(b \in {\cal Q}_{\lambda}(S,{\cal F})\) and that \(||b||= ||a||.\)
 Our aim is thus to 
find elements  \(\zeta_k, k=1,2,...,M\) in \(K\), such that 
\begin{equation}\label{beta2}
b= \sum_{k=1}^M\beta_k\zeta_k 
\end{equation}
and also
\begin{equation}\label{xi2}
\sum_{k=1}^{M} \beta_k ||\xi_k - \zeta_k|| = ||a-b||.
\end{equation}

Recall that we have assumed that \(a \not = b\) and hence 
\(||a-b||\not = 0.\) We define 
\[
\Delta= ||a-b||/2.
\]
Let us also define \(a_1 \in {\cal Q}_{\lambda}(S, {\cal F})\) 
by 
\(
a_1 = \sum_{k=1}^{M-1} \beta_k \xi_k.\)
Clearly \(||a_1||=||a||- \beta_M.\) 

Now suppose that we can find a probability measure \(\zeta_M \in K\),
such  that, if we define 
\begin{equation}\label{b1definition}
b_1 = b - \beta_M\zeta_M,
\end{equation}
 then
\begin{equation}\label{conditionb11}
b_1\in {\cal Q}_{\lambda}(S,{\cal F})
\end{equation}
and 
\begin{equation}\label{conditionb12}
||a-b||=||a_1-b_1||+ \beta_M ||\xi_M- \zeta_M||.
\end{equation}
From (\ref{conditionb11}) and the definition of \(b_1\) it then follows 
that \(||b_1||=||b||-\beta_M = ||a||-\beta_M = ||a_1||\) and then,
using the induction hypothesis, it follows that we can find
probability measures \(\zeta_k, k=1,2,...,M-1,\) such that
\(b_1= \sum_{k=1}^{M-1}\beta_k\zeta_k\)
and 
\begin{equation}\label{diffequality}
\sum_{k=1}^{M-1}\beta_k||\xi_k - \zeta_k|| = ||a_1-b_1||,
\end{equation}
and consequently, by using (\ref{conditionb12}) and
(\ref{diffequality}),
it follows that  
\[||a-b||= 
\sum_{k=1}^{M-1}\beta_k||\xi_k - \zeta_k|| + \beta_M||\xi_M-\zeta_M|| =
\sum_{k=1}^{M}\beta_k||\xi_k - \zeta_k|| 
\]
and hence  (\ref{beta2})
and (\ref{xi2}) hold with \(N=M\).

To determine a vector  \(\zeta_M \in K\) 
such  that, if we define \(b_1 \) by (\ref{b1definition}),
then (\ref{conditionb11}) and (\ref{conditionb12}) hold,
we proceed as follows.
 
First, let  \(F_1, F_2 \in {\cal F}\) be such that \(F_2 = S\setminus F_1\) and such that
\(a(F\cap F_1) \geq b(F\cap F_1)\)  for all \(F\in {\cal F}\) satisfying \(F\subset F_1\)
and such that 
\(
a(F\cap F_2) < b(F\cap F_2)\)  for all \(F\in {\cal F}\) satisfying \(F\subset F_2\).
(\(F_1, F_2\) constitutes a Hahn decomposition.)
We write \({\cal F}_1=\{F \in {\cal F}: F \subset F_1\}\) and 
\({\cal F}_2=\{F \in {\cal F}: F \subset F_2\}.\)

Next define a measure \(c \in {\cal Q}_{\lambda}(S, {\cal F})\) by
\begin{equation}\label{cdefinition}
c(F)= ((a-a_1)\wedge (a-b))(F\cap F_1),  \; F \in {\cal F},
\end{equation}
and set 
\[\Delta_0 = c(F_1).\]
Obviously \(\Delta_0 \leq \Delta\).
We now  define \(\zeta_M\) as follows:
\[
\zeta_M(F)= \xi_M(F) - c(F)/\beta_M, \; 
\;if\; F \in
{\cal F}_1 ,\]
\[
\zeta_M(F)= \xi_M(F) +(\Delta_0/\Delta)(b(F)-a(F))/\beta_M, \; 
\;if\; F \in
{\cal F}_2.\]

We have to verify  that \(\zeta_M \in K\).
We first show that \(\zeta_M \in {\cal Q}_{\lambda}(S, {\cal F})\).
For \(F \in {\cal F}_1 \) we find, from the definition of \(c\),
(see (\ref{cdefinition})), that
\[
\zeta_M(F)= \xi_M(F) - c(F)/\beta_M =
(a(F)-a_1(F)- c(F))/\beta_M \geq 0 
\]
and, if \(F \in {\cal F}_2\), then obviously 
\(\zeta_M(F) \geq 0\). 
Hence \(\zeta_M \in {\cal Q}(S, {\cal
  F})\). Since \(a, b, c\) and \(\xi_M\)
 belong to \({\cal Q}_{\lambda}(S, {\cal F})\),
it follows that also \(\zeta_M \in {\cal Q}_{\lambda}(S, {\cal
  F})\). 

To prove that \(\zeta_M \in K\), we need to show that 
\(\zeta_M(S)=1\).  Since
\[\zeta_M(F_1)=\xi_M(F_1)-\Delta_0/\beta_M\] 
and
\[
\zeta_M(F_2)=\xi_M(F_2)+(\Delta_0/\Delta)(b(F_2)-a(F_2))/\beta_M 
= \xi_M(F_2) + \Delta_0/\beta_M,\]
we find that 
\(\zeta_M(S)=\xi_M(F_1)+\xi_M(F_2)=
 1\),
and hence \(\zeta_M \in K\).
We also find that 
\[
||\xi_M-\zeta_M||=\xi_M(F_1)-\zeta_M(F_1) + \zeta_M(F_2)-
\xi_M(F_2) = \]
\begin{equation}\label{zetaxi}
c(F_1)/\beta_M+ c(F_1)/\beta_M = 2\Delta_0/\beta_M .
\end{equation}
Furthermore, 
if \(b_1\) is defined by (\ref{b1definition}), we find, that 
if \(F \in {\cal F}_1 \), then
\[b_1(F)=b(F)-\beta_M\xi_M(F)+c(F)=b(F)-a(F)+a_1(F)+c(F)=\]
\[
b(F)+a_1(F)+  ((a-a_1)\wedge (a-b))(F)  -a(F) = \]
\[ b(F)+a_1(F)-(a_1 \vee b)(F) \geq 0,\]
and, if  \(F \in {\cal F}_2 \), then, since \(\Delta_0 \leq \Delta\), we obtain
\[b_1(F)=b(F)-\beta_M\xi_M(F)-(b(F)-a(F))\Delta_0/\Delta \geq\]
\[
 b(F)-a(F)+a_1(F) -(b(F)-a(F)) \geq a_1(F).\]
Hence (\ref{conditionb11}) is satisfied.

It thus remains to show that (\ref{conditionb12}) is satisfied.
Since 
\[
b_1(F) = b(F)+a_1(F)-(a_1 \vee b)(F) \leq a_1(F),\]
if \(F \in {\cal F}_1\), and,
as we just showed, \(b_1(F) \geq a_1(F) \), if \(F \in {\cal F}_2,\)
we find 
\[
||a_1-b_1|| = a_1(F_1)-b_1(F_1)+ b_1(F_2)-a_1(F_2)=\]
\[
a(F_1)-\beta_M\xi_M(F_1) - b(F_1) + \beta_M\xi_M(F_1) +
c(F_1) + \]
\[b(F_2)-\beta_M\xi_M(F_2)-
(\Delta_0/\Delta)(b(F_2)-a(F_2)) - a(F_2) + \beta_M \xi_M(F_2) =\]
\[
a(F_1)-b(F_1) + \Delta_0 + b(F_2)-a(F_2) + \Delta_0 = 2\Delta + 2
\Delta_0
\]
and since \(||a-b||= 2\Delta \) and 
\(\beta_M||\xi_M - \zeta_M||= 2 \Delta_0 \) 
because of (\ref{zetaxi}), the equality (\ref{conditionb12}) holds
and thereby the proof of the lemma is completed. \(\;\Box\)

Using Lemma \ref{finitecase} and Lemma \ref{lowerbound} it is now 
easy to conclude the proof of 
Theorem \ref{barycentertheorem}. Thus 
let \(x, y \in K\) and suppose
\(\mu \in {\cal Q}^r(K|x)\). 
What we want to prove is that to every \(\epsilon > 0\)
we can find a measure \(\nu \in {\cal Q}^r(K|y)\) such that
\[
d_K(\mu,\nu) <r ||x-y||+\epsilon.
\]

Thus, let \(\epsilon > 0\) be given.
From the general theory of measures we know, since \((K,{\cal E})\) is 
a complete, separable, metric space, that we can find  a measure 
\(\mu_1 \in {\cal Q}^r(K,{\cal E})\)  of the form
\(
\mu_1= \sum_{k=1}^N\beta_k\delta_{\xi_k}
\)
such that \(d_K(\mu,\mu_1)<\epsilon/2\), where 
thus \(\xi_k, k=1,2,...,N\) belong to \(K\) and \(\beta_k > 0\) for 
\(k=1,2,...,N\). 
From Lemma \ref{lowerbound} now follows that we have
\[
\epsilon/2 > d_K(\mu,\mu_1) \geq ||rx- {\overline b}(\mu_1)||, 
\]
and from Lemma \ref{finitecase} follows that we can find a 
measure \(\nu\in {\cal Q}^r(K|y)\),
such that 
\[d_K(\mu_1,\nu)=||{\overline b}(\mu_1) - ry||.\]
From the triangle inequality then follows, that 
\[
d_K(\mu,\nu)\leq d_K(\mu,\mu_1)+d_K(\mu_1,\nu)<\epsilon/2 + 
||{\overline b}(\mu_1)-ry||\leq 
\]
\[\epsilon/2 + ||{\overline b}(\mu_1)-rx|| + r||x-y||
\leq \epsilon/2 + \epsilon/2 + r||x-y||.
\]
Hence,
\(
d_K(\mu,\nu)< r||x-y|| + \epsilon
\)
and thereby Theorem \ref{barycentertheorem} is proved. \(\;\Box\)

By using  Theorem \ref{kunita} and Theorem \ref{barycentertheorem}
we obtain the following corollary.
\begin{corr}\label{nearbarycentercorollary}
Let \({\cal H}=\{(S,{\cal F},\delta_0),(p,\lambda),(A,{\cal A},\varrho), 
(m,\tau)\}\)  
be a regular  HMM 
 and let \({\bf P}\) be the filter kernel.
\newline
A. Suppose  \({\cal H}\) is strongly ergodic
with limit measure \(\pi\).
Then, to every \(\eta > 0\) and every finite set \({\cal M}\) of elements in \(K\),
 we can find an integer \(N\), such that for every 
\(x  \in {\cal M}\) there exists
a probability \(\nu_x  \in {\cal P}(K|\pi)\), such that 
for every \(u \in Lip[K]\) 
\begin{equation}\label{nearbary1}
|\langle u,\delta_x {\bf P}^N\rangle -
\langle u, \nu_x\rangle| < \eta\gamma(u).
\end{equation}
B. If furthermore  \({\cal H}\) is uniformly ergodic, then, to every 
\(\eta > 0\), we can find an integer \(N\), such that, for every 
\(x \in K\), there exists a measure \(\nu_x \in {\cal P}(K|\pi)\), such that
the inequality (\ref{nearbary1})
 holds
for  every \(u \in Lip[K]\). 
\end{corr}
{\bf Proof}. Suppose \({\cal H}\) is strongly ergodic and
that 
 \({\cal M}=\{x_i, i=1,2,...,M\}\) is a finite set of elements in \(K\).
  From Theorem \ref{kunita} 
 follows, that  to every \(\eta > 0\),  we can find an integer \(N\), such that,  for every  \(x_i \in {\cal M}\),
\begin{equation}\label{barycenterinequality}
\delta_{TV}({\overline b}(\delta_{x_i}{\bf P}^n), \pi) < \eta, \;if \; n \geq N.
\end{equation}
From Theorem \ref{barycentertheorem} then follows that,
to every \(x_i \in {\cal M}\), we can find a probability measure 
\(\nu_i \in {\cal P}(K|\pi)\), such that 
\(d_K(\delta_{x_i} {\bf P}^N,\nu_i)<\eta\),
from which follows, that (\ref{nearbary1})  holds, 
if \(u \in Lip[K]\) and \(x \in {\cal M}\). Thereby part A is proved.

Next suppose  that \({\cal H}\) is uniformly ergodic with limit measure \(\pi\).
 From Theorem \ref{kunita} 
 follows that, to every 
\(\eta > 0\),  we can find an integer \(N\), such that,  for all \(x \in K \),
\begin{equation}\label{barycenterinequality2}
\delta_{TV}({\overline b}({\delta_x}{\bf P}^n), \pi) < \eta, \;if \; n \geq N.
\end{equation}
From Theorem \ref{barycentertheorem} then follows that, 
to every \(x \in K\), we can find a probability measure 
\(\nu_x \in {\cal P}(K|\pi)\), such that
\(d_K(\delta_{x} {\bf P}^N,\nu_x)<\eta\),
from which follows that
(\ref{nearbary1}) holds for all \(x \in K\),
if \(u \in Lip[K]\). Thereby Part B of the corollary is 
also proved. \(\Box\)

\section{Verifying the shrinking property}\label
{sectionverifyingshrinkingproperty}
From Lemma \ref{equicontinuitylemma1} we know that the filter kernel of 
 a regular HMM is Lipschitz equicontinuous. Therefore, by 
 Theorem \ref{auxiliarytheorem},  in order to prove 
 Theorem \ref{maintheorem}, it suffices
 to prove, that the filter kernel  of the  HMM under consideration
in Theorem \ref{maintheorem} 
 has the shrinking property, 
and, if also the HMM is uniformly ergodic, then the filter 
kernel  has the strong shrinking 
property.

We first prove the following lemma. 
\begin{lem}\label{shrinkingestimate}
Let \({\cal H}=\{(S,{\cal F},\delta_0),(p,\lambda),(A,{\cal A},\varrho), 
(m,\tau)\}\)  be a regular HMM which is strongly ergodic with limit measure 
\(\pi\). As usual, let \({\bf P}\) denote the filter kernel induced 
by \({\cal H}\). Suppose Condition E holds. Then,
\newline
A:  for every  \(\rho > 0\), there 
exists a number \(\alpha>0\) and an integer \(N\), such that for 
any two probability measures \(\mu\) and \(\nu\) in \({\cal P}(K|\pi)\)
\[
|\langle u, \mu{\bf P}^n\rangle-
\langle u, \nu{\bf P}^n\rangle| 
\leq \alpha \gamma(u)\rho + (1- \alpha)osc({\bf T}^{n-N}u),
\]
if  \( u \in Lip[K]\) and  \( n\geq N;\)
\newline
B:
 for every  \(\rho > 0\), there 
exists a number \(\alpha>0\) and an integer \(N\), such that for
any two probability measures \(\mu\) and \(\nu\) in \({\cal P}(K|\pi)\) 
and  any \(\kappa >0\), there exists a compact set \(F\)  such that
\[
|\langle u, \mu{\bf P}^n\rangle-
\langle u, \nu{\bf P}^n\rangle| \leq \alpha \gamma(u)\rho + \kappa osc(u) +
 (1- \alpha)osc_F({\bf T}^{n-N}u),\]
if  \( u \in Lip[K]\) and  \( n\geq N.\)
\end{lem}
{\bf Proof}.  Let \(\rho > 0\) be given and let 
\(\mu, \nu \in {\cal P}(K|\pi)\). Since Condition E is satisfied, there exist a number 
\(\alpha > 0\) and an integer \(N\) - independent of \(\mu \) and \(\nu\) - and
a coupling \({\tilde \mu}\) of \(\mu{\bf P}^N \) and \(\nu {\bf P}^N \) such that, 
if
\(D_{\rho}=\{(x,y)\in K^2: \delta_{TV}(x,y) 
<\rho/3\}\),
then  
\[{\tilde \mu}_N(D_{\rho}) \geq \alpha. \]
Hence, if \(u \in Lip[K]\),  \(n \geq N\) and we set \(v={\bf T}^{n-N}u\),  we find
\[
|\langle u,\mu {\bf P}^n\rangle -
\langle u,\nu {\bf P}^n\rangle | = 
|\int_{K}{\bf T}^{n}u(z)\mu(dz) - 
\int_{K}{\bf T}^{n}u(z)\nu(dz)| =\]
\[
|\int_{K}{\bf T}^{n-N}u(z)\mu{\bf P}^N(dz) - 
\int_{K}{\bf T}^{n-N}u(z)\nu{\bf P}^N(dz)| =\]
\begin{equation}\label{couplingequality}
|\int_{K\times K} (v(z)-v(z')) {\tilde \mu}(dz,dz')|.
\end{equation}

Next set 
\[
B_1= \{(z,z') \in K^2: \delta_{TV}(z,z')<\rho/3\},
\]
and
\[
B_2=\{(z,z')\in K^2 : \delta_{TV}(z,z_1)\geq \rho/3\}.
\]
Using the fact that \(\gamma({\bf T}^mu)\leq 3 \gamma(u), \forall m \geq 1 \) because of 
Lemma \ref{equicontinuitylemma1}, and that 
\begin{equation}\label{impinequality}
b \min \{\epsilon,\Theta\} + (1-b) \Theta \leq
a \epsilon + (1-a)\Theta,
\end{equation}
if 
\[ 0 <a \leq b \leq 1 \;, \;\epsilon > 0\; and 
\;\Theta > 0,
\]
 we obtain
\[
|\int_{K\times K} (v(z)-v(z')) {\tilde \mu}(dz,dz')| \leq
|\int_{B_1} (v(z) - v(z')){\tilde\mu}(dz,dz')| +\]
\[|\int_{B_2} (v(z) - v(z')){\tilde\mu}(dz,dz')|\leq
\min \{osc(v),\gamma(v)(\rho/3)\}{\tilde \mu}(B_1) + \]
\[osc(v) (1-{\tilde \mu}(B_1)) \leq 
 \gamma(v)(\rho/3) \alpha + (1-\alpha)osc(v) \leq \]
\[ \gamma(u)\rho \alpha +
(1-\alpha)osc({\bf T}^{n-N}u),\]
which combined with (\ref{couplingequality})  implies that
\[
|\langle u, \mu{\bf P}^n\rangle-
\langle u, \nu{\bf P}^n\rangle| \leq \alpha \gamma(u)\rho + (1- \alpha)osc({\bf T}^{n-N}u)
\]
and hence part A is proved.

Next let \(\kappa > 0\) also be given.  Since \((K,{\cal E})\) is a complete,
 separable, metric space, there exists a compact set \(F \in {\cal E}\) such that 
\begin{equation}\label{kappainequality}
{\tilde \mu}((K\setminus F)\times (K\setminus F)) \leq \kappa.
\end{equation}
Further,  define 
\[
B_3= \{
(z,z') \in K\times K: \delta_{TV}(z,z')<\rho/3, z\in F, z' \in F\},
\]
\[B_4=\{(z,z')\in K\times K : \delta_{TV}(z,z')\geq \rho/3, z\in F, z' \in F\}\]
and
\[
 B_5= K\times K\setminus (B_3\cup B_4). \]
Then,
\[
|\int_{K\times K} (v(z)-v(z')) {\tilde \mu}(dz,dz')|
\leq
|\int_{B_3} (v(z) - v(z')){\tilde\mu}(dz,dz')|+\]
\[
|\int_{B_4} (v(z) - v(z')){\tilde\mu}(dz,dz')|+
|\int_{B_5} (v(z) - v(z')){\tilde\mu}(dz,dz')| \leq
\]
\[\min \{osc_F(v),\gamma(v)(\rho/3)\}{\tilde \mu}(B_3) + 
osc_F(v) (1-{\tilde\mu}(B_3)) + osc(v){\tilde \mu}(B_5) \]
and by using 
(\ref{impinequality}),  (\ref{kappainequality}), the fact that \(\gamma(v) \leq 3 \gamma(u)\) because of
Lemma 
\ref{equicontinuitylemma1} and the fact that \(osc({\bf T}^nu) \leq osc(u)\) for all 
integers \(n\geq 1\),
 we find that 
\[
|\int_{K\times K} (v(z)-v(z')) {\tilde \mu}(dz,dz')| \leq
\alpha \gamma(u)\rho  + (1-\alpha) osc_F (v) + \kappa osc(u)
\]
which together with   (\ref{couplingequality}) and the fact 
that \(v = {\bf T}^{n-N}u\) implies that
\[
|\langle u, \mu{\bf P}^n\rangle-
\langle u, \nu{\bf P}^n\rangle| 
\leq
\alpha \gamma(u)\rho  + (1-\alpha) osc_F ({\bf T}^{n-N}u) + \kappa osc(u)
\]
and hence Part B is proved. \(\;\Box\)

To complete the proof of Theorem \ref{maintheorem} it suffices to 
prove the following two propositions.
\begin{prop}\label{strongshrinkingpropertyprop}
Let \({\cal H}=\{(S,{\cal F},\delta_0),(p,\lambda),(A,{\cal A},\varrho), 
(m,\tau)\}\)  be a regular HMM which is uniformly ergodic with limit measure 
\(\pi\). Suppose Condition E holds. Then the strong shrinking property
holds.
\end{prop}

\begin{prop}\label{shrinkingpropertyprop}
Let \({\cal H}=\{(S,{\cal F},\delta_0),(p,\lambda),(A,{\cal A},\varrho), 
(m,\tau)\}\)  be a regular HMM which is strongly ergodic with limit measure 
\(\pi\). Suppose Condition E holds. Then the shrinking property holds.
 \end{prop}
{\bf Proofs}. We first prove Proposition \ref{strongshrinkingpropertyprop}.
Let \(\rho > 0 \) be given. What we want to prove is
that we can find an integer \(N\) and a number \(\alpha > 0\)  
such that,
if the integer \(n\geq N\), 
then, for all \(u \in Lip[K]\),
\[
osc(T^{n}u)\leq
\alpha \rho \gamma(u) + 
(1-\alpha)
osc(T^{n-N}u).
\]

Set \(\rho_1= \rho /6\).  
From Part A of  Lemma \ref{shrinkingestimate} we know that we can find a number
 \(\alpha>0\) and an integer \(N_2\), such that, for 
any two probability measures \(\mu\) and \(\nu\) in \({\cal
  P}(K|\pi)\),
we have
\begin{equation}\label{N2estimate}
|\langle u, \mu{\bf P}^n\rangle-
\langle u, \nu{\bf P}^n\rangle| 
\leq \alpha \gamma(u)\rho_1 + (1- \alpha)osc({\bf T}^{n-N_2}u),
\end{equation}
if  \( u \in Lip[K]\) and  \( n\geq N_2\). 
Since we have assumed that \({\cal H}\)
 is uniformly ergodic with limit measure \(\pi\), it follows from
Part B of   Corollary 
\ref{nearbarycentercorollary}, 
that we can find an integer \(N_1\), such that 
for any two probability measures  \(x\) and \(y\) in \(K\) 
there exists probability measures
\( \nu_x\) and \(\nu_y\) in \({\cal P}(K|\pi)\),
such that for all \(u \in Lip[K]\)
\begin{equation}\label{N1estimatex}
|\langle u,\delta _x {\bf P}^{N_1}\rangle -
\langle u, \nu_x\rangle| < \alpha(\rho_1/2)\gamma(u)
\end{equation}
and
\begin{equation}\label{N1estimatey}
|\langle u,\delta _y {\bf P}^{N_1}\rangle -
\langle u, \nu_y\rangle| < \alpha(\rho_1/2)\gamma(u).
\end{equation}

Now, set \(N=N_1+ N_2 \), let \(n \geq N\), set \(m=n-N_1\)  and 
let  \(x\) and \(y\)
be two arbitrary probability measures in \(K\).  From 
(\ref{N1estimatex}), (\ref{N1estimatey}), (\ref{N2estimate}) and  
(\ref{gammainequality})
 follows that  
\[
|{\bf T}^nu(x)-{\bf T}^nu(y)| =
|\langle {\bf T}^mu, \delta_x {\bf P}^{N_1}\rangle -
\langle {\bf T}^mu, \delta_y {\bf P}^{N_1}\rangle|\leq 
\]
\[
|\langle {\bf T}^mu, \nu_x\rangle -
\langle {\bf T}^mu, \nu_y\rangle | + \gamma({\bf T}^mu)\alpha\rho_1 \leq\]
\[
|\langle u, \nu_x{\bf P}^m\rangle -
\langle u, \nu_y {\bf P}^m\rangle | + 3\gamma(u)\alpha\rho_1 \leq\]
\[
 \alpha \gamma({\bf T}^mu)\rho_1 + (1- \alpha)osc({\bf T}^{m-N_2}u) +
\gamma(u)\alpha\rho/2 \leq\]
\[
\alpha 3 \gamma(u)\rho/6 + (1-\alpha)
osc({\bf T}^{n-N}u) + \gamma(u)\alpha \rho/2
=\alpha \gamma(u) \rho + (1-\alpha) osc({\bf T}^{n-N}u).\] 
Hence
\[
osc({\bf T}^nu) \leq \alpha \gamma(u)\rho + (1- \alpha)osc({\bf T}^{n-N}u),
\]
and hence the strong shrinking property holds
and thereby Proposition \ref{strongshrinkingpropertyprop}  is proved.

We now prove Proposition \ref{shrinkingpropertyprop}.
Let \(\rho > 0 \) be given. What we want to prove  is, that we can find 
a number \(\alpha > 0\),  such that  for any nonempty, compact set \( E \in {\cal E}\) ,
any \(\eta > 0\) and any \(\kappa > 0\), we can find a nonempty compact set 
\(F\) and an integer \(N\), such that
\begin{equation}\label{finalestimate}
osc_E({\bf T}^{n}u)\leq
\eta \gamma(u) + \kappa osc(u)+ 
\alpha \rho \gamma(u) + 
(1-\alpha)
osc_F({\bf T}^{n-N}u)
\end{equation}
for all \(u \in Lip[K]\). 

Thus, let also \(E \in {\cal E}\) be a given, nonempty, compact set, and let 
also   \(\eta > 0\) and \(\kappa > 0\) be given. Set \(\eta_1=\eta/12 \).
Since \(E\) is a nonempty, compact set  in a metric space,  we can find a finite set 
\({\cal M}=\{x_i, i=1,2,...,M\}\) of elements in \(K\) such that 
\[
\sup_{x \in E} \min \{\delta_{TV}(x,x_i): x_i \in {\cal M}\} < \eta_1.
\]
Since \({\cal M}\) is a finite set, it follows from part A of Corollary 
\ref{nearbarycentercorollary}, that there exists an integer \(N_1\) such that
for every \(x_i\) in \({\cal M}\) there exists a measure 
\(\nu_i \in {\cal P}(K|\pi)\) such that 
\begin{equation}\label{nearbary}
|\langle u,\delta_{x_i} {\bf P}^{N_1}\rangle -
\langle u, \nu_i\rangle| < \eta_1 \gamma(u).
\end{equation}
Set \({\cal V}= \{\nu_1, \nu_2,...,\nu_M\}\).

  From Part B of  Lemma \ref{shrinkingestimate} we know, that we can choose 
 \(\alpha>0\) and the integer \(N_2\)  in such a way,  that if 
\(\nu_i\) and \(\nu_j\) belong to \({\cal V}\), then  there exists 
 a compact set \(F_{i,j}\in {\cal E}\), such that
\[
|\langle u, \nu_i{\bf P}^m\rangle-
\langle u, \nu_j{\bf P}^m\rangle| < \alpha \gamma(u)\rho + \kappa osc(u) +
 (1- \alpha)osc_{F_{i,j}}({\bf T}^{m-N_2}u),
\]
if  \( u \in Lip[K]\) and  \( m\geq N_2.\)

By defining \(F=\cup_{1\leq i < j \leq M} F_{i,j}\) it clearly follows that we also 
have 
\[
|\langle u, \nu_i{\bf P}^m\rangle-
\langle u, \nu_j{\bf P}^m\rangle|= \]
\begin{equation}\label{Finequality}
| \langle {\bf T}^mu,\nu_i\rangle - \langle {\bf T}^mu, \nu_j \rangle| 
 < \alpha \gamma(u)\rho + \kappa osc(u) +
 (1- \alpha)osc_{F}({\bf T}^{m-N_2}u)
\end{equation}
if  \( u \in Lip[K]\) , \( m\geq N_2\) and \(\nu_i, \nu_j \in {\cal V}\).

Now set \(N=N_1+ N_2 \), let \(n \geq N\), set \(m=n-N_1\)  and let  \(x\) and \(y\)
be two arbitrary probability measures in \(E\). Let \(x_i \in {\cal M}\)
satisfy \(\delta_{TV}(x,x_i) < \eta_1\) and let  \(x_j \in {\cal M}\) 
satisfy \(\delta_{TV}(y,x_j) < \eta_1\). 
From the triangle inequality then  follows that
\begin{equation}\label{Eestimate}
|{\bf T}^nu(x)-{\bf T}^nu(y)| \leq
|{\bf T}^nu(x_i)-{\bf T}^nu(x_j)| + 2\eta_1\gamma({\bf T}^nu).
\end{equation}
 From (\ref{nearbary}) and the triangle inequality follows also that
\begin{equation}\label{closeestimate}
|{\bf T}^nu(x_i)-{\bf T}^nu(x_j)| \leq
| \langle {\bf T}^mu,\nu_i\rangle - \langle{\bf  T}^mu, \nu_j \rangle| + 
2 \eta_1 \gamma({\bf T}^mu).
\end{equation}
By  combining 
(\ref{Eestimate}),
(\ref{closeestimate}) and
 (\ref{Finequality}) 
 we find
\[
|{\bf T}^nu(x)-{\bf T}^nu(y)| \leq
\]
\[
2\eta_1\gamma({\bf T}^nu) + 2 \eta_1 \gamma({\bf T}^mu)+ 
\alpha \gamma(u)\rho + \kappa osc(u) +
 (1- \alpha)osc_{F}({\bf T}^{m-N_2}u). 
\]
Since  \(x\) and \(y\) are arbitrarily chosen in the given set \(E\), 
and  \(\gamma({\bf T}^nu) \leq 3\gamma(u) \), for all \(n\geq 1\),
it follows that 
\[
osc_E({\bf T}^nu) \leq 12 \eta_1\gamma(u) +  \alpha \gamma(u)\rho_1 +
\kappa osc(u) +  (1-\alpha)osc_{F}({\bf T}^{m-N_2}u)  
\]
and, since \(\eta_1=\eta/12\)  and
 \(m-N_2 = n-N\), we find
that
\[
osc_E({\bf T}^nu) \leq  \eta\gamma(u) +  \alpha \gamma(u)\rho +
\kappa osc(u) +  (1-\alpha)osc_{F}({\bf T}^{n-N}u), 
\]
if \(u \in Lip[K] \), which was what we wanted to prove.
Thereby the  proof of  Theorem \ref{maintheorem} is completed. \(\Box\)
\newline
{\bf Remark}. Consider the following condition. {\em Condition \({\cal E}\):
There exists \(z \in K\), such that for every \(\rho > 0\) and every open set \(O\) containing 
\(z\), there exists an element \(x\in K\) such that 
\[
\limsup_{N\rightarrow \infty} \frac{1}{N} \sum_{n=1}^N {\bf P}^n(x,O) > 0.
\]
}
From the proof of Proposition 2.1 of \cite{Sza06} and Lemma 5.1, it follows that,
if a HMM is strongly ergodic and Condition \({\cal E}\) holds, 
then there exists an element \(z \in K\) such that \(\{P^n(z, \cdot), n=1,2,...\}\) is
 a tight sequence. Therefore, if we could verify Condition \({\cal E}\) 
then we could replace the 
conclusion {\em "weakly contracting"} by the conclusion {\em "weakly ergodic"}
  in the first part of Theorem \ref{maintheorem} and we could  omit the second part. \(\Box\)

\section{The random mapping associated to a HMM}
\label{sectionrandommapping}
The purpose of the remaining part of the paper is to introduce some further conditions, which in concrete applications probably will be easier to verify than Condition E. We will conclude the
paper with two simple - and rather concrete examples.

In this section we shall introduce a notion which we call {\em the random mapping associated 
to a regular HMM}. The motivation for this is twofold. One reason  is that we obtain
useful notations. 
The other reason is   that by introducing  random mappings
we build a bridge 
between the theory of filtering processes and 
the theory of random systems with complete connections.

Let 
 \({\cal H}=\{(S,{\cal F}, \delta_0),
(p,\lambda),   (A,{\cal A},\varrho), (m,\tau)\}\) be a
regular HMM.  As usual, let \(K={\cal P}_{\lambda}(S, {\cal F})\)
and  let 
\({\cal E}\) denote the Borel field on \(K\) 
induced by the total variation distance.
Furthermore, as defined in Section \ref{sectionregularhmm}, let  
 \(M_a:{\cal Q}_{\lambda}(S,{\cal F})\rightarrow 
{\cal Q}_{\lambda}(S,{\cal F}) \) be defined by 
\(
M_a(x)(F) =\int_S\int_Fm(s,t,a)x(ds)\lambda(ds)
\),
let  \(g:K\times A \rightarrow [0,\infty)\) be defined by
\(
g(x,a)=||xM_a||\), let 
 \(G:K\times {\cal A} \rightarrow [0,\infty )\) be defined by
\(G(x,B)=\int_B g(x,a)\tau(da) \;
\)
 and let \(h: K\times A \rightarrow K \) be defined  by
\(
h(x,a)= xM_a/||xM_a||
\) 
if \(||xM_a|| > 0\)
and  \(h(x,a)=x \) if \(||xM_a||=0\).
(See (\ref{gdefinition00}),  (\ref{Gdefinition00}), 
 (\ref{hdefinition00}) and(\ref{hdefinition000}).)
\begin{definition}\label{definitionrandom mapping0}
Let \({\cal H}=\{(S,{\cal F}, \delta_0),
(p,\lambda),   (A,{\cal A},\varrho), (m,\tau)\}\) be a
regular HMM and let  \(g:K\times A \rightarrow [0,\infty)\), 
\(\;G:K\times {\cal A} \rightarrow [0,1]\) and \(h:K\times A\rightarrow K\) 
be defined by (\ref{gdefinition00}),
 (\ref{Gdefinition00}), (\ref{hdefinition00}) and 
(\ref{hdefinition000}) respectively. 
We call the 4-tuple
\[
\{(K,{\cal E}), (A, {\cal A}), (g,\tau), h\}
\]
the {\bf random mapping associated to \({\cal H}\)} and we call 
\(G\) the tr.pr.f generated by  \((g,\tau)\). \(\;\Box\)
\end{definition}

Next, for \(x\in K\),  as in Section \ref{sectionregularhmm},
we set \(A_x^+=\{a \in A: ||xM_a|| > 0\}\)  and, if   
\(E \in {\cal E}\), we set \(B(x,E) = \{a \in A: h(x,a) \in E\} \).
 From the definition of the filter kernel  \({\bf P}\)
(see (\ref{filter})) we find that
\[ {\bf P}(x,E) =\int_{A^+_{x}} I_E(\frac{xM_a}{||xM_a||})||xM_a||\tau(da) =
\int_{A^+_{x}} I_E(h(x,a))||xM_a||\tau(da) =\]
\begin{equation}\label{PGequality0}
\int_{B(x,E)}g(x,a)\tau(da) = G(x, B(x,E)) 
\end{equation}
and, if \(u\in B[K]\), we find that
\begin{equation}\label{trm0}
\langle u, \mu{\bf P} \rangle = \langle {\bf T}u, \mu\rangle =
\int_{K}\int_A u(h(x,a))g(x,a)\tau(da)\mu(dx).
\end{equation}
{ \bf Historical remark}.
The random mapping associated to a regular HMM can be considered as a 
{\em random system with complete connections}. (See e.g 
\cite{IG90} for the definition of a random system with complete
connections.) As mentioned in the introduction other names for the concept 
random system with complete connections
are {\em learning model}  and 
{\em 
iterated function system with place-dependent probabilities}.
The terminology {\em random mapping} is inspired by the notion "random function" 
used  in the paper \cite{DF99} by P Diaconis and D Freedman. 

That there is  a random mapping - or a random system with complete connections - associated
to a regular HMM,  is not a new observation. 
Already in  1957,  Blackwell proves a theorem  
(\cite{Bla57}, Theorem 2) 
for random systems with complete connections,  
which he applies to the filtering process 
he is considering.
(Theorem 2 of \cite{Bla57} was in fact proved already 1937 
by W Doeblin and R Fortet in the classical paper \cite{DF37}.) 
In section 2.3.3.1 of the book \cite{IT69} from 1969 
the connection between partially observed Markov chains (HMMs) 
and random systems with complete connections  is described and 
also  in the book \cite{IG90}  
this connection is  mentioned at several places.
In the
paper \cite{Kai73} from 1973,  a  HMM with finite state space is considered and  it is proved that
the associated random mapping 
 is  a so called 
{\em distance diminishing model} as defined by F Norman in  
Chapter 2 of \cite{Nor72}, 
if
the tr.pr.m of the hidden Markov chain
is  {\em strictly positive}; from this fact it follows
 that   
 the filtering process converges in distribution with 
{\em geometric
convergence rate}. In the paper 
 \cite{Ant12} from 2012 by C Anton Popescu a similar result is proved.
 The connection between filtering processes and random systems with 
complete connections
is also utilized  in \cite{Kai75}.  \(\;\Box\)

Our next aim is to define the {\em Vasershtein coupling} of the random mapping.
 associated to a regular  HMM.

As before, let \(K^2= K\times K\), \({\cal E}^2={\cal  E}\times {\cal E}\), 
\(A^2=A\times A\) and \( {\cal A}^2 = {\cal A}\otimes {\cal A}\).
Let \(D=\{(a,b)\in A^2:a=b \}\). The set \(D\) is measurable, since 
\((A,{\cal A},\varrho)\) is a complete, separable, metric space.
For \(x,y \in K\), define \(C_1(x,y) = \{a: g(x,a) \geq g(y,a)\}\),
define \(C_2(x,y)=A\setminus C_1(x,y)\)
and define  \(C^2(x,y)= \{(a,b)\in A^2: a \in C_1(x,y), \;b \in C_2(x,y)\}\).
For \(B \in {\cal A}^2\),
we define \(\Pi (B) = \{a \in A: (a,a) \in B\}\). 
That \(A_1(x,y)\) and \(A_2(x,y)\) are measurable is obvious since the function 
\(g\) is continuous, and that \(\Pi (B) \in {\cal A}\) follows from the fact
that the set \(D\)  is measurable together with the fact that 
the mapping \(\vartheta:A \rightarrow A^2 \) defined by 
\(\vartheta(a)= (a,a) \) is measurable. 

Next define \({\check g}:K\times K\times A \rightarrow [0, \infty)\) by
\({\check g}(x,y,a) = \min\{g(x,a), g(y,a)\}\) 
and for \(x,y\in K\) define 
\(\Delta(x,y) = \int _A( g(x,a)-{\check g}(x,y,a))\tau(da)/2.\)
We define 
\({\tilde G}_V:K^2\times {\cal A}^2 \rightarrow [0,1]\)  by
\[
{\tilde G}_V((x,y), B) = \int_{\Pi(B)} {\check g}((x,y),a)\tau(da) + \]
\begin{equation}\label{GVdefinition}
\int \int_{B\cap C^2(x,y)}
(g(x,a)-{\check g}(x,y,a))(g(x,b)- {\check g}(x,y,b))\tau(da)\tau(db)/\Delta(x,y)
\end{equation}
where the last term is omitted if \(\Delta(x,y)=0\).

That \({\tilde G}_V\) is a tr.p.f from \((K^2, {\cal E}^2)\) to 
\((A^2, {\cal A}^2)\) is easily verified and that  
\({\tilde G}_V((x,y), \cdot)\) is a coupling of 
\(G(x, \cdot) \) and \( G(y, \cdot) \) for all \( x,y \in K\) 
where thus \(G:K\times {\cal E}\rightarrow [0, \infty) \) is the 
tr.pr.f generated by \((g,\tau)\), 
is easily checked - and well-known. (See \cite{Lin92}, Section I.5.)
We call \({\tilde G}_V:K^2\times {\cal A}^2 \rightarrow [0,1]\) the
 {\em Vasershtein coupling} of 
\((g,\tau)\) or of \(G\).
 
Next, 
define 
\({\tilde h}:(K\times A)\times (K\times A) \rightarrow K^2 \) by
\[
{\tilde h}((x,a),(y,b))= (h(x,a), h(y,b)).\]
Since \(h:K\times A \rightarrow K\) is measurable, so is 
\({\tilde h}:(K\times A)\times (K\times A) \rightarrow K^2 \).
We  call the 4-tuple 
\(\{(K^2, {\cal E}^2), (A^2, {\cal A}^2), {\tilde G}_{V}, {\tilde h}\}\)
the {\em Vasershtein coupling} of the random mapping 
\(\{(K,{\cal E}), (A, {\cal A}), (g,\tau), h\}\).
\newline
{\bf Remark}. The original paper using the  Vasershtein coupling 
is \cite{Vas69}. For an  early application of the Vasershtein coupling
to random systems with complete connections see \cite{Kai81}, where the 
Vasershtein coupling is used when proving 
the central limit theorem for the so called  state sequence of a
random system with complete connections. See also \cite{Kai94}, 
Sections 5-8 for other applications.
 \(\;  \Box \)

 Next, for \((x,y) \in K^2\) and 
 \({\tilde E} \in {\cal E}^2\), we set 
\[
{\tilde B}((x,y),{\tilde E}) = 
\{(a,b) \in A^2 : (h(x,a),h(y,b))\in {\tilde E}\},\]
and we  define 
\({\bf {\tilde P}}_V: K^2\times {\cal E}^2 \rightarrow [0,1] \)
 by 
\begin{equation}\label{Ptildedefinition}
{\bf {\tilde P}}_V((x,y), {\tilde E}) = 
{\tilde G}_V((x,y), {\tilde B}((x,y),{\tilde E})).
\end{equation}
Since 
\({\tilde G}_V:K^2\times {\cal A}^2 \rightarrow [0,1]\)  
is a tr.pr.f and 
\({\tilde h}:(K\times A)\times (K\times A) \rightarrow K^2 \) is measurable,
it follows from  Lemma 1.41 of \cite{Kal02} that 
 \({\bf {\tilde P}}_V\) is a tr.pr.f on \((K^2, {\cal E}^2)\).
That    
\({\bf {\tilde P}}_V((x,y), \cdot)\) is a coupling of 
\({\bf  P}(x, \cdot) \) and \({\bf  P}(y, \cdot) \)
for every \(x,y\) in \(K\),
follows easily from the fact that 
\({\tilde G}_V(x,y,\cdot)\) is a 
coupling of \(G(x, \cdot) \) and \(G(y,\cdot)\).
Therefore, if \(\mu, \nu \in {\cal P}(K,{\cal E})\)
and we define \({\tilde \mu} \in {\cal P}(K^2, {\cal E}^2)\) as the 
product measure of \(\mu\) and \(\nu\), it 
follows that
 \({\tilde \mu}{\bf {\tilde P}}_V \) is a coupling of 
\(\mu {\bf P }\) and \(\nu {\bf P}\). 

For sake of convenience we call \({\bf {\tilde P}}_V\) the 
{\em V-coupling} of  \({\bf P}\) induced  by the
coupling \({\tilde G}_V\) and we call
 \({\tilde \mu}{\bf {\tilde P}}_V \) the {\em V-coupling} of 
\(\mu {\bf P}\) and \(\nu {\bf P}\) induced by the coupling \({\tilde G}_V\). 

An important property of the Vasershtein coupling \({\tilde G}_V\)
  is described  in the next proposition. 

\begin{prop}\label{propositionGV0}
 Let 
\({\cal H}=\{(S,{\cal F}, \delta_0),(p,\lambda),   (A,{\cal A},\varrho), (m,\tau)\}\)
be a regular HMM, let 
\(\{(K,{\cal E}), (A, {\cal A}), (g,\tau), h\}\) be the associated
random mapping 
and let \({\tilde G}_V: K^2\times {\cal A}^2  \rightarrow
[0, \infty)\) be the Vasershtein coupling of \((g, \tau)\).  

Let \(K_0 \in {\cal E}\), let \(B \in {\cal A}\), and 
let \( 0 <  \beta, \eta < \infty\)  be such that \(\tau(B)=\beta >0 \)
and 
\(g(x,a) \geq \eta , \forall x \in K_0, \;\forall a \in B\).
Then 
\[{\tilde G}_V((x,y), \{(a,a): a \in B\})\geq \eta \beta ,\;\; \forall x,y \in K_0.\]
\end{prop}
{\bf Proof}.
Let \(x,y \in K_0\). From the definition (\ref{GVdefinition})
 of \( {\tilde G}_V \) it follows that 
\[
{\tilde G}_V((x,y),  \{(a,a): a \in B\}) = \int_B \min \{g(x,a), g(y,a)\}\tau(da) \geq \eta \beta  \]
which was to be proved. \(\;\Box\)
 
Next, let again 
\({\cal H}=\{(S,{\cal F}, \delta_0),(p,\lambda),   (A,{\cal A},\varrho), (m,\tau)\}\)
be a regular HMM. Set  
\(
A^1=A\;\),  \({\cal A}^1={\cal A}\)  and, for \(n=2,3,...\), define \(A^n\) and  \({\cal A}^n\) recursively by 
\( \; A^{n+1}=A^1\times A^n\)
and \({\cal A}^{n+1}={\cal A}^1 \otimes {\cal A}^n\).
For \((a_1,a_2,...,a_n) \in A^n\) we often write
\(
a^n=(a_1,a_2,...,a_n)
\)
and, if  \((a_1,a_2,...,a_n) \in {\cal A}^n\), we write 
\[
M_{a_1}M_{a_2}...M_{a_n}=M^n_{a^n},
\]
where thus \(M_{a}\) for \(a \in A\) is defined by (\ref{Mamap}).

It will be convenient to introduce the following mappings. For 
\(n=1,2,...\), we  define
\(h^{(n)}:K\times A^{n} \rightarrow K, n=1,2,...\) 
by
\begin{equation}\label{h(n)definition1}
h^{(n)}(x,a^n)=\frac{xM^n_{a^n}}{||xM^n_{a^n}||},
\; if\;  ||xM^n_{a^n}|| > 0
\end{equation}
and by 
\begin{equation}\label{h(n)definition2}
h^{(n)}(x,a^n)=x,
\; if\;  ||xM^n_{a^n}||=0,
\end{equation}
and we define
\(
g^{(n)}:K\times A^n \rightarrow [0, \infty), n=1,2,...,\)
by
\begin{equation}\label{g(n)definition}  
g^{(n)}(x,a^n) = ||xM^n_{a^n}||.
\end{equation} 
We denote  the {\em n-product measure} of 
\(\tau \in {\cal Q}^{\infty}(A,{\cal A})\) by \(\tau^n\),
 we write \(\tau^1= \tau\)
and, for \(n=1,2,...\), we define
the tr.pr.f \(G^{(n)}\) from \((K,{\cal E})\) to \((A^n, {\cal A}^n)\)
by
\[
G^{(n)}(x, B)=\int_{B}g^{(n)}(x,a^n)\tau^n(da^n). 
\]
Clearly \(\{(K,{\cal E}), (A^n, {\cal A}^n), (g^{(n)},\tau^n),
h^{(n)}\}\) 
is the random mapping 
associated to the \(nth\)  iteration \({\cal H}^n\) of \({\cal H}\). (See
Section \ref{sectionregularhmm} for the definition of an iterated HMM.)

Next, let us 
for \(x \in K\) and \(E \in {\cal E}\), define
\[ 
B^{n}(x,E)= \{a^n \in A^n: h^{(n)}(x,a^n) \in E\},\; n=2,3,...\,  .
\]
From  (\ref{PTrelations}) and  
(\ref{PGequality0}) follows  that, for \(n=2,3,...\), 
\begin{equation}\label{PGEequality}
{\bf P}^n(x,E)= {\bf P}^{(n)}(x,E) =G^{(n)}(x,B^{n}(x,E)),
\end{equation}
a representation of \({\bf P}^n\), which we will have use of below.

 We shall now introduce a slightly stronger condition than 
Condition E, a condition  which is 
formulated by using the functions \(g^{(n)}\) and \(h^{(n)}\) of  
the random mapping associated to the \({\cal H}^n\).

\begin{definition} Let 
 \({\cal H}=\{(S,{\cal F}, \delta_0),(p,\lambda),   (A,{\cal A},\varrho), (m,\tau)\}\) be a
strongly ergodic,  regular HMM with limit measure \(\pi\). 
 For \(n=1,2,... \) let
\({\cal H}^n \) denote the \(nth\) iteration of \({\cal H}\) and let 
\(\{(K,{\cal E}), (A^n, {\cal A}^n), (g^{(n)}, \tau^n), h^{(n)}\}\) denote
the random mapping associated to \({\cal H}^n\).
\newline
{\bf Condition E1}:
To every \(\rho > 0\), there exists an integer \(N\),  a set \(K_0 \in {\cal E}\), 
a set \(B \in {\cal A}^N\) and positive constants \(\xi\), \(\beta\) and   \(\eta\),
such that  
\newline
1) \[\mu(K_0) \geq \xi, \;\; \forall \mu \in {\cal P}(K|\pi),\]
\newline
2)
 \[
\tau^N(B) \geq \beta,\]
3)
 if \(x \in  K_0\)  and \( a^N \in B\), then
\begin{equation}\label{gNinequality}
g^{(N)}(x,a^N) \geq \eta ,
\end{equation}
4) 
if \(x, y \in K_0\)  and \(a^N \in B\), then 
\begin{equation}\label{hNinequality}
||h^{(N)}(x,a^N) - h^{(N)}(y,a^N)||< \rho. 
\end{equation}
 \end{definition}

\begin{lem}\label{E1impliesE} Let 
 \({\cal H}=\{(S,{\cal F}, \delta_0),(p,\lambda),   (A,{\cal A},\varrho), (m,\tau)\}\) be a
strongly ergodic,  regular HMM with limit measure \(\pi\). 
Then  Condition E1 implies Condition E.
\end{lem}
{\bf Proof}.  Let \(\rho > 0\) be given. Choose the integer \(N\), 
the set \(K_0 \in {\cal E}\),
the set \(B \in A^N\), the constants \(\xi > 0\), \(\beta > 0\) and \(\eta > 0\)
such that hypotheses 1), 2), 3)  and 4) of Condition E1 hold.
 
Let \(\mu\) and \(\nu\)
belong to \({\cal P}(K|\pi)\).
 Let 
\( {\cal G}^{(N)}=\{(K, {\cal E}, (A^N, {\cal A}^N), (g^{(N)},\tau^N), h^{(N)}\} \)
be the random mapping associated to the \(Nth\) iteration of \({\cal  H}\)
and let 
\({\tilde {\cal G}}^{(N)}= 
\{(K^2, {\cal E}^2), (A^{2N}, {\cal A}^{2N}), {\tilde G}^{(N)}_{V}, {\tilde h}^{(N)}\}\)
be the Vasershtein  coupling of the random mapping \({\cal G}^{(N)}\).
Set \({\tilde B}=\{(a^N,b^N)\in A^N\times A^N: a^N=b^N, \;a^N \in B\}\).

Since  \(g^{(N)}(x,a^N) \geq \eta \),  if \(x\in K_0\) 
and \(a^N \in B\) , and also  \(\tau^N(B)\geq \beta \),
it follows from  Proposition \ref{propositionGV0}  that 
\[
 {\tilde G}^{(N)}_{V}((x,y),{\tilde B}) \geq \eta \beta
\]
if  \(x,y \in K_0\).
Now let 
\[D_{\rho} = \{(z_1,z_2) \in K\times K: ||z_1- z_2|| < \rho \},\]
let 
\[
{\tilde A}^N(D_{\rho}) = 
\{(a^N,b^N) \in A^N\times A^N: (h^{(N)}(x,a^N),h^{(N)}(y,b^N))\in D_{\rho}\}\]
and   let
\({\bf {\tilde P}}^{(N)}_V\)
 be the V-coupling of 
\({\bf P}^{(N)}\) induced  by  \({\tilde G}^{(N)}_{V}\).
From the definition of the V-coupling of a filter kernel (see 
(\ref{Ptildedefinition})) and the fact that \({\tilde B} \subset {\tilde A}^N(D_{\rho}) \),
it follows that 
\[
{\bf {\tilde P}}^{(N)}_V((x,y), D_{\rho}) = 
{\tilde G}_V^{(N)}((x,y),  {\tilde A}^N(D_{\rho}))\geq 
{\tilde G}^{(N)}_V((x,y), {\tilde B})\geq \beta \eta.
\]
Hence, if we define \({\tilde \mu} = \mu \otimes \nu\)
and set \(\alpha=  \xi^2\beta\eta\), then
\[
{\tilde \mu}{\bf {\tilde P}}^{(N)}_V(D_{\rho}) \geq \xi^2 \beta\eta = \alpha, \]
since \({\tilde \mu}(K_0\times K_0) \geq \xi^2\) and therefore, since 
\(
{\tilde \mu}{\bf {\tilde P}}^{(N)}\) is a coupling of 
\(\mu {\bf P}^N \) and \(\nu {\bf P}^N \), it follows that  Condition E holds. \(\Box\)

{\bf Remark 1}. 
Suppose that \({\cal H}=\{(S,{\cal F}, \delta_0), (p, \lambda), 
(A, {\cal A},\varrho), (m,\tau)\}\)  is a strongly ergodic, regular
HMM with finite state space, finite observation space,
stationary measure \(\pi\) and such that the hidden Markov chain is irreducible.
Suppose also  that Condition KR  is satisfied. (See (\ref{kr}).)
Using Condition KR and the fact that the hidden Markov chain
is an aperiodic, irreducible Markov chain on a finite state space,
it is  not difficult to prove that to every \(\rho > 0\)
we can find an integer \(N\) and a sequence \(b_1,b_2,...,b_N\) of 
elements in \(A\) such that the product 
\[\Lambda_N=\prod_{n=1}^N M(b_n)\]
 of stepping matrices is such that there
exist an element \(i \in S\) and a number \(\eta_1 > 0\)
 such that 
\newline
1) the {\em (i,i)th} element of 
the matrix \(\Lambda_N\) satisfies 
\[
(\Lambda_N)_{i,i} = \eta_1,
\]
and
\newline
 2)  if \(x,y \in K\) are  such that \((x)_i\geq (\pi)_i/2 \) and 
\((y)_i \geq (\pi)_i/2\) then
\[ ||\frac{x\Lambda_N}{||x\Lambda_N||} - \frac{y\Lambda_N}{||y\Lambda_N||}|| < \rho.
\]
(Note that \((\pi)_j > 0\) for all \(j \in S\) since the hidden Markov chain 
is irreducible.)

Therefore, if 
 \(\{(K,{\cal E}), (A^N, {\cal A}^N), (g^{(N)},\tau^N), h^{(N)}\}\) denotes the  
random mapping associated to the \(Nth\) iterate of \({\cal H}\),
 and we define \(B \subset A^N\) by
\(B=\{(b_1,b_2,...,b_N)\}\)
then clearly \(\tau^N(B)=1\), since 
we assume that \(\tau\) is the counting measure, when the observation 
space is finite. If we define \(K_0 \subset K\) by
\[K_0=\{x\in K: (x)_{i} \geq(\pi)_i/2 \}
\]
and set \((\pi)_i/2= \xi\),
we find that 
\(\mu(K_0) \geq \xi\)
because of Lemma 
\ref{barycenterlemma}. 
Furthermore, if we set \(\eta =  \eta_1\xi \), we find that, if 
\(x \in K_0\) and  \(a^N \in B\), then 
 \[g^{(N)}(x,a^N)=||xM(a^N)|| \geq \xi\eta_1 =\eta \]
and, if  also \( y \in K_0\),  then also
\[
||h^{(N)}(x,a) - h^{(N)}(y,a^N)||< \rho. 
\]
Hence Condition E1 is satisfied. \(\;\Box\)

{\bf Remark 2}. 
Suppose \({\cal H}=\{(S,{\cal F}, \delta_0), (p, \lambda), 
(A, {\cal A},\varrho), (m,\tau)\}\)  is a strongly ergodic, regular,
HMM with denumerable state space, denumerable observation space,
stationary measure \(\pi\) and an  irreducible hidden Markov chain. 

Suppose also the following condition, introduced in \cite{Kai11}, holds.
\newline
{\bf Condition B}: {\em  For every \(\rho>0\) there exists an element \(i_0 \in S\) such that
 if \(C \subset K \) is a compact set such that 
\begin{equation}\label{Cinequality}
\mu(C \cap \{x:(x)_{i_0} \geq (\pi)_{i_0}/2\}) \geq (\pi)_{i_0}/3,\;\;
\forall \mu \in {\cal P}(K|\pi),
\end{equation}
then we can find an integer \(N\) and a sequence \(b_1, b_2, ...., b_N\) 
such that if we let \(M(b_n)\),  for \( n=1,2,...,N\) denote
 the stepping matrix associated
to \(b_n\) and define  
\[
\Lambda_N= M(b_1)M(b_2)...M(b_N)\]
then 
\[
|| \delta_{i_0}\Lambda_N|| > 0
\]
and if \( x \in C\cap \{x : (x)_{i_0} \geq (\pi)_{i_0}/2\}\) then also
\[
||\frac{ x\Lambda_N}{||x\Lambda_N||} -
\frac{ \delta_{i_0}\Lambda_N}{|| \delta_{i_0}\Lambda_N||} ||
< \rho.
\]
}

We shall now show that Condition E1 is satisfied. Thus let  \(\rho > 0\).   Set \(\rho_1=\rho/2\).
Let \(i_0 \in S \) and the compact set \(C \subset K\) be such that 
 (\ref{Cinequality}) holds. That such a set
exists  for any choice of \(i_0\) follows from Lemma 5.7 of \cite{Kai11}.
Let \(K_0 \in {\cal E}\) be defined by
\(K_0=C \cap \{x:(x)_{i_0} \geq (\pi)_{i_0}/2\}.\)
Note that \((\pi)_{i_0} > 0\)  since the hidden Markov chain 
is irreducible.
From Condition B
follows that we can find 
 an integer \(N\) and a sequence \(b_1, b_2, ...., b_N\) 
such that if we  define  
\(
\Lambda_N= M(b_1)M(b_2)...M(b_N),\)
then \(
|| \delta_{i_0}\Lambda_N|| > 0\) and if \(x \in K_0\) then 
\[
||\frac{ x\Lambda_N}{||x\Lambda_N||} -
\frac{ \delta_{i_0}\Lambda_N}{|| \delta_{i_0}\Lambda_N||} ||
< \rho_1.
\]

Now let 
 \(\{(K,{\cal E}), (A^N, {\cal A}^N), (g^{(N)},\tau^N), h^{(N)}\}\) denote the  
random mapping associated to the \(Nth\) iterate of \({\cal H}\),
and  define \(B \subset A^N\) by
\(B=\{(b_1,b_2,...,b_N)\}\),
Then clearly \(\tau^N(B)=1\), since 
we assume that \(\tau\) is the counting measure when the observation 
space is denumerable.  Moreover, if
we define  \(\xi = (\pi)_{i_0}/3\),
then  
\[\mu(K_0)\geq \xi,\;\; \forall \mu \in {\cal P}(K|\pi).\]  
Therefore, if we define
\[\eta= ||\delta_{i_0}\Lambda_N||(\pi)_{i_0}/2,\] 
we find that if \(a^N \in B\) and \(x \in K_0\), then
 \(g^{(N)}(x,a^N) = ||x\Lambda_N|| \geq \eta\)
and  
\[
||h^{(N)}(x,a^N) - h^{(N)}(y,a^N)||=
||\frac{ x\Lambda_N}{||x\Lambda_N||} -
\frac{ y\Lambda_N}{||y\Lambda_N||}||\leq\]
\[
||\frac{ x\Lambda_N}{||x\Lambda_N||} -
\frac{ \delta_{i_0}\Lambda_N}{|| \delta_{i_0}\Lambda_N||} || +
||\frac{ y\Lambda_N}{||y\Lambda_N||} -
\frac{ \delta_{i_0}\Lambda_N}{|| \delta_{i_0}\Lambda_N||} || \leq 2\rho_1= \rho.
\]
Hence all the hypotheses determining  Condition E1  are fulfilled and 
hence Condition E1 holds.
\(\;\Box\)

\section{Estimates of  iterations of integral kernels}\label{sectioncontracting0} 

In order to verify Condition E1, we want to find conditions regarding 
a HMM 
\newline
 \(\{(S, {\cal F}, \delta_0), (p, \lambda), (A, {\cal A},
\varrho), (m,\tau)\}\) 
such that for every \(\rho > 0\) there exists an  integer \(N\), 
and a 
subset \(K_0 \in {\cal E}\) such that 
\[
||\frac{xM^n_{a^n}}{||xM^n_{a^n}||} -
\frac{yM^n_{a^n}}{||yM^n_{a^n}||}|| < \rho
\]
if \(x,y\) belong to \(K_0\),
 where thus \(xM^n_{a_n}\) is the measure 
in \({\cal Q}_{\lambda}(S,{\cal F})\) defined by
\[
xM^n_{a^n}(F) = \int_S\int_F m^n(s,t,a^n)\lambda(dt)x(ds)
\]
and   
\(m^n:S\times S \times A^n\rightarrow [0, \infty)\) is defined  recursively 
by \(m^1(s,t,a) =m(s,t,a)\) and 
\[
m^{n+1}(s,t,a^{n+1}) = 
\int_{S}m^n(s,\sigma, a^{n}) m(\sigma,t,a_{n+1})\lambda(d\sigma).
\]

In this section we shall prove a theorem 
in which an estimate for a class of 
nonnegative kernels is stated.

Let as usual \((S,{\cal F},\delta)\) be a complete, separable, metric space and 
let \(\lambda \) be a positive, \(\sigma -finite\) measure on
\((S,{\cal F})\). We define the set \(D_{\lambda}[S]\) as the set of all 
non-negative, measurable functions defined on \(S\times S\). If \(k\in 
D_{\lambda}[S]\) is such that
\[
\sup \{ \int_S k(s,t)\lambda(dt):  s\in S\} < \infty,
\]
 we call \(k\) a {\em density kernel}. Recall that 
\({\cal Q}(S, {\cal F})\) denotes the set of nonnegative, finite,  
 measures on \((S, {\cal F})\).
\begin{definition}
Let \(k \in D_{\lambda}[S]\). We say that 
\(k\) has {\bf rectangular support}  if there exist \(F \in {\cal F}\)
and \(G \in {\cal F}\) such that
\(\lambda(F)> 0
\;and \;
\;\lambda(G) > 0, \)  
and such that if \((s,t) \in F\times G\) then 
\(
k(s,t)>0
\)
and if \((s,t) \not \in F \times G \)
then \(k(s,t)=0.\)
We call \(F \times G\) the rectangular support of \(k\). 
\end{definition}
{\bf Remark}. In case \(S\) is a finite set, then the notion {\em  rectangular support}
is equivalent to the notion {\em subrectangular matrix}  presented in the introduction. 
\(\Box\)

The following theorem is a generalisation of Lemma 6.2 of \cite{Kai75}.

\begin{thm}\label{kerneltheorem0}
Let \(k_m,\;m=1,2,...,n, \;\; n \geq 1\), be density kernels
belonging to \(D_{\lambda}[S]\) 
having rectangular supports \(F_m\times G_m, \;m=1,2,...,n ,\)
where thus \(\lambda(F_m)\lambda(G_m) >0, \;m=1,2,...,n\).
Let \(K_m:S\times {\cal F}\rightarrow [0,\infty)\)
be defined by
\[
K_m(s,E) = \int_E k_m(s,t)\lambda(dt),
\]
and, 
for \(m=1,2,...,n,\) define \(K^{m,n}:S\times 
{\cal F} \rightarrow [0, \infty) \) recursively by \(K^{n,n}=K_n\) and 
\begin{equation}\label{Kiteration}
 K^{m-1, n}(s, E)= \int_Sk_{m-1}(s, t)K^{m,n}(t,E)\lambda(dt),\; m=n, n-1,...,2.
\end{equation}
Set \(K^n = K^{1,n}\) 
and, for \(x\in {\cal P}(S, {\cal F})\), let 
\(xK^n \in {\cal Q}(S, {\cal F})\)
be defined by \(xK^n(E) = \int_S K^n(s,E)x(ds)\). 

Now, suppose that there exist
 numbers 
\(\kappa_m \geq 1\) such that
for \(1 \leq m \leq n\),
\begin{equation}\label{crossratiobound}
\sup\{ \frac{k_m(s_1,t_1)k_m(s_2,t_2)} 
{k_m(s_2,t_1)k_m(s_1,t_2)}: s_1, s_2 \in F_m , \;t_1,t_2 \in G_m\}
\leq \kappa_m^2 .
\end{equation}
Suppose also, that 
\begin{equation}\label{supportcondition}
K^n(s,S) >0
\end{equation} 
for all \(s \in F_1\).

Then,
if  \(x, y \in {\cal Q}(S,{\cal F})\)
are such that \(x(F_1)>0\)  and also \(y(F_1) > 0\),
and \(n\geq 1\), it follows that 
\begin{equation}\label{Knestimate} 
|| \frac{xK^n}{||xK^n||} - \frac{yK^n}{||yK^n||} || 
\leq 2 \prod_{m=1}^{n}\frac{(\kappa_m-1)}{(\kappa_m +1)}. 
\end{equation}

\end{thm}
{\bf Proof}.  We first state the following lemma.
\begin{lem}\label{contractionlemma0}
Let \(n\geq 1\), let \(k_m, K_m, K^{m,n}, \;m=1,2,...,n, \;\;\) and \(K^n\)
  be defined, - and have
the same properties -,  as in 
Theorem \ref{kerneltheorem0}.
Then 
\begin{equation}\label{kernelinequality3}
\sup\{|\frac{K^n(s_1,E)}{K^n(s_1,G_n)} -\frac{
    K^n(s_2,E)}{K^n(s_2,G_n)}|: s_1,s_2 \in F_1, \;E \in {\cal F}\}
\leq  \prod_{m=1}^{n}\frac{(\kappa_m -1 )}{(\kappa_m +1 )}. \;\Box
\end{equation}
\end{lem}
{\bf Proof of Lemma \ref{contractionlemma0}}. The lemma is a simple 
consequence of the following proposition, which is a special version of 
a result due to E Hopf from 1963. (See Theorem 1 in \cite{Hop63}.)
\begin{prop}\label{Hopfresult} 
Let  \((S,{\cal F},\delta)\) be a complete, separable, metric space and 
let \(\lambda \) be a positive, \(\sigma -finite\) measure on
\((S,{\cal F})\) and let  \(k\in 
D_{\lambda}[S]\) be 
density kernel with rectangular support \(F\times G\).
Suppose that
there exists a number \(\kappa \geq 1\) such that
\[\sup\{ \frac{k(s_1,t_1)k(s_2,t_2)} 
{k(s_2,t_1)k(s_1,t_2)}: s_1, s_2 \in F , \;t_1,t_2 \in G\}
\leq \kappa^2.
\]
Let \(u,v \in B[S, {\cal F}]\) be nonnegative functions   such that 
\(\sup\{ \frac{v(t)}{u(t)} : t \in G \} < \infty.\)
Define \(u_1  :S \rightarrow [0,\infty)\) and 
\(v_1: S\rightarrow [0,\infty)\) by
\(u_1(s)=\int_Sk(s,t)u(t)\lambda(dt)\) and 
\(v_1(s)=\int_Sk(s,t)v(t)\lambda(dt)\). 
Then 
\[
osc_F ( \frac{v_1}{u_1}) \leq \frac{\kappa-1}{\kappa+1} 
osc_G(\frac{v}{u}).  \;\Box
 \]
\end{prop}

By applying Proposition \ref{Hopfresult}
 we find that for every \(E \in {\cal F}\)
\[osc(\frac{K_n(\cdot, E)}{K_n(\cdot, G_n)}) 
\leq \frac{\kappa_n - 1}{\kappa_n + 1}\]
and then, using the integral representation (\ref{Kiteration})
and  Proposition \ref{Hopfresult}, 
the inequality  (\ref{kernelinequality3}) follows easily by induction. 
\(\;\Box.\)

To conclude the proof of Theorem \ref{kerneltheorem0} we argue as
follows. (The argument is inspired by an argument in \cite{FK60}.)

Let \(x,y\in {\cal Q}(S,{\cal F})\) be such that both \(x(F_1)>0\) and 
\(y(F_1) > 0\).
We write \(K^n = U\). What we want to prove is that, 
if \(n\geq 1\), then 
\[ 
|| \frac{xU}{||xU||} - \frac{yU}{||yU||} || 
\leq 2\prod_{m=1}^{n}\frac{(\kappa_m - 1)}{(\kappa_m +1 )},
\]

Let  \(E \in {\cal F}\).  Then
\(xU(E)/||xU||\) can be written
\[
xU(E)/||xU|| = \int_{F_1}\frac{U(s,E)}{xU(G_n)}x(ds) =
\int_{F_1}\frac{U(s,E)}{U(s,G_n)} \alpha(ds),
\]
where thus 
\[
\alpha(ds) = \frac{ U(s,G_n)}{xU(G_n)}x(ds).\]
Evidently \(\alpha \in {\cal P}(S,{\cal F})\).  

In a similar manner we can write 
\[
yU(E)/||yU|| = 
\int_{F_1}\frac{U(s,E)}{U(s,G_n)} \beta(ds),
\]
where thus \(\beta \in {\cal P}(S,{\cal F}) \) is defined by
\[
\beta(ds) = \frac{U(s,G_n)}{yU(G_n)}y(ds).
\]
Hence, by using the 
inequality
(\ref{oscinequality}), we find
\[
  |\frac{xU(E)}{||xU||} -\frac{yU(E)}{||yU||}| =
|\int_{F_1}\frac{U(s,E)}{U(s,G_n)}\alpha(ds) -
\int_{F_1}\frac{U(s,E)}{U(s,G_n)}\beta(ds)|
\leq \]
\begin{equation}\label{Uestimate}
 \sup\{ \frac{U(s_1,E)}{U(s_1,G_n)}-\frac{U(s_2,E)}{U(s_2,G_n)}: 
s_1,s_2 \in F_1 \} (1/2)||\alpha-\beta ||
 \end{equation}
and since \(||\alpha - \beta || \leq 2 \) and (\ref{Uestimate})
holds for all \(E\in {\cal F}\), it follows from 
Lemma \ref{contractionlemma0} that  (\ref{Knestimate})
holds. \(\Box\)

We shall next prove a theorem for HMMs   based on 
  Theorem \ref{kerneltheorem0}. First however we introduce yet  another
condition.
\begin{definition}\label{conditionP}
Let 
\({\cal H}\)= 
\(\{(S,{\cal F}, \delta_0\},(p, \lambda), (A,{\cal A},\varrho ),  (m,\tau)\}\)
 be a strongly ergodic, regular HMM, with stationary measure \(\pi\). 
If  
there exists a set \(F_0 \in {\cal F}\),
and a set
\(B_0 \in {\cal A}\),
such that 
\newline
1)
\[\pi(F_0) > 0, \] 
2)
\[\tau( B_0) > 0,\]
3) there exist positive numbers \(d_0\), \(D_0\) and \(\beta_0\), such that 
for every \(  a \in B_0\)  there exists a 
subset \(F_1(a) \in {\cal F}\), such that 
\newline
(a) 
\[F_1(a) \subset F_0\]
(b)
\[ \lambda(F_1(a))\geq \beta_0
\]
(c)
\[
d_0 \leq  m(s,t,a) \leq D_0, \;\forall (s,t) \in F_0 \times F_1(a)
\]
(d)
\[
m(s,t,a) = 0,\; \; \forall (s,t) \in F_0 \times( F_0 \setminus F_1(a)),
\]
then we say that \({\cal H}\) satisfies {\bf Condition P}.
\end{definition}
{\bf Remark}. The idea to formulate a condition like Condition P 
comes  
from the paper \(\cite{KR06}\) by Kochman and Reeds 
 and their proof of the fact that 
Condition A of the paper \cite{Kai75} implies their ``
rank 1 condition''. 
Condition P, as  introduced above, 
is a rather straight forward generalisation of a condition
introduced in \cite{Kai09}, section 9.

\begin{thm}\label{theoremconditionP}
Let 
\({\cal H}\)= 
\(\{(S,{\cal F}, \delta_0\},(p, \lambda), (A,{\cal A},\varrho ),  (m,\tau)\}\)
 be a strongly ergodic, regular HMM, with stationary measure \(\pi\). 
Suppose Condition P is satisfied.
Then Condition E1  is satisfied. \(\;\Box\)
\end{thm}
{\bf Proof}.  Let \(F_0\), \(B_0\),
 \(m:S\times S\times A \rightarrow [0, \infty)\),
\(d_0, D_0\) , \(\eta_0\)  and \(F_1(a), a \in B_0\) be chosen such that the 
hypotheses of Condition P are satisfied. 

Let \(\rho > 0\) be given.
What we want to prove is that there exist an integer \(N\),
a set \(K_0\), 
a number \(\xi > 0\),   a set \(B \in
{\cal A}^N\),
a number \(\beta >0\)  and  a number \(\eta > 0\),
such that 
\newline
(i): \[\mu(K) \geq \xi, \;\;\; \forall \mu \in {\cal P}(K|\pi),\]
(ii): \[\tau^N(B)\geq \beta,\]
(iii):
for all \( x \in K_0\) and all \(a^N \in B\), 
\[||xM^N_{a^N}|| \geq \eta,\] and 
\newline
(iv): for all \(x, y \in K_0\) and \(a^N \in B\)
then 
\begin{equation}\label{rhocontraction}
||\frac{xM^n_{a^n}}{||xM^n_{a^n}||} -
\frac{yM^n_{a^n}}{||yM^n_{a^n}||}|| < \rho.
\end{equation}

The choice of \(K_0\) is simple; we simply set
\(K_0= \{x\in K: x(F_0)\geq \pi(F_0)/2\}\), where thus  
\(F_0\) is the set determined by Condition P.  Since \(\pi(F_0) >0\)
it follows from Lemma \ref{barycenterlemma} that if we set 
 \(\xi = \pi(F_0)/2\), then
\(\mu(K_0) \geq \xi\) if \(\mu \in {\cal P}(K|\pi)\)
and hence
hypothesis 1) of Condition E1 is fulfilled. 

Next, set \(\kappa=D_0/d_0\) where thus \(d_0\)
and \(D_0\) are the constants occurring in hypothesis 3) of Condition P.
From the hypotheses of Condition P it  follows, that, if 
\(a \in B_0\) and we define 
\(m_a \in D_{\lambda}[S]\) by \(m_a(s,t)=m(s,t,a)I_{F_0}(s)\), then
\(m_a\) has the rectangular support 
\(F_0\times F_1(a)\) and \(m_a\) also  satisfies
\begin{equation}\label{mcrossratiobound}
\sup\{ \frac{m_a(s_1,t_1)m_a(s_2,t_2)} 
{m_a(s_2,t_1)m_a(s_1,t_2)}: s_1, s_2 \in F_0 , \;t_1,t_2 \in F_1(a)\}
\leq \kappa^2 .
\end{equation}

We  now simply define the integer \(N\)
by 
\begin{equation}\label{N1definition}
N = \min \{ n:2(\frac{\kappa - 1}{\kappa + 1})^n\}<\rho \},
\end{equation}
and we define the set  \(B\) in \({\cal A}^{N}\)  by
\(
B= B_1\times B_2\times ... \times B_N,
\)
where \(B_i=B_0, i=1,2,...,N\).

By defining \(\beta= \tau(B_0)^N\) we find that \(\tau^N(B) =\beta >
0\) and hence hypothesis 2) of Condition E1 is fulfilled.

Next, let \(x \in K_0\) and \(a^N \in B\). Then 
\(
||xM^N_{a^N}|| = \int_{S}\int_S m^N(s,t,a^N)x(ds)\lambda(dt).
\)
 From condition  3) of Condition P  follows that,
if \(s \in F_0\), then 
\[
\int_S m^N(s,t,a^N)\lambda(dt) \geq d_0^N\prod_{i=1}^N \lambda(F(a_i))
\geq d_0^N \beta_0^N.
\]
Therefore, if we define 
\[
\eta =(\pi(F)/2)  d_0^N\beta_0^N
\]
and use the fact that \(x(F)\geq\pi(F)/2\) if \(x \in K_0\), we find that
\[
||xM^N_{a^N}|| \geq \int_F\int_S m^N(s,t,a^N)x(ds)\lambda(dt) \geq \eta.\]
Hence hypothesis 3) of Condition E1 is fulfilled.

It remains to show, that, if \(x,y \in K_0\) 
and \((a_1,a_2,...,a_N)=a^N \in B\),
then (\ref{rhocontraction}) holds. 
But this follows immediately from Theorem \ref{kerneltheorem0} 
and the definition of the integer \(N\). Hence also hypothesis 4)  
of Condition E1 is fulfilled and hence Condition E1 is satisfied
\(\Box\)

\section{Examples}\label{sectionexamples}
Our first example is obtained by making a denumerable partition of the state space.  
\begin{example}\label{examplepartition}
 Let \({\cal H}_1=\{(S,{\cal F},\delta_0),(p,\lambda), (A,{\cal
  A},\varrho),(m,\tau)\}\) be a regular HMM, such that \(A\) is a
denumerable set and 
such that 1) 
for each \(a \in A \) there exists a set \(S_a \in {\cal F}\) such
that \(\lambda(S_a)>0\), 2) 
\(\cup_a S_a= S\) and 3) for each \(a \in A\) 
\[
m(s,t,a) = p(s,t)I_{S_a}, 
\]
where as usual \(I_F\) denotes the indicator function of a set \(F \subset S\).
\end{example}

\begin{thm}\label{theorempartitionexample} 
Let \({\cal H}_1 =
\{(S,{\cal F},\delta_0),(p,\lambda), (A,{\cal A},\varrho),(m,\tau)\}\) 
be the HMM defined  in  Example \ref{examplepartition} and let \({\bf P}\)
denote the induced filter kernel.
Suppose that 
\newline
a) the hidden Markov chain determined by the tr.pr.f \(P\) is 
 strongly ergodic with stationary measure \(\pi\); 
\newline
b) there exist an element \(a_0 \in A\) and two positive numbers
\(d_0,
D_0\)  satisfying \(d_0 \leq D_0\), such that \(\pi(S_{a_0}) > 0\)
and
\newline
\[ d_0 \leq p(s,t) \leq D_0, \;\forall (s,t) \in S_{a_0}\times S_{a_0}. \]
Then the filter kernel \({\bf P}\) is weakly ergodic.
  \(\;\Box\)
\end{thm}
{\bf Proof}. We shall first verify that
the hypotheses of Condition P are fulfilled.

  First, 
let \(F_0=S_{a_0}\). By assumption \(\pi(S_{a_0})>0\) 
and therefore it obviously follows that \(\pi(F_0)>0\). Hence hypothesis 1) of 
Condition P is satisfied with this choice of \(F_0\).

Next set \(B=\{a_0\}\). Since \(\tau\)  
is the counting measure \(\tau(B)=1 >0\); 
hence hypothesis 2) of Condition P holds.

Now let \(F_1(a_0)=F_0\). Evidently \(F_1(a_0) \subset F_0\). Since
\(\pi(F_0)>0\) and 
\[\pi(F_0)=\int_{F_0} p(s,t)\pi(ds)\lambda(dt) \leq D_0\lambda(F_0)\pi(F_0),\]
it follows that \(\lambda(F_0)>0\). 
Hence conditions 3a) and 3b)  of Condition P are  satisfied.  

Further,   since \(m(s,t,a_0)=p(s,t) \) if 
\((s,t) \in F_0\times F_0\) and \(m(s,t,a_0) = 0 \) if \((s,t) \in F_0 \times 
(S\setminus F)\), it is clear that conditions 3c) and 3d) of Condition P hold.
Hence Condition P is satisfied.
From Theorem \ref{theoremconditionP},
Lemma \ref{E1impliesE}   and Theorem \ref{maintheorem} follows that 
the filter kernel is weakly contracting. 
If furthermore the Markov chain is uniformly ergodic then hypothesis 3) of Theorem \ref{maintheorem}
is fulfilled and 
the filter kernel is weakly ergodic. 

In order to prove that the filter kernel is weakly ergodic without this extra assumption, we shall
use a  result in \cite{Sza06}. We shall   show that 
the following condition is satisfied.
\newline
{\bf  Condition} \({\cal E}\)1: {\em There exists an element \(x_0 \in K\) such that 
for every \(\epsilon > 0\) 
\[\liminf_{n \rightarrow \infty} {\bf P}^n(x, B(x_0,\epsilon)) > 0, \;\forall x \in K \]
 where \(B(x_0,\epsilon) = \{y\in K:\delta_{TV}(x_0,y)< \epsilon\}.\)}

Once we have verified Condition \({\cal E}\)1,   if follows from 
 Proposition 2.1  of \cite{Sza06} and Lemma \ref{equicontinuitylemma1},
 that \(\{{\bf P}^n(z, \cdot), n=1,2,...\}\) is a tight 
sequence, since obviously  Condition \({\cal E}\)1  implies  Condition \({\cal E}\) of \cite{Sza06}.
(Condition \({\cal E}\) is also formulated at the end of 
Section \ref{sectionverifyingshrinkingproperty}.)
That the filter kernel is weakly  ergodic  follows then from  
Theorem \ref{maintheorem}, since hypothesis 2) of
Theorem \ref{maintheorem}  is fulfilled.

To verify  Condition \({\cal E}\)1 we argue as follows. Set \(F_0= S_{a_0}\), 
define \(k:F_0\times F_0 \rightarrow [0,\infty)\) by \(k(s,t)=p(s,t)\) and define
\(\kappa=D_0/d_0\).
Since \(d_0 \leq p(s,t) \leq D_0\) if \((s,t) \in F_0\times F_0\), it follows that there
exists a positive function \(q:F_0 \rightarrow (0,\infty)\)
satisfying \(\int_{F_0}q(t)\lambda(dt)=1\)  and a number \(\beta>0\)
such that \(\int_{F_0} K(s,t) q(t)\lambda(dt) = \beta q(s)\). (See e.g \cite{Hop63}.)
Moreover, if we define \(x_0 \in K\) by
\[x_0(F)= \int_F q(t)\lambda(dt)\]
it follows from Theorem \ref{kerneltheorem0} that for any \(x \in K\) such that \(x(F_0)>0 \)  
\[ 
|| \frac{xK^n}{||xK^n||} - x_0 || 
\leq 2( \frac{\kappa-1}{\kappa +1})^{n}. 
\]

Now let \(\epsilon > 0\) be given , and let \(x\in K\) be chosen arbitrary.
Define \(N_0\) by 
\[N_0 = \min \{n:  2( \frac{\kappa-1}{\kappa +1})^{n-1} 
 < \epsilon\}\]
and define 
\[\alpha  = d_0^{N_0}\lambda(F_0)^{N_0}.\]
From the definition of \(N_0\) follows that if \(\mu \in {\cal P}(K, {\cal E})\)
satisfies \[\mu(\{x: x(F_0)>\pi(F_0)/3\}) \geq \pi(F_0)/3) \]
then 
\begin{equation}\label{firstestimate}
\mu{\bf P}^{N_0}(\{z: \delta_{TV}(z,x_0)< \epsilon \}) \geq \alpha\pi(F_0)/3.
\end{equation}

Next, by Lemma \ref{barycenterlemma} it follows easily, that 
 if \(\mu \in {\cal P}(K, {\cal E})\)
satisfies \(\delta_{TV}({\overline b}(\mu),\pi)<\pi(F_0)/6\), then
\(\mu(\{z: z(F_0) \geq \pi(F_0)/3\})\geq  \pi(F_0)/3\).
and from Theorem \ref{kunita}  follows that we can  choose \(N_1\) 
so large that if \( n\geq N_1\) then 
\begin{equation}\label{secondestimate}
\delta_{TV}({\overline b}(x{\bf P}^n, \pi) < \pi(F_0)/6.
\end{equation}
Finally, by combining (\ref{firstestimate}) and  (\ref{secondestimate})
we conclude that if \(n \geq N_1+N_0\) then 
\({\bf P}^n(x, B(x_0,\epsilon)) \geq (\pi(F_0)/3)^2\alpha \)
and  hence Condition \({\cal E}\)1 is satisfied.
\( \Box\)

\begin{example}\label{secondexample} Let 
\(
{\cal H}_2= \{(S,{\cal F}, \delta_0),(p, \lambda), (A, {\cal A},
\varrho), (m, \tau)\}\) be a HMM with densities such that the 
probability density kernel 
\(m:S\times S \times A \rightarrow [0, \infty) \) can be written
\[
m(s,t,a)=p(s,t)q(t,a),\]
where thus \(q:S\times A \rightarrow [0,\infty)\)  is a measurable
function satisfying
\[
\int_{A}q(t,a)\lambda(dt)= 1 , \;\forall t \in S.
\]
We assume that \(\sup \{p(s,t):s,t \} < \infty\) and \( \sup \{q(s,t):s,t \} < \infty\).

For each \(a \in A\) set \(S_+(a)= \{t:q(t,a) > 0\}\). We assume that 
\(
\lambda(S_+(a)) > 0, \;\forall a \in A.
\)
We also assume that the density function \(q\) is such that 
 for every \(\epsilon > 0\), we can find an \(\eta > 0\), such that,
if \(\varrho(a,b) < \eta \), then 
\begin{equation}\label{epsilon1example2}
\lambda(S_+(a)\, \Delta \, S_+(b)) < \epsilon,
\end{equation}
where 
\[S_+(a)\,\Delta\, S_+(b) = (S_+(a)\setminus S_+(b))
\cup (S_+(b) \setminus S_+(a)),
\] 
and  
\begin{equation}\label{epsilon2example2}
|q(t,a)-q(t,b)|<\epsilon, \;\forall t \in S_+(a)\cap S_+(b).  \;\Box
\end{equation}
\end{example}

\begin{prop}\label{propexamplefuzzy}
Let \({\cal H}_2\) be defined as in Example \ref{secondexample}.
Then \({\cal H}_2\) is regular.
\end{prop}

{\bf Proof}.  
What we need to prove is that  
\({\overline M} : {\cal Q}_{\lambda}(S,{\cal F})\times A \rightarrow  
{\cal Q}_{\lambda}(S,{\cal F})\) is a {\em continuous} function where 
thus \({\overline M}\) is defined by
\[
{\overline M}(x,a)(F) = \int_S \int_F p(s,t)q(t,a)\lambda(dt)x(ds).
\]

That 
\({\overline M}:{\cal Q}_{\lambda}(S,{\cal F})\times A \rightarrow  
{\cal Q}_{\lambda}(S,{\cal F})\) 
is a continuous function in the first variable follows easily from the 
boundedness condition regarding the probability density kernel \(q\).

That 
\({\overline M}:{\cal Q}_{\lambda}(S,{\cal F})\times A \rightarrow  
{\cal Q}_{\lambda}(S,{\cal F})\) 
also is continuous in the second variable follows easily from
(\ref{epsilon1example2}) and (\ref{epsilon2example2}) together 
with the hypothesis that both \(q:S\times A\rightarrow [0,\infty)\)
and \(p:S \times S \rightarrow [0,\infty)\) are uniformly
bounded. Since the proof is elementary we omit the details.
\(\;\Box\)

\begin{thm}\label{theoremexample2a}
 Let \({\cal H}_2 =\{(S,{\cal F},\delta_0),(p,\lambda), (A,{\cal
  A},\varrho),(m,\tau)\}\) be the HMM defined  in  
Example \ref{secondexample} and let \({\bf P}\) denote 
the induced filter kernel.
Suppose that 
\newline
a) the HMM \({\cal H}_2\) is 
{\bf strongly ergodic} with stationary measure \(\pi\); 
\newline
b) there exists  a set \(F_0 \in {\cal F}\) and a set  \(B_0 \in {\cal A}\) 
such that
\newline
(i) 
\[\pi(F_0) > 0,\] 
(ii) 
\[
\tau(B_0) > 0.
\]
(iii):
\[ S_+(a) \subset F_0, \;\forall a \in B_0\]
(iv):  there exists a constant \(\beta_0\) such that 
\[\lambda (S_+(a))  \geq \beta_0 , \;\forall a \in B_0\]
(v):  there exists a constant \(c_0\) 
\[
\{t: 0 < q(t,a) < c_0\} = \emptyset , \;\forall a \in B_0
\]
(vi):  there exists a constant \(c_1 >0\) such that 
\[
p(s,t)\geq c_1, \; (s,t) \in F_0\times F_0.
\]
Then the filter kernel \({\bf P}\) is weakly contracting.
If furthermore  \({\cal H}_2\) is uniformly ergodic,
then  the filter kernel
 \({\bf P}\) 
 is weakly ergodic.  \(\;\Box\)
\end{thm}

{\bf Proof}.  
It suffices to verify that the HMM \({\cal H}_2\) satisfies 
the conditions 1)-3) of Condition P. 
We shall verify the hypotheses
of Condition P, when \(F_0, B_0,\beta_0\) are chosen as in 
the hypotheses of Theorem \ref{theoremexample2a} and 
\(F_1(a)=S_+(a)\).

Since the sets \(F_0\) and \(B_0\) are such that \(\pi(F_0)>0\) and
\(\tau(B_0)>0\),  conditions 1) and 2) of Condition P are satisfied.
Since \(\lambda(F_1(a))\geq \beta_0\) for all \(a \in B_0\)
because of hypothesis \((iv)\), 
it follows that  condition 3b) of Condition P is satisfied.
From hypothesis  \((iii)\) we know that \(F_1(a)\subset F_0\)
and  from hypothesis \((v)\) follows, that, if \(t \in F_1(a)\)
and \(a \in B_0\),
 then \(q(t,a)\geq  c_0\) and from hypothesis \((vi)\) we know that 
\(p(s,t)\geq c_1,\;\)if \(s,t \in F_0\).

From the assumptions we have made in Example \ref{secondexample}
 regarding the HMM \({\cal H}_2\), we know that there exist constants
\(C_0\) and \(C_1\) such that \(\sup_{t,a}q(t,a) \leq C_0\) and
\(\sup_{s,t}p(s,t) \leq C_1\). 
 Hence, if we define 
\(d_0=c_0c_1\) and \(D_0=C_0C_1\)  and recall that 
\(m(s,t,a)=p(s,t)q(t,a), \) 
we find that, if \(a \in B_0\) then
\(
d_0 \leq m(s,t,a) \leq D_0\) if \( (s,t) \in F_0\times F_1(a),  
\)
and that 
\(
m(s,t,a)=0 \), if \(s \in F_0\) and \(t \not \in F_1(a).
\)
Hence, also condition 3) of Condition P is satisfied,
and hence Condition P is satisfied.

The conclusions of the theorem now follows from Theorem \ref{theoremconditionP},
Lemma \ref{E1impliesE} and Theorem \ref{maintheorem}.
\( \Box\)

\section{Acknowledgements} I want to thank Sten Kaijser,
Fuzhou Gong 
and Luc{\'a}{\v s} Mal{\'y}  
for valuable discussions and Svante Janson for an important reference.

\end{document}